\documentclass[a4paper, 12pt]{amsart}
\usepackage{xcolor}
\usepackage{amsmath}
\usepackage{amsthm}
\usepackage{amssymb}
\usepackage{comment}
\usepackage{graphicx}
\usepackage{subcaption}
\usepackage{float}
\usepackage{MnSymbol}
\usepackage{enumerate}
\usepackage{tikz}
\usetikzlibrary{intersections, calc}
\usetikzlibrary{arrows,shapes,positioning}
\usetikzlibrary{decorations.markings}
\tikzstyle arrowstyle=[scale=1]
\tikzstyle directed=[postaction={decorate,decoration={markings,
    mark=at position 0.65 with {\arrow[arrowstyle]{latex}}}}]
\tikzstyle reverse directed=[postaction={decorate,decoration={markings,
    mark=at position 0.5 with {\arrowreversed[arrowstyle]{latex};}}}]
\usetikzlibrary{decorations.pathreplacing}
\theoremstyle{definition}

\newtheorem{proposition}{Proposition}[section]
\newtheorem{definition}{Definition}[section]

\newtheorem{lemma}{Lemma}[section]
\newtheorem{cor}{Corollary}[section]

\newcommand{\odd}{\textrm{odd}}
\newcommand{\even}{\textrm{even}}
\newcommand{\fpl}{\mathfrak{fpl}}

\def\Sym#1{\mathfrak{S}_{#1}}
\def\ASM#1{\mathcal{A}_{#1}}
\title{
	The characterizations of an Alternating Sign Matrix using a triplet
}
\date{
		\today
}
\author{
	Toyokazu Ohmoto 
}
\begin{document}
\maketitle
\begin{center}
	Division of Mathematics and Physics, \\
	Graduate School of Natural Science and Technology, \\
	Okayama University	\\
	Tsushima-naka, Okayama, 700-8530, Japan	\\
	\textit{e-mail address}: ohmoto@s.okayama-u.ac.jp
\end{center}
\begin{abstract}
	Alternating Sign Matrix(ASM for short) is a square matrix
	which is consist of $0$, $1$ and $-1$.
	In this paper,
	we characterize an ASM
	by showing a bijection between alternating sign matrix and six vertex model,
	and a bijection between six vertex model and height function.
	In order to show these bijections,
	we define a triplet $({a}_{i, j}, {c}_{i, j}, {r}_{i, j})$ for each entry of an ASM.
	We also define a \textsl{track} for each index of height function, and 
	state more properties of height function. 
\end{abstract}
\section{Introduction}\label{section:Introduction}
In this paper, we mainly deal with the Alternating sign matricies
which were adovocated by Mills-Robbins-Rumsey
\cite{MRR1}.
The Alternating sign matrix is a square matrix
consisting of $0$, $1$ and $-1$
which satisfies several conditions,
and is often abbreviated as ASM.
Among the proofs of their enumeration,
the proof by Zeilberger \cite{DZ}
and the proof by Kuperberg \cite{Ku}
are famaous. 
In the proof by Kuperberg \cite{Ku},
the six vertex models derived from statistical physics is used.
The six vertex model is described as a map
which gives direction to the edges of the lattice graph.
On the other hands,
there is a related object called Fully pucked loop model,
which is described as a map
such that it associates one of two colors with each edges of a lattice graph.
In order to get some characterizations of ASMs which has degree $n$,
a graph consisting of $n \times n$ vertices called interiror vertex
and $4n$ vertices called boundary vertices,
where each interior vertecies is adjacent to $4$ edges
is used.
Wieland~\cite{Wieland_2000} states
a characterization using a operation which called gyration
defined for FPLs.
As another approach,
a $(n + 1) \times (n + 1)$ matrix which called height function
is defined for a ASM of size $n$,
and it define partially order for ASMs of size $n$.
In this paper,
we focus on the conditions of ASM 
for the sum of a row
(resp. column),
and define a triplet $\left( {a}_{i, j}, {c}_{i, j}, {r}_{i, j} \right)$
for each entries of a ASM.
Then
we describe a correspondence ASMs, six vertex models and height functions
by using the triplets $\left( {a}_{i, j}, {c}_{i, j}, {r}_{i, j} \right)$.
We will also introduce the gyration.
\section{Alternating sign matrix}
In this section we define alternating sign matrices and related objects
called six vertex model and fully packed loop model.
There are bijections between these objects, and we explain those bijections
in this section.
\subsection{Alternating sign matrix (ASM)}
As it is well-known,
each element $\sigma$ in $\Sym{n}$
corresponds to its permutation matrix $P_\sigma=(\delta_{i, \sigma(j)})_{1\leq i,j\leq n}$.
For example, the permutation
$
	\begin{pmatrix}
		1	&	2	&	3	&	4	&	5	\\
		2	&	4	&	1	&	5	&	3
	\end{pmatrix}
\in\Sym{5}
$
corresponds to the following square matrix:
\[
	\begin{pmatrix}
		0	&	0	&	1	&	0	&	0	\\
		1	&	0	&	0	&	0	&	0	\\
		0	&	0	&	0	&	0	&	1	\\
		0	&	1	&	0	&	0	&	0	\\
		0	&	0	&	0	&	1	&	0
	\end{pmatrix}.
\] 
A permutation matrix
$A = {\left( {a}_{i, j} \right)}_{1 \leq i, j \leq n}$
 of size $n$
is characterized by the properties (i) $a_{i, j} \in \{ 0, 1 \}$ 
for $1 \leq i, j \leq n$ and 
$\sum_{i = 1}^{n} {a}_{i, j}=\sum_{j = 1}^{n} {a}_{i, j}=1$ for each $i,j$.
Now we define the alternating sign matrices
whose entries consist of $0, 1$ or $- 1$ 
and regarded as an extension of permutation matrices. 
\begin{definition}\label{def:ASM}
	A square matrix $A = {\left( {a}_{i, j} \right)}_{1 \leq i, j \leq n}$
	of size $n$ is called 
	{\sl alternating sign matrix} (or {\sl ASM} shortly)
	if it satisfies the following conditions:
	\begin{subequations}
		\begin{align}
			&a_{i, j} \in \{ 0, 1, -1 \} & (1 \leq i, j \leq n),  \\
			&\sum\limits_{j = 1}^{k} a_{i, j},\quad\sum\limits_{j = 1}^{k} a_{i, j} 
				\in \{ 0, 1 \} & (1 \leq i,j,k \leq n), \\
			&\sum\limits_{i = 1}^{n} a_{i, j} =
				\sum\limits_{j = 1}^{n} a_{i, j} 
					= 1& (1 \leq i, j \leq n). 
		\end{align}
	\end{subequations}
Let $\ASM{n}$ denote the set of all ASM's of size $n$.
\end{definition}
The following is an example of ASM of size $4$: 
	\[
		\begin{pmatrix}
			0	&	1	&	0	&	0	\\
			1	&	- 1	&	1	&	0	\\
			0	&	1	&	0	&	0	\\
			0	&	0	&	0	&	1
		\end{pmatrix}
	\]. 
\subsection{Six vertex model and Fully packed loop model}\label{subsection:SVM_FPL}
Now we define the six vertex model and fully packed loop model
for the sake we provide a rigorous approach.
First we give the definition of a planary (simple and finite) graph which we name 
$
	{L}_{m, n} = 
		\left(
			V({L}_{m, n}), E({L}_{m, n})
		\right)
$.
The vertex set
$
	V({L}_{m, n}) = {V}_{0}(m, n) \sqcup {V}_{1}(m, n)
$
and the edge set
$
	E({L}_{m, n}) = {E}_{0}(m, n) \sqcup {E}_{1}(m, n)
$
are composed of the two kinds of sets, respectively.
Each vertex set is defined as

	\begin{equation*} 
		V_{0}(m, n) := \left\{
			(i, j) \in \mathbb{Z}^{2} 
		\, \middle| \, 
			1 \leq i \leq m,
			1 \leq j \leq n
		\right\} ,
	\end{equation*}
	which is called the set of  
	{\sl interior vertices}
and 
	\begin{equation*}
		V_{1}(m, n):= \left\{
			(i, j) 
		\, \middle| \, 
			i \in [m], \:
			j \in \{ 0, n + 1 \} 
		\right\}
				\sqcup
		\left\{
			(i, j) 
		\, \middle| \,
			i \in \{ 0, m + 1 \} , \:
			j \in [n]
		\right\} ,
	\end{equation*}
	which is called the set of  
	{\sl boundary vertices}. 
In graph theory a pair of vertices is called an edge.
Here we call a pair of the form
$\left\{(i, j),(i, j + 1)\right\}$ with $(i, j),(i, j + 1)\in V({L}_{m, n})$
 a {\sl horizontal edge},
and a pair of the form 
$\left\{(i, j),(i + 1, j)\right\}$ with $(i, j),(i + 1, j)\in V({L}_{m, n})$
 a {\sl vertical edge}.
The edge set $E(L_{m, n})$ of the graph $L_{m, n}$ is, by definition,
the set of all horizontal and vertical edges.
We also define ${E}_{0}(m, n)$ as the set of all edges,
whose endpoints are both interior vertices,
and ${E}_{1}(m, n)$ as the set of edges such that one of the endpoints is
a boundary vertex.
More precisely
	\begin{align}
		E_{0}(m, n) := &	
		\left\{
			\{(i, j), (i, j + 1) \}
		\, \middle| \, 
			1 \leq i \leq m, \, 
			1 \leq j \leq n - 1
		\right\} \nonumber \\
		\label{E0}
			 & \sqcup 
			 \left\{
				\{ (i, j), (i + 1, j) \}
			\, \middle| \, 
				1 \leq i \leq m - 1, \, 
				1 \leq j \leq n
			\right\} ,
	\end{align}
	and
	\begin{align}
		E_{1}(m, n)	:= & 
		\left\{
			\{ (i, j), (i, j + 1) \}
		\, \middle| \, 
			1 \leq i \leq m, \:
			j \in \{ 0, n \} 
		\right\} \nonumber \\
		\label{E1}
			& \sqcup
		\left\{
			\{ (i, j), (i + 1, j) \}
		\, \middle| \, 
			i \in \{ 0, m \} , \:
			1 \leq j \leq n
		\right\} .
	\end{align}
	We call an element of ${E}_{1}(m, n)$ 
	a {\sl boundary edge}.
	Especially, we denote ${L}_{n, n}$ as ${L}_{n}$. 
	We can regard the graph $L_{n}$ as the subset of $xy$-plane
	as in Figure~\ref{pic:L_3}.
	We label the boundary edges counterclockwise 
	with 
	$
		{e}_{1} = \left\{ (1, 0), (1, 1) \right\} , 
		{e}_{2}, \ldots 
	$. 
	starting from 
	${e}_{1}$. 
For example, Figure~\ref{pic:L_3} is the graph ${L}_{3}$,
in which the double circle dots are the boundary vetices.
\begin{figure}[htbp]
		\centering
		\begin{tikzpicture}[scale=0.75]
			\coordinate (base_1) at (0, -1);
			\coordinate (base_2) at (1, 0);
			\draw[->] (0, 0) -- ($4.5*(base_1)$);
				\node at ($4.5*(base_1)$)[below]{$i$};
			\draw[->] (0, 0) -- ($4.5*(base_2)$);
				\node at ($4.5*(base_2)$)[right]{$j$};
			\foreach \i in {1, 2, 3}{
				\foreach \j in {1, 2, 3}{
					\coordinate (v_\i_\j) at ($\i*(base_1) + \j*(base_2)$);
				}
			}
			\coordinate (b_1) at ($1*(base_2)$);
			\foreach \k/\i in {2/1, 3/2, 4/3}{
				\coordinate (b_\k) at ($\i*(base_1)$);
			}
			\foreach \k/\j in {5/1, 6/2, 7/3}{
				\coordinate (b_\k) at ($4*(base_1) + \j*(base_2)$);
			}
			\foreach \k/\i in {8/3, 9/2, 10/1}{
				\coordinate (b_\k) at ($\i*(base_1) + 4*(base_2)$);
			}
			\foreach \k/\j in {11/3, 12/2}{
				\coordinate (b_\k) at ($\j*(base_2)$);
			}
			\foreach \i in {1, 2, 3}{
				\foreach \j/\jj in {1/2, 2/3}{
					\draw ($\i*(base_1) + \j*(base_2)$) -- ($\i*(base_1) + \jj*(base_2)$);
				}
			}
			\foreach \j in {1, 2, 3}{
				\foreach \i/\ii in {1/2, 2/3}{
					\draw ($\i*(base_1) + \j*(base_2)$) -- ($\ii*(base_1) + \j*(base_2)$);
				}
			}
			\draw (v_1_1) -- (b_1);
			\foreach \k/\i in {2/1, 3/2, 4/3}{
				\draw (b_\k) -- (v_\i_1);
			}
			\foreach \k/\j in {5/1, 6/2, 7/3}{
				\draw (v_3_\j) -- (b_\k);
			}
			\foreach \k/\i in {8/3, 9/2, 10/1}{
				\draw (b_\k) -- (v_\i_3);
			}
			\foreach \k/\j in {11/3, 12/2}{
				\draw (v_1_\j) -- (b_\k);
			}
			\foreach \x in {1, 2, 3}{
				\foreach \y in {1, 2, 3}{
					\fill (v_\x_\y) circle (2pt);
				}
			}
			\foreach \k in {1, 2, ..., 12}{
				\fill[fill=white, draw=black] (b_\k) circle (4pt);
				\fill (b_\k) circle (2pt);
			}
			\node at (b_1)[above]{${e}_{1}$};
			\foreach \k in {2, 3, 4}{
				\node at (b_\k)[left]{${e}_{\k}$};
			}
			\foreach \k in {5, 6, 7}{
				\node at (b_\k)[below]{${e}_{\k}$};
			}
			\foreach \k in {8, 9, 10}{
				\node at (b_\k)[right]{${e}_{\k}$};
			}
			\foreach \k in {11, 12}{
				\node at (b_\k)[above]{${e}_{\k}$};
			}
		\end{tikzpicture}
		\caption{${L}_{3}$}
		\label{pic:L_3}
\end{figure}
\subsubsection{Six vertex model}
A six vertex model is a type of statistical mechanics model 
in which the Boltzmann weights are associated with each vertex in the model.
If this model has six possible states for each vertex,
we call it a {\sl six vertex model}.
A state of our model is obtained by giving a direction to each edge of $L_{m, n}$.
Recall that an edge of $L_{m, n}$ is denoted by an unorderd pair 
$\{ u, v \}$ of vertices.
We use ordered pair $(u,v)$ to denoted the direction from $u$ to $v$,
and $(v,u)$ to denote the opposite one.
We call a map which
associate $(u, v)$ or $(v, u)$ to $\left\{ u, v \right\}$
for each edge $\left\{ u, v \right\} \in E\left( {L}_{m, n} \right)$
\textsl{orientation of ${L}_{m, n}$}.
Here, if the orientation $\varphi$ satisfies $\varphi(\{ u, v \} ) = (u, v)$,
we say
$u$ is a {\sl source}\/  
or the directed edge $(u, v)$ goes {\sl out}\/ of $u$.
On the other hand, we say $v$ is a {\sl sink}\/
or  $(u, v)$ comes {\sl in}\/ $v$.
A {\sl state}\/ is, by definition,
a way asigning a direction to each edge.
Since
each interior vertex $v$ has exactly $4$ adjacent edges,
there are $2^4$ ways to orient these edges,
We say a state  
{\sl $2$-in-$2$-out} 
if there are $2$ edges in
and $2$ edges out for every vertex. 
Then 
we call a orientation
$\varphi$ 
a {\sl state} of 
{\sl six vertex model} on ${L}_{m, n}$
if each interior vertex is $2$-in-$2$-out. 
	When $\varphi$ is a state of six vertex model on ${L}_{m, n}$, 
	we call 
	$
		{\varphi}{|}_{
			{E}_{1}(m, n)
		}
	$  
	{\sl boundary condition} 
	of $\varphi$. 
We usually fix a boundary condition. 
Then the boundary condition on ${L}_{n}$ like Figure~\ref{pic:bc1} 
is called 
{\sl open boundary condition}. 
We denote the set of all six vertex model on ${L}_{n}$
which has open boundary condition
as $\mathcal{SV}(n)$.
\begin{figure}[htbp]
	\centering
	\begin{minipage}{0.45\hsize}
		\centering
		\begin{tikzpicture}[scale=0.75]
			\coordinate (base_1) at (0, -1);
			\coordinate (base_2) at (1, 0);
			\foreach \i in {1, 2, 3}{
				\foreach \j in {1, 2, 3}{
					\coordinate (v_\i_\j) at ($\i*(base_1) + \j*(base_2)$);
				}
			}
			\coordinate (b_1) at ($1*(base_2)$);
			\foreach \k/\i in {2/1, 3/2, 4/3}{
				\coordinate (b_\k) at ($\i*(base_1)$);
			}
			\foreach \k/\j in {5/1, 6/2, 7/3}{
				\coordinate (b_\k) at ($4*(base_1) + \j*(base_2)$);
			}
			\foreach \k/\i in {8/3, 9/2, 10/1}{
				\coordinate (b_\k) at ($\i*(base_1) + 4*(base_2)$);
			}
			\foreach \k/\j in {11/3, 12/2}{
				\coordinate (b_\k) at ($\j*(base_2)$);
			}
			\foreach \i in {1, 2, 3}{
				\foreach \j/\jj in {1/2, 2/3}{
					\draw ($\i*(base_1) + \j*(base_2)$) -- ($\i*(base_1) + \jj*(base_2)$);
				}
			}
			\foreach \j in {1, 2, 3}{
				\foreach \i/\ii in {1/2, 2/3}{
					\draw ($\i*(base_1) + \j*(base_2)$) -- ($\ii*(base_1) + \j*(base_2)$);
				}
			}
			\draw[directed] (v_1_1) -- (b_1);
			\foreach \k/\i in {2/1, 3/2, 4/3}{
				\draw[directed] (b_\k) -- (v_\i_1);
			}
			\foreach \k/\j in {5/1, 6/2, 7/3}{
				\draw[directed] (v_3_\j) -- (b_\k);
			}
			\foreach \k/\i in {8/3, 9/2, 10/1}{
				\draw[directed] (b_\k) -- (v_\i_3);
			}
			\foreach \k/\j in {11/3, 12/2}{
				\draw[directed] (v_1_\j) -- (b_\k);
			}
			\foreach \x in {1, 2, 3}{
				\foreach \y in {1, 2, 3}{
					\fill (v_\x_\y) circle (2pt);
				}
			}
			\foreach \k in {1, 2, ..., 12}{
				\fill[fill=white, draw=black] (b_\k) circle (4pt);
				\fill (b_\k) circle (2pt);
			}
			\node at (b_1)[above]{${e}_{1}$};
			\foreach \k in {2, 3, 4}{
				\node at (b_\k)[left]{${e}_{\k}$};
			}
			\foreach \k in {5, 6, 7}{
				\node at (b_\k)[below]{${e}_{\k}$};
			}
			\foreach \k in {8, 9, 10}{
				\node at (b_\k)[right]{${e}_{\k}$};
			}
			\foreach \k in {11, 12}{
				\node at (b_\k)[above]{${e}_{\k}$};
			}
		\end{tikzpicture}
		\caption{
			An example of
			boundary condition
		}
		\label{pic:bc1}
	\end{minipage}
	\begin{minipage}{0.45\hsize}
		\centering
		\begin{tikzpicture}[scale=0.75]
			\coordinate (base_1) at (0, -1);
			\coordinate (base_2) at (1, 0);
			\foreach \i in {1, 2, 3}{
				\foreach \j in {1, 2, 3}{
					\coordinate (v_\i_\j) at ($\i*(base_1) + \j*(base_2)$);
				}
			}
			\coordinate (b_1) at ($1*(base_2)$);
			\foreach \k/\i in {2/1, 3/2, 4/3}{
				\coordinate (b_\k) at ($\i*(base_1)$);
			}
			\foreach \k/\j in {5/1, 6/2, 7/3}{
				\coordinate (b_\k) at ($4*(base_1) + \j*(base_2)$);
			}
			\foreach \k/\i in {8/3, 9/2, 10/1}{
				\coordinate (b_\k) at ($\i*(base_1) + 4*(base_2)$);
			}
			\foreach \k/\j in {11/3, 12/2}{
				\coordinate (b_\k) at ($\j*(base_2)$);
			}
				\foreach \j/\jj in {1/2}{
					\draw[directed] ($1*(base_1) + \j*(base_2)$) -- ($1*(base_1) + \jj*(base_2)$);
				}
				\foreach \j/\jj in {2/3}{
					\draw[reverse directed] ($1*(base_1) + \j*(base_2)$) -- ($1*(base_1) + \jj*(base_2)$);
				}
				\foreach \j/\jj in {1/2, 2/3}{
					\draw[reverse directed] ($2*(base_1) + \j*(base_2)$) -- ($2*(base_1) + \jj*(base_2)$);
				}
				\foreach \j/\jj in {1/2, 2/3}{
					\draw[directed] ($3*(base_1) + \j*(base_2)$) -- ($3*(base_1) + \jj*(base_2)$);
				}
				\foreach \i/\ii in {2/3}{
					\draw[directed] ($\i*(base_1) + 1*(base_2)$) -- ($\ii*(base_1) + 1*(base_2)$);
				}
				\foreach \i/\ii in {1/2}{
					\draw[reverse directed] ($\i*(base_1) + 1*(base_2)$) -- ($\ii*(base_1) + 1*(base_2)$);
				}
				\foreach \i/\ii in {1/2, 2/3}{
					\draw[directed] ($\i*(base_1) + 2*(base_2)$) -- ($\ii*(base_1) + 2*(base_2)$);
				}
				\foreach \i/\ii in {1/2, 2/3}{
					\draw[reverse directed] ($\i*(base_1) + 3*(base_2)$) -- ($\ii*(base_1) + 3*(base_2)$);
				}
			\draw[directed] (v_1_1) -- (b_1);
			\foreach \k/\i in {2/1, 3/2, 4/3}{
				\draw[directed] (b_\k) -- (v_\i_1);
			}
			\foreach \k/\j in {5/1, 6/2, 7/3}{
				\draw[directed] (v_3_\j) -- (b_\k);
			}
			\foreach \k/\i in {8/3, 9/2, 10/1}{
				\draw[directed] (b_\k) -- (v_\i_3);
			}
			\foreach \k/\j in {11/3, 12/2}{
				\draw[directed] (v_1_\j) -- (b_\k);
			}
			\foreach \x in {1, 2, 3}{
				\foreach \y in {1, 2, 3}{
					\fill (v_\x_\y) circle (2pt);
				}
			}
			\foreach \k in {1, 2, ..., 12}{
				\fill[fill=white, draw=black] (b_\k) circle (4pt);
				\fill (b_\k) circle (2pt);
			}
			\node at (b_1)[above]{${e}_{1}$};
			\foreach \k in {2, 3, 4}{
				\node at (b_\k)[left]{${e}_{\k}$};
			}
			\foreach \k in {5, 6, 7}{
				\node at (b_\k)[below]{${e}_{\k}$};
			}
			\foreach \k in {8, 9, 10}{
				\node at (b_\k)[right]{${e}_{\k}$};
			}
			\foreach \k in {11, 12}{
				\node at (b_\k)[above]{${e}_{\k}$};
			}
		\end{tikzpicture}
		\caption{
			An example of a state of
			Six vertex model
		}
		\label{pic:svm1}
	\end{minipage}
\end{figure}
\subsubsection{Fully packed loop model}
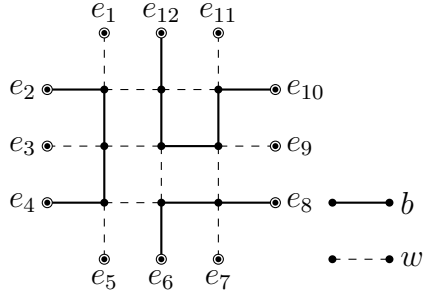
\begin{figure}[htbp]
		\centering
		\begin{tikzpicture}[scale=0.75]
			\coordinate (base_1) at (0, -1);
			\coordinate (base_2) at (1, 0);
			\foreach \x in {1, 2, 3}{
				\foreach \y in {1, 2, 3}{
					\coordinate (v_\x_\y) at ($\x*(base_1) + \y*(base_2)$);
				}
			}
			\foreach \k/\j in {11/3, 12/2, 1/1}{
				\coordinate (b_\k) at ($\j*(base_2)$);
					\node at (b_\k)[above]{${e}_{\k}$};
			}
			\foreach \k/\i in {2/1, 3/2, 4/3}{
				\coordinate (b_\k) at ($\i*(base_1)$);
					\node at (b_\k)[left]{${e}_{\k}$};
			}
			\foreach \k/\j in {5/1, 6/2, 7/3}{
				\coordinate (b_\k) at ($4*(base_1) + \j*(base_2)$);
					\node at (b_\k)[below]{${e}_{\k}$};
			}
			\foreach \k/\i in {8/3, 9/2, 10/1}{
				\coordinate (b_\k) at ($\i*(base_1) + 4*(base_2)$);
					\node at (b_\k)[right]{${e}_{\k}$};
			}
				\foreach \k/\j in {1/1, 11/3}{
					\draw[dashed] (v_1_\j) -- (b_\k);
				}
				\foreach \k/\j in {12/2}{
					\draw[thick] (v_1_\j) -- (b_\k);
				}
				\foreach \k/\i in {3/2}{
					\draw[dashed] (b_\k) -- (v_\i_1);
				}
				\foreach \k/\i in {2/1, 4/3}{
					\draw[thick] (b_\k) -- (v_\i_1);
				}
				\foreach \k/\j in {5/1, 7/3}{
					\draw[dashed] (v_3_\j) -- (b_\k);
				}
				\foreach \k/\j in {6/2}{
					\draw[thick] (v_3_\j) -- (b_\k);
				}
				\foreach \k/\i in {9/2}{
					\draw[dashed] (b_\k) -- (v_\i_3);
				}
				\foreach \k/\i in {8/3, 10/1}{
					\draw[thick] (b_\k) -- (v_\i_3);
				}
					\foreach \j/\jj in {1/2, 2/3}{
						\draw[dashed] (v_1_\j) -- (v_1_\jj);
					}
					\foreach \j/\jj in {1/2}{
						\draw[dashed] (v_2_\j) -- (v_2_\jj);
					}
					\foreach \j/\jj in {2/3}{
						\draw[thick] (v_2_\j) -- (v_2_\jj);
					}
					\foreach \j/\jj in {1/2}{
						\draw[dashed] (v_3_\j) -- (v_3_\jj);
					}
					\foreach \j/\jj in {2/3}{
						\draw[thick] (v_3_\j) -- (v_3_\jj);
					}
					\foreach \i/\ii in {1/2, 2/3}{
						\draw[thick] (v_\i_1) -- (v_\ii_1);
					}
					\foreach \i/\ii in {2/3}{
						\draw[dashed] (v_\i_2) -- (v_\ii_2);
					}
					\foreach \i/\ii in {1/2}{
						\draw[thick] (v_\i_2) -- (v_\ii_2);
					}
					\foreach \i/\ii in {2/3}{
						\draw[dashed] (v_\i_3) -- (v_\ii_3);
					}
					\foreach \i/\ii in {1/2}{
						\draw[thick] (v_\i_3) -- (v_\ii_3);
					}
				\foreach \x in {1, 2, 3}{
					\foreach \y in {1, 2, 3}{
						\fill (v_\x_\y) circle (2pt);
					}
				}
				\foreach \k in {1, 2, ..., 12}{
					\fill[fill=white, draw] (b_\k) circle (2.5pt);
					\fill (b_\k) circle (1.5pt);
				}
			\coordinate (b1) at ($3*(base_1) + 5*(base_2)$);
			\coordinate (b2) at ($3*(base_1) + 6*(base_2)$);
			\coordinate (w1) at ($4*(base_1) + 5*(base_2)$);
			\coordinate (w2) at ($4*(base_1) + 6*(base_2)$);
			\draw[thick] (b1) -- (b2);
			\draw[dashed] (w1) -- (w2);
			\foreach \x in {b1, b2, w1, w2}{
				\fill (\x) circle (2pt);
			}
			\node at (b2)[right]{$b$};
			\node at (w2)[right]{$w$};
		\end{tikzpicture}
		\caption{
			an example of FPL on ${L}_{3}$
		}
		\label{pic:fpl1}
\end{figure}

We define the following map 
$
	\psi \colon E({L}_{m, n}) \rightarrow \{ b, w \}
$.  
Here, 
$b$
comes from 
black and 
$w$ 
comes from 
white. 
Let $v$ be an interior vertex, 
then we say $v$ is 
{\sl
$2$-$2$-colored} 
if $2$ out of $4$ edges which is incident to $v$ are $b$, 
and the rest are $w$. 
We also call a 
$\psi$ 
{\sl
fully packed loop model} on ${L}_{m, n}$ 
if each interior vertex is $2$-$2$-colored, 
and it is abbreviated as FPL. 
Figure \ref{pic:fpl1} is an example of FPL on ${L}_{3}$.
In Figure \ref{pic:fpl1}, 
we draw edge $e$ with solid line
(resp. dashed line)
when $e$ has color $b$
(resp. $w$)
for each edge $e \in E({L}_{3})$.

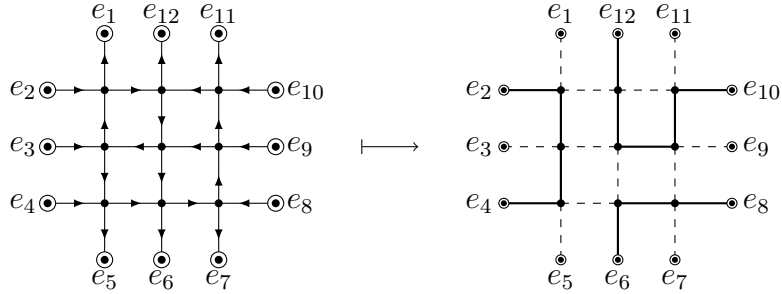
\begin{figure}[htbp]
	\centering
		\begin{tikzpicture}[scale=0.75]
			\coordinate (base_1) at (0, -1);
			\coordinate (base_2) at (1, 0);
			\foreach \i in {1, 2, 3}{
				\foreach \j in {1, 2, 3}{
					\coordinate (v_\i_\j) at ($\i*(base_1) + \j*(base_2) - 6*(base_2)$);
				}
			}
			\coordinate (b_1) at ($1*(base_2) - 6*(base_2)$);
			\foreach \k/\i in {2/1, 3/2, 4/3}{
				\coordinate (b_\k) at ($\i*(base_1) - 6*(base_2)$);
			}
			\foreach \k/\j in {5/1, 6/2, 7/3}{
				\coordinate (b_\k) at ($4*(base_1) + \j*(base_2) - 6*(base_2)$);
			}
			\foreach \k/\i in {8/3, 9/2, 10/1}{
				\coordinate (b_\k) at ($\i*(base_1) + 4*(base_2) - 6*(base_2)$);
			}
			\foreach \k/\j in {11/3, 12/2}{
				\coordinate (b_\k) at ($\j*(base_2) - 6*(base_2)$);
			}
				\foreach \j/\jj in {1/2}{
					\draw[directed] (v_1_\j) -- (v_1_\jj);
				}
				\foreach \j/\jj in {2/3}{
					\draw[reverse directed] (v_1_\j) -- (v_1_\jj);
				}
				\foreach \j/\jj in {1/2, 2/3}{
					\draw[reverse directed] (v_2_\j) -- (v_2_\jj);
				}
				\foreach \j/\jj in {1/2, 2/3}{
					\draw[directed] (v_3_\j) -- (v_3_\jj);
				}
				\foreach \i/\ii in {2/3}{
					\draw[directed] (v_\i_1) -- (v_\ii_1);
				}
				\foreach \i/\ii in {1/2}{
					\draw[reverse directed] (v_\i_1) -- (v_\ii_1);
				}
				\foreach \i/\ii in {1/2, 2/3}{
					\draw[directed] (v_\i_2) -- (v_\ii_2);
				}
				\foreach \i/\ii in {1/2, 2/3}{
					\draw[reverse directed] (v_\i_3) -- (v_\ii_3);
				}
			\draw[directed] (v_1_1) -- (b_1);
			\foreach \k/\i in {2/1, 3/2, 4/3}{
				\draw[directed] (b_\k) -- (v_\i_1);
			}
			\foreach \k/\j in {5/1, 6/2, 7/3}{
				\draw[directed] (v_3_\j) -- (b_\k);
			}
			\foreach \k/\i in {8/3, 9/2, 10/1}{
				\draw[directed] (b_\k) -- (v_\i_3);
			}
			\foreach \k/\j in {11/3, 12/2}{
				\draw[directed] (v_1_\j) -- (b_\k);
			}
			\foreach \x in {1, 2, 3}{
				\foreach \y in {1, 2, 3}{
					\fill (v_\x_\y) circle (2pt);
				}
			}
			\foreach \k in {1, 2, ..., 12}{
				\fill[fill=white, draw=black] (b_\k) circle (4pt);
				\fill (b_\k) circle (2pt);
			}
			\node at (b_1)[above]{${e}_{1}$};
			\foreach \k in {2, 3, 4}{
				\node at (b_\k)[left]{${e}_{\k}$};
			}
			\foreach \k in {5, 6, 7}{
				\node at (b_\k)[below]{${e}_{\k}$};
			}
			\foreach \k in {8, 9, 10}{
				\node at (b_\k)[right]{${e}_{\k}$};
			}
			\foreach \k in {11, 12}{
				\node at (b_\k)[above]{${e}_{\k}$};
			}
		\draw[|->] ($2*(base_1) - 0.5*(base_2)$) -- ($2*(base_1) + 0.5*(base_2)$);
			\foreach \x in {1, 2, 3}{
				\foreach \y in {1, 2, 3}{
					\coordinate (w_\x_\y) at ($\x*(base_1) + \y*(base_2) + 2*(base_2)$);
				}
			}
			\foreach \k/\j in {11/3, 12/2, 1/1}{
				\coordinate (p_\k) at ($\j*(base_2) + 2*(base_2)$);
					\node at (p_\k)[above]{${e}_{\k}$};
			}
			\foreach \k/\i in {2/1, 3/2, 4/3}{
				\coordinate (p_\k) at ($\i*(base_1) + 2*(base_2)$);
					\node at (p_\k)[left]{${e}_{\k}$};
			}
			\foreach \k/\j in {5/1, 6/2, 7/3}{
				\coordinate (p_\k) at ($4*(base_1) + \j*(base_2) + 2*(base_2)$);
					\node at (p_\k)[below]{${e}_{\k}$};
			}
			\foreach \k/\i in {8/3, 9/2, 10/1}{
				\coordinate (p_\k) at ($\i*(base_1) + 4*(base_2) + 2*(base_2)$);
					\node at (p_\k)[right]{${e}_{\k}$};
			}
				\foreach \k/\j in {1/1, 11/3}{
					\draw[dashed] (w_1_\j) -- (p_\k);
				}
				\foreach \k/\j in {12/2}{
					\draw[thick] (w_1_\j) -- (p_\k);
				}
				\foreach \k/\i in {3/2}{
					\draw[dashed] (p_\k) -- (w_\i_1);
				}
				\foreach \k/\i in {2/1, 4/3}{
					\draw[thick] (p_\k) -- (w_\i_1);
				}
				\foreach \k/\j in {5/1, 7/3}{
					\draw[dashed] (w_3_\j) -- (p_\k);
				}
				\foreach \k/\j in {6/2}{
					\draw[thick] (w_3_\j) -- (p_\k);
				}
				\foreach \k/\i in {9/2}{
					\draw[dashed] (p_\k) -- (w_\i_3);
				}
				\foreach \k/\i in {8/3, 10/1}{
					\draw[thick] (p_\k) -- (w_\i_3);
				}
					\foreach \j/\jj in {1/2, 2/3}{
						\draw[dashed] (w_1_\j) -- (w_1_\jj);
					}
					\foreach \j/\jj in {1/2}{
						\draw[dashed] (w_2_\j) -- (w_2_\jj);
					}
					\foreach \j/\jj in {2/3}{
						\draw[thick] (w_2_\j) -- (w_2_\jj);
					}
					\foreach \j/\jj in {1/2}{
						\draw[dashed] (w_3_\j) -- (w_3_\jj);
					}
					\foreach \j/\jj in {2/3}{
						\draw[thick] (w_3_\j) -- (w_3_\jj);
					}
					\foreach \i/\ii in {1/2, 2/3}{
						\draw[thick] (w_\i_1) -- (w_\ii_1);
					}
					\foreach \i/\ii in {2/3}{
						\draw[dashed] (w_\i_2) -- (w_\ii_2);
					}
					\foreach \i/\ii in {1/2}{
						\draw[thick] (w_\i_2) -- (w_\ii_2);
					}
					\foreach \i/\ii in {2/3}{
						\draw[dashed] (w_\i_3) -- (w_\ii_3);
					}
					\foreach \i/\ii in {1/2}{
						\draw[thick] (w_\i_3) -- (w_\ii_3);
					}
				\foreach \x in {1, 2, 3}{
					\foreach \y in {1, 2, 3}{
						\fill (w_\x_\y) circle (2pt);
					}
				}
				\foreach \k in {1, 2, ..., 12}{
					\fill[fill=white, draw] (p_\k) circle (2.5pt);
					\fill (p_\k) circle (1.5pt);
				}
		\end{tikzpicture}
	\caption{
		An example of assignment a six vertex model to a FPL
	}
	\label{pic:SVM_to_FPL_example}
\end{figure}

Here, let $\psi$ be a FPL on ${L}_{m, n}$ and 
${e}_{i}$ be the $i$-th boundary edge, 
then
\[
	\tau =
		\left(
			\psi({e}_{1}), \psi({e}_{2}), \ldots , \psi({e}_{2m + 2n})
		\right)
\]
is called 
boundary condition of $\psi$. 
We also denote $\mathfrak{fpl}(n, \tau)$ 
the set of all FPL on ${L}_{n}$ with a fixed boundary condition $\tau$, 
when $\tau \in { \{ b, w \} }^{4n}$ is given. 
Now, we denote the following boundary conditions on ${L}_{n}$ as 
${\tau}_{+}$ and ${\tau}_{-}$ respectively: 
	\begin{subequations}
 		\begin{align} \label{eq:tp}
 			{\tau}_{+} 
 				&:= \left( b, w, b, w, \cdots , b, w \right) ,
 		\\
 		\label{eq:tm}
 			{\tau}_{-} 
 				&:= \left( w, b, w, b, \cdots , w, b \right). 
 		\end{align}
 	\end{subequations}
These conditions have $b$ and $w$ arranged alternately. 
We usually 
consider $\mathfrak{fpl}(n, {\tau}_{-})$. 
\begin{figure}[htbp]
	\begin{minipage}{0.45\hsize}
		\centering
		\begin{tikzpicture}[scale=0.75]
			\coordinate (base_1) at (0, -1);
			\coordinate (base_2) at (1, 0);
			\foreach \x in {1, 2, 3}{
				\foreach \y in {1, 2, 3}{
					\coordinate (v_\x_\y) at ($\x*(base_1) + \y*(base_2)$);
				}
			}
			\foreach \k/\j in {11/3, 12/2, 1/1}{
				\coordinate (b_\k) at ($\j*(base_2)$);
					\node at (b_\k)[above]{${e}_{\k}$};
			}
			\foreach \k/\i in {2/1, 3/2, 4/3}{
				\coordinate (b_\k) at ($\i*(base_1)$);
					\node at (b_\k)[left]{${e}_{\k}$};
			}
			\foreach \k/\j in {5/1, 6/2, 7/3}{
				\coordinate (b_\k) at ($4*(base_1) + \j*(base_2)$);
					\node at (b_\k)[below]{${e}_{\k}$};
			}
			\foreach \k/\i in {8/3, 9/2, 10/1}{
				\coordinate (b_\k) at ($\i*(base_1) + 4*(base_2)$);
					\node at (b_\k)[right]{${e}_{\k}$};
			}
				\foreach \k/\j in {1/1, 11/3}{
					\draw[thick] (v_1_\j) -- (b_\k);
				}
				\foreach \k/\j in {12/2}{
					\draw[dashed] (v_1_\j) -- (b_\k);
				}
				\foreach \k/\i in {3/2}{
					\draw[thick] (b_\k) -- (v_\i_1);
				}
				\foreach \k/\i in {2/1, 4/3}{
					\draw[dashed] (b_\k) -- (v_\i_1);
				}
				\foreach \k/\j in {5/1, 7/3}{
					\draw[thick] (v_3_\j) -- (b_\k);
				}
				\foreach \k/\j in {6/2}{
					\draw[dashed] (v_3_\j) -- (b_\k);
				}
				\foreach \k/\i in {9/2}{
					\draw[thick] (b_\k) -- (v_\i_3);
				}
				\foreach \k/\i in {8/3, 10/1}{
					\draw[dashed] (b_\k) -- (v_\i_3);
				}
				\foreach \i in {1, 2, 3}{
					\foreach \j in {1, 2, 3}{
						\fill (v_\i_\j) circle (2pt);
					}
				}
				\foreach \k in {1, 2, ..., 12}{
					\fill[fill=white, draw] (b_\k) circle (2.5pt);
					\fill (b_\k) circle (1.5pt);
				}
			\coordinate (b1) at ($5.5*(base_1) + 0.5*(base_2)$);
			\coordinate (b2) at ($5.5*(base_1) + 1.5*(base_2)$);
			\coordinate (w1) at ($5.5*(base_1) + 2.5*(base_2)$);
			\coordinate (w2) at ($5.5*(base_1) + 3.5*(base_2)$);
			\draw[thick] (b1) -- (b2);
			\draw[dashed] (w1) -- (w2);
			\foreach \x in {b1, b2, w1, w2}{
				\fill (\x) circle (2pt);
			}
			\node at (b2)[right]{$b$};
			\node at (w2)[right]{$w$};
		\end{tikzpicture}
		\caption{
			${\tau}_{+}$ when $n = 3$
		}
		\label{pic:tau_plus_n3}
	\end{minipage}
	\begin{minipage}{0.45\hsize}
		\centering
		\begin{tikzpicture}[scale=0.75]
			\coordinate (base_1) at (0, -1);
			\coordinate (base_2) at (1, 0);
			\foreach \x in {1, 2, 3}{
				\foreach \y in {1, 2, 3}{
					\coordinate (v_\x_\y) at ($\x*(base_1) + \y*(base_2)$);
				}
			}
			\foreach \k/\j in {11/3, 12/2, 1/1}{
				\coordinate (b_\k) at ($\j*(base_2)$);
					\node at (b_\k)[above]{${e}_{\k}$};
			}
			\foreach \k/\i in {2/1, 3/2, 4/3}{
				\coordinate (b_\k) at ($\i*(base_1)$);
					\node at (b_\k)[left]{${e}_{\k}$};
			}
			\foreach \k/\j in {5/1, 6/2, 7/3}{
				\coordinate (b_\k) at ($4*(base_1) + \j*(base_2)$);
					\node at (b_\k)[below]{${e}_{\k}$};
			}
			\foreach \k/\i in {8/3, 9/2, 10/1}{
				\coordinate (b_\k) at ($\i*(base_1) + 4*(base_2)$);
					\node at (b_\k)[right]{${e}_{\k}$};
			}
				\foreach \k/\j in {1/1, 11/3}{
					\draw[dashed] (v_1_\j) -- (b_\k);
				}
				\foreach \k/\j in {12/2}{
					\draw[thick] (v_1_\j) -- (b_\k);
				}
				\foreach \k/\i in {3/2}{
					\draw[dashed] (b_\k) -- (v_\i_1);
				}
				\foreach \k/\i in {2/1, 4/3}{
					\draw[thick] (b_\k) -- (v_\i_1);
				}
				\foreach \k/\j in {5/1, 7/3}{
					\draw[dashed] (v_3_\j) -- (b_\k);
				}
				\foreach \k/\j in {6/2}{
					\draw[thick] (v_3_\j) -- (b_\k);
				}
				\foreach \k/\i in {9/2}{
					\draw[dashed] (b_\k) -- (v_\i_3);
				}
				\foreach \k/\i in {8/3, 10/1}{
					\draw[thick] (b_\k) -- (v_\i_3);
				}
				\foreach \i in {1, 2, 3}{
					\foreach \j in {1, 2, 3}{
						\fill (v_\i_\j) circle (2pt);
					}
				} 
				\foreach \k in {1, 2, ..., 12}{
					\fill[fill=white, draw] (b_\k) circle (2.5pt);
					\fill (b_\k) circle (1.5pt);
				}
			\coordinate (b1) at ($5.5*(base_1) + 0.5*(base_2)$);
			\coordinate (b2) at ($5.5*(base_1) + 1.5*(base_2)$);
			\coordinate (w1) at ($5.5*(base_1) + 2.5*(base_2)$);
			\coordinate (w2) at ($5.5*(base_1) + 3.5*(base_2)$);
			\draw[thick] (b1) -- (b2);
			\draw[dashed] (w1) -- (w2);
			\foreach \x in {b1, b2, w1, w2}{
				\fill (\x) circle (2pt);
			}
			\node at (b2)[right]{$b$};
			\node at (w2)[right]{$w$};
		\end{tikzpicture}
		\caption{
			${\tau}_{-}$ when $n = 3$
		}
		\label{pic:tau_minus_n3}
	\end{minipage}
\end{figure}
It is well-known that 
there is the correopondence between
the 
six vertex model on ${L}_{n}$ which has the open boundary condition
and 
the FPL model on ${L}_{n}$ with ${\tau}_{-}$. 
Here, we explain the well-known correspondence 
from a state $\varphi$ to a FPL $\psi$. 
Firstly, 
we shall define parity of vertex. 
When $(i, j) \in V({L}_{n})$, 
we call $(i, j)$ odd if 
$i + j$ is odd, 
even otherwise.  
Notice that 
each edge of ${L}_{n}$ is incident to an odd vertex and an even vertex. 
For each edge $\left\{ u, v \right\} \in E \left( {L}_{n} \right)$
which goes out of $u$
(i.e., 
$
	{\varphi}(\{ u, v \} ) = (u, v)
$),
we set $\psi(\{ u, v \} ) = b$
(resp. $w$)
if $u$ is odd vertex
(resp. even).  
In Figure \ref{pic:SVM_to_FPL_example},
we draw odd vertex (resp. even)
as $\circ$ (resp. $\bullet$).
\subsubsection{Plaquette}
We define 
subgraph 
of ${L}_{n}$ 
called plaquette.
For $0 \leq i, j \leq n$, 
we define 
$
	{\alpha}_{i, j} = \left(
		V({\alpha}_{i, j}), E({\alpha}_{i, j})
	\right)
$ 
as follows: 
\begin{subequations}
\begin{align}
	V({\alpha}_{i, j})	&:=
		\left\{
			(i, j), (i, j + 1), (i + 1, j), (i + 1, j + 1)
		\right\} 
		\cap V({L}_{n}), 
		\\
	E({\alpha}_{i, j})	&:=
		\left\{
			e \in E({L}_{n})
		\; \middle| \;
			e \textrm{ incdents } u, v
				\; (u, v \in V({\alpha}_{i, j}))
		\right\} . 
\end{align}
\end{subequations}
We call ${\alpha}_{i, j}$ interior plaquette if 
any vertex of ${\alpha}_{i, j}$ is interior vertex, 
boundary plaquette otherwise. 
Now, ${\alpha}_{i, j}$ is interior plaquette if and only if 
$1 \leq i, j \leq n - 1$. 
We also define parity of plaquette ${\alpha}_{i, j}$. 
We call ${\alpha}_{i, j}$ odd if $i + j$ is odd, even otherwise. 
\subsection{The correpondence between six vertex model and ASM}
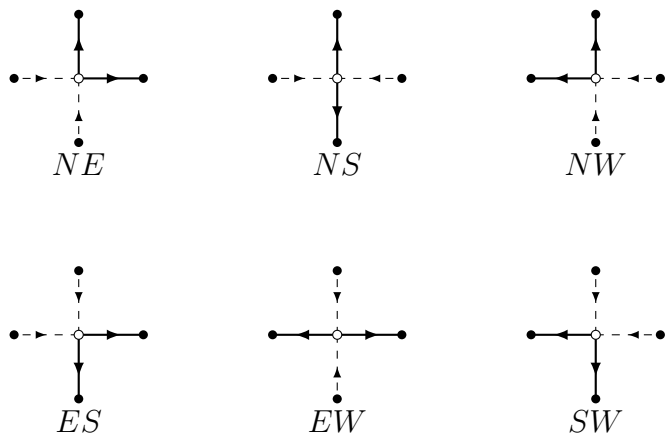
\begin{figure}[htbp]
		\centering
		\begin{tikzpicture}[scale=0.85]
			\coordinate (base_1) at (0, -1);
			\coordinate (base_2) at (1, 0);
			\coordinate (a) at ($1*(base_1) + 1*(base_2)$);
			\coordinate (b) at ($1*(base_1) + 5*(base_2)$);
			\coordinate (c) at ($1*(base_1) + 9*(base_2)$);
			\coordinate (d) at ($5*(base_1) + 1*(base_2)$);
			\coordinate (e) at ($5*(base_1) + 5*(base_2)$);
			\coordinate (f) at ($5*(base_1) + 9*(base_2)$);
			\foreach \p in {a, b, c, d, e, f}{
				\coordinate (\p_1) at ($(\p) - (base_1)$);
				\coordinate (\p_2) at ($(\p) + (base_2)$);
				\coordinate (\p_3) at ($(\p) + (base_1)$);
				\coordinate (\p_4) at ($(\p) - (base_2)$);
			}
			\foreach \x in {1, 2}{
				\draw[thick, black, directed] (a) -- (a_\x);
			}
			\foreach \x in {3, 4}{
				\draw[dashed, reverse directed] (a) -- (a_\x);
			}
			\node at (a_3)[below]{$NE$};
			\foreach \x in {1, 3}{
				\draw[thick, black, directed] (b) -- (b_\x);
			}
			\foreach \x in {2, 4}{
				\draw[dashed, reverse directed] (b) -- (b_\x);
			}
			\node at (b_3)[below]{$NS$};
			\foreach \x in {1, 4}{
				\draw[thick, black, directed] (c) -- (c_\x);
			}
			\foreach \x in {2, 3}{
				\draw[dashed, reverse directed] (c) -- (c_\x);
			}
			\node at (c_3)[below]{$NW$};
			\foreach \x in {2, 3}{
				\draw[thick, black, directed] (d) -- (d_\x);
			}
			\foreach \x in {1, 4}{
				\draw[dashed, reverse directed] (d) -- (d_\x);
			}
			\node at (d_3)[below]{$ES$};
			\foreach \x in {2, 4}{
				\draw[thick, black, directed] (e) -- (e_\x);
			}
			\foreach \x in {1, 3}{
				\draw[dashed, reverse directed] (e) -- (e_\x);
			}
			\node at (e_3)[below]{$EW$};
			\foreach \x in {3, 4}{
				\draw[thick, black, directed] (f) -- (f_\x);
			}
			\foreach \x in {1, 2}{
				\draw[dashed, reverse directed] (f) -- (f_\x);
			}
			\node at (f_3)[below]{$SW$};
			\foreach \x in {a, b, c, d, e, f}{
				\fill[fill=white, draw] (\x) circle (2pt);
				\foreach \y in {1, 2, 3, 4}{
					\fill (\x_\y) circle (2pt);
				}
			}
		\end{tikzpicture}
		\caption{The 6 possible choise of an interior vertex}
		\label{pic:vertex_config_SVM}
	\end{figure}

It is well known that there is a bijection
$\mathcal{SV}(n)$ to $\mathcal{A}_{n}$.  
We shall state the correpondence between six vertex model and ASM
in this section.
Notice that
each interior vertex $(i, j)$ has $6$ possible choise
when $(i, j)$ is $2$-in-$2$-out. 
Here, we respectively set the $4$ edges which is incident to interior vertex $(i, j)$
as $N = \left\{ (i, j), (i - 1, j) \right\}$, $E = \left\{ (i, j), (i, j + 1) \right\}$,
$S = \left\{ (i, j), (i + 1, j) \right\}$, $W = \left\{ (i, j), (i, j - 1) \right\}$.
Then, the $6$ possible choise can be expressed as
$NE$, $NS$, $NW$, $ES$, $EW$, $SW$
by specifying two edges goes out of $(i, j)$.
In Figure \ref{pic:vertex_config_SVM},
we illustrate the $6$ possible choise of an interior vertex.

We shall asign a state of six vertex model on ${L}_{n}$
to a square matrix of degree $n$
by using 
the choise of the interiror vertex $(i, j)$
to determine 
the $(i, j)$-entry 
for each $i, j$ which satisfy $1 \leq i, j \leq n$.
\begin{definition}\label{def:SV_to_ASM}
Let $n$ be an positive integer. 
We define a map 
$
	f \colon \mathcal{SV}(n) \rightarrow \mathcal{A}_{n}
$.
When we put $f \left( \varphi \right)$ as ${ \left( {a}_{i,j} \right) }_{1 \leq i, j \leq n}$
for a $\varphi \in \mathcal{SV}(n)$,
we define ${a}_{i, j}$
as follows:
${a}_{i, j} = 1$ if the choice of $(i, j)$ is $NS$,
${a}_{i, j} = -1$ if the choice of $(i, j)$ is $EW$,
${a}_{i, j} = 0$ otherwise
for each $1 \leq i, j \leq n$.
\end{definition}
We shall show the matrix 
${ \left( {a}_{i, j} \right) }_{1 \leq i, j \leq n}$
which is obtained by definition \ref{def:SV_to_ASM}
satisfies the conditions of ASM.
First, we fix $i \in [n]$ 
arbitraily.
There exists $j \in [n]$ such that
both of the horizontal edges 
$\left\{ (i, j), (i, j - 1) \right\}$ and $\left\{ (i, j), (i, j + 1) \right\}$ 
come in or go out of $(i, j)$.
In fact,
if $\left\{ (i, j), (i, j - 1) \right\}$ come in (resp. go out of) $(i, j)$
and $\left\{ (i, j), (i, j + 1) \right\}$ go out of (resp. come in) $(i, j)$
for each $1 \leq j \leq n$,
it contradicts the boundary condition for 
$\left\{ (i, n), (i, n + 1) \right\}$
(resp. $\left\{ (i, 0), (i, 1) \right\}$).
Moreover,
both of the horizontal edge 
$\left\{ (i, j), (i, j - 1) \right\}$and $\left\{ (i, j), (i, j + 1) \right\}$ 
come in $(i, j)$
for the smallest and the largest $j$
such that 
both of the horizontal edge 
$\left\{ (i, j), (i, j - 1) \right\}$and $\left\{ (i, j), (i, j + 1) \right\}$ 
come in or go out of $(i, j)$.
Then
we remark that
the choise of $(i, j)$ is $NS$ (resp. $EW$)
if both of the horizontal edge 
$\left\{ (i, j), (i, j - 1) \right\}$ and $\left\{ (i, j), (i, j + 1) \right\}$
come in (resp. go out of) $(i, j)$,
and ${a}_{i, j}$ is equal to $0$
if either of the horizontal edge 
$\left\{ (i, j), (i, j - 1) \right\}$ and $\left\{ (i, j), (i, j + 1) \right\}$
comes in $(i, j)$.
Since
a vertex 
such that
both of the horizontal edges  
come in the vertex 
and
a vertex 
such that
both of the horizontal edges 
go out of the vertex 
appears 
except a vertices
such that
either of the horizontal edges
come in the vertex 
in the $i$-th row,
the matrix ${ \left( {a}_{i, j} \right) }_{1 \leq i, j \leq n}$ satisfies
$\sum_{1 \leq l \leq n} {a}_{i, l} = 1$ and
$\sum_{1 \leq l \leq j} {a}_{i, l} \in \left\{ 0, 1 \right\}$
for each $1 \leq j \leq n$.
Then the matrix
${ \left( {a}_{i, j} \right) }_{1 \leq i, j \leq n}$
satisfies the condition on row of ASM.
The same goes for columns.

Now, we denote $\sum_{1 \leq k \leq i} {a}_{k, j}$ as ${c}_{i, j}$
and $\sum_{1 \leq l \leq j} {a}_{i, l}$ as ${r}_{i, j}$
for each $(i, j)$.
Then, 
we shall show that
the choice of $(i, j)$ can be determined by the trilpet $({a}_{i,j}, {c}_{i, j}, {r}_{i, j})$
for each interior veretex $(i, j)$.
First, we state that
there are exactly 6 possible value for $({a}_{i, j}, {c}_{i, j}, {r}_{i, j})$.
When ${a}_{i, j} = 0$,
$(i, j)$ can be devided into following the four cases:
\begin{enumerate}[(i)]
	\item	\label{case:NE}
		There exists ${i}_{1}$ greater than $i$ 
		and ${j}_{1}$ greater than $j$
		such that
		${a}_{{i}_{1}, j} = {a}_{i, {j}_{1}} = 1$,
		${a}_{k, j} = 0$ for $i \leq k < {i}_{1}$ and
		${a}_{i, l} = 0$ for $j \leq l < {j}_{1}$,
	\item	\label{case:NW}
		There exists ${i}_{1}$ greater than $i$ 
		and ${j}_{0}$ less than $j$
		such that
		${a}_{{i}_{1}, j} = {a}_{i, {j}_{0}} = 1$,
		${a}_{k, j} = 0$ for $i \leq k < {i}_{1}$ and
		${a}_{i, l} = 0$ for ${j}_{0} < l \leq j$,
	\item	\label{case:ES}
		There exists ${i}_{0}$ less than $i$ 
		and ${j}_{1}$ greater than $j$
		such that
		${a}_{{i}_{0}, j} = {a}_{i, {j}_{1}} = 1$,
		${a}_{k, j} = 0$ for ${i}_{0} < k \leq i$ and
		${a}_{i, l} = 0$ for $j \leq l < {j}_{1}$,
	\item	\label{case:SW}
		There exists ${i}_{0}$ less than $i$ 
		and ${j}_{0}$ less than $j$
		such that
		${a}_{{i}_{0}, j} = {a}_{i, {j}_{0}} = 1$,
		${a}_{k, j} = 0$ for ${i}_{0} < k \leq i$ and
		${a}_{i, l} = 0$ for ${j}_{0} < l \leq j$.
\end{enumerate}
In each four cases,
the triplet$({a}_{i, j}, {c}_{i, j}, {r}_{i, j})$ equals
$(0, 0, 0)$ when the case \eqref{case:NE},
$(0, 0, 1)$ when the case \eqref{case:NW},
$(0, 1, 0)$ when the case \eqref{case:ES} and 
$(0, 1, 1)$ when the case \eqref{case:SW}.
On the other hand,
the triplet $({a}_{i, j}, {c}_{i, j}, {r}_{i, j})$ must be $(1, 1, 1)$ (resp. $(-1, 0, 0)$)
when ${a}_{i, j} = 1$ (resp. ${a}_{i, j} = -1$)
because
$1$ apper first except $0$, and
$1$ and $-1$ alternatly appear except $0$
for each row and column in ASM.
Then, the $6$ possible value for the triplet $({a}_{i, j}, {c}_{i, j}, {r}_{i, j})$ are 
$(0, 0, 0)$, $(0, 0, 1)$, $(0, 1, 0)$, $(0, 1, 1)$, $(1, 1, 1)$ and $(-1, 0, 0)$.

Now, 
the vertical edge
$\left\{ (i, j), (i - 1, j) \right\}$ goes out of (resp. comes in) $(i, j)$ and
$\left\{ (i, j), (i + 1, j) \right\}$ comes in (resp. goes out of) $(i, j)$
if there  exists ${i}_{1}$ greater (resp. ${i}_{0}$ less) than $i$
such that
${a}_{{i}_{1}, j} = 1$ (resp. ${a}_{ {i}_{0}, j } = 1$) and
${a}_{k, j} = 0$ for $i \leq k < {i}_{1}$ (resp. ${i}_{0} < k \leq i$)
because
the vertical edge $\left\{ ({i}_{1}, j), ({i}_{1} - 1, j) \right\}$ 
(resp. $\left\{ ({i}_{0}, j), ({i}_{0} + 1, j) \right\}$)
goes out of $({i}_{1}, j)$
(resp. $({i}_{0}, j)$)
and
either of the vertical edges 
$\left\{ (k, j), (k - 1, j) \right\}$ and $\left\{ (k, j), (k + 1, j) \right\}$
comes in $(k, j)$
for $i \leq k < {i}_{1}$
(resp. ${i}_{0} < k \leq i$).
On the other hand, 
the horizontal edge
$\left\{ (i, j), (i, j - 1) \right\}$ comes in (goes out of) of $(i, j)$ and
$\left\{ (i, j), (i, j + 1) \right\}$ goes out of (comes in) $(i, j)$
if there  exists ${j}_{1}$ greater (resp. ${j}_{0}$ less) than $j$
such that
${a}_{i, {j}_{1}} = 1$ (resp. ${a}_{ i, {j}_{0} } = 1$) and
${a}_{i, l} = 0$ for $j \leq l < {j}_{1}$ (resp. ${j}_{0} < l \leq j$)
because
the horizontal edge $\left\{ (i, {j}_{1}), (i, {j}_{1} - 1) \right\}$ 
(resp. $\left\{ (i, {j}_{0}), (i, {j}_{0} + 1) \right\}$)
comes in $(i, {j}_{1})$
(resp. $({i}_{0}, j)$)
and
either of the vertical edges 
$\left\{ (i, l), (i, l - 1) \right\}$ and $\left\{ (i, l), (i, l + 1) \right\}$
comes in $(i, l)$
for $j \leq l < {j}_{1}$
(resp. ${j}_{0} < l \leq j$).
Therefore
the choice of $(i, j)$ is 
$NE$ when 
$({a}_{i, j}, {c}_{i, j}, {r}_{i, j}) = (0, 0, 0)$,
$NW$ when 
$({a}_{i, j}, {c}_{i, j}, {r}_{i, j}) = (0, 0, 1)$,
$ES$ when 
$({a}_{i, j}, {c}_{i, j}, {r}_{i, j}) = (0, 1, 0)$, 
$SW$ when 
$({a}_{i, j}, {c}_{i, j}, {r}_{i, j}) = (0, 1, 1)$,
$NS$ when 
$({a}_{i, j}, {c}_{i, j}, {r}_{i, j}) = (1, 1, 1)$ and
$EW$ when
$({a}_{i, j}, {c}_{i, j}, {r}_{i, j}) = (-1, 0, 0)$.
Since the triplet $({a}_{i, j}, {c}_{i, j}, {r}_{i, j})$ are consitent for each $1 \leq i, j \leq n$
if the image are consistent, 
the map $f \colon \mathcal{SV}(n) \rightarrow \mathcal{A}_{n}$ is injective.
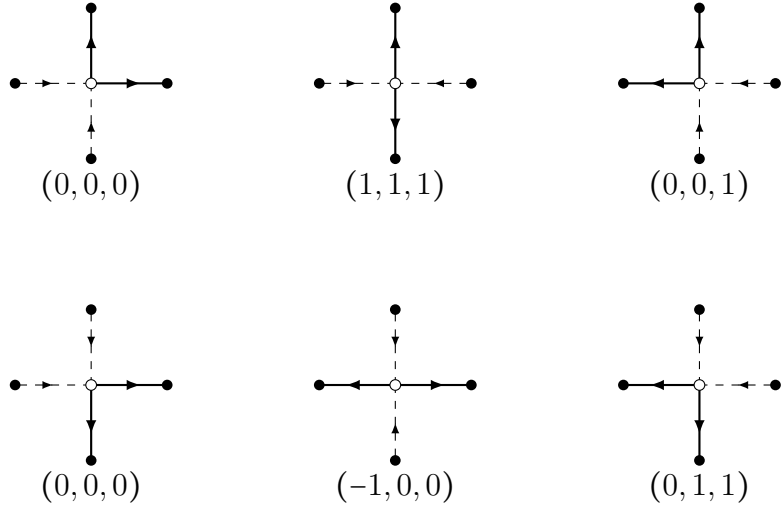
\begin{figure}[htbp]
	\centering
	\begin{tikzpicture}
		\coordinate (base_1) at (0, -1);
			\coordinate (base_2) at (1, 0);
			\coordinate (a) at ($1*(base_1) + 1*(base_2)$);
			\coordinate (b) at ($1*(base_1) + 5*(base_2)$);
			\coordinate (c) at ($1*(base_1) + 9*(base_2)$);
			\coordinate (d) at ($5*(base_1) + 1*(base_2)$);
			\coordinate (e) at ($5*(base_1) + 5*(base_2)$);
			\coordinate (f) at ($5*(base_1) + 9*(base_2)$);
			\foreach \p in {a, b, c, d, e, f}{
				\coordinate (\p_1) at ($(\p) - (base_1)$);
				\coordinate (\p_2) at ($(\p) + (base_2)$);
				\coordinate (\p_3) at ($(\p) + (base_1)$);
				\coordinate (\p_4) at ($(\p) - (base_2)$);
			}
			\foreach \x in {1, 2}{
				\draw[thick, black, directed] (a) -- (a_\x);
			}
			\foreach \x in {3, 4}{
				\draw[dashed, reverse directed] (a) -- (a_\x);
			}
			\node at (a_3)[below]{$(0, 0, 0)$};
			\foreach \x in {1, 3}{
				\draw[thick, black, directed] (b) -- (b_\x);
			}
			\foreach \x in {2, 4}{
				\draw[dashed, reverse directed] (b) -- (b_\x);
			}
			\node at (b_3)[below]{$(1, 1, 1)$};
			\foreach \x in {1, 4}{
				\draw[thick, black, directed] (c) -- (c_\x);
			}
			\foreach \x in {2, 3}{
				\draw[dashed, reverse directed] (c) -- (c_\x);
			}
			\node at (c_3)[below]{$(0, 0, 1)$};
			\foreach \x in {2, 3}{
				\draw[thick, black, directed] (d) -- (d_\x);
			}
			\foreach \x in {1, 4}{
				\draw[dashed, reverse directed] (d) -- (d_\x);
			}
			\node at (d_3)[below]{$(0, 0, 0)$};
			\foreach \x in {2, 4}{
				\draw[thick, black, directed] (e) -- (e_\x);
			}
			\foreach \x in {1, 3}{
				\draw[dashed, reverse directed] (e) -- (e_\x);
			}
			\node at (e_3)[below]{$(-1, 0, 0)$};
			\foreach \x in {3, 4}{
				\draw[thick, black, directed] (f) -- (f_\x);
			}
			\foreach \x in {1, 2}{
				\draw[dashed, reverse directed] (f) -- (f_\x);
			}
			\node at (f_3)[below]{$(0, 1, 1)$};
			\foreach \x in {a, b, c, d, e, f}{
				\fill[fill=white, draw] (\x) circle (2pt);
				\foreach \y in {1, 2, 3, 4}{
					\fill (\x_\y) circle (2pt);
				}
			}
	\end{tikzpicture}
	\caption{The triplet $({a}_{i, j}, {c}_{i, j}, {r}_{i, j})$ which correspond the choise of $(i, j)$}
\end{figure}

Next, we shall show a state of six vertex model in $\mathcal{SV}(n)$ can be determined 
by using the triplet of each entries
for any ASM.
Let us take an ASM ${ \left( {a }_{i, j} \right) }_{1 \leq i, j \leq n}$ arbitrarily.
In this manner,
the direction of a horizontal (resp. vertical) edage 
$\left\{ (i, j), (i, j + 1) \right\}$ (resp. $\left\{ (i, j), (i + 1, j) \right\}$)
can be determined in two ways from 
$\left( {a}_{i, j}, {c}_{i, j}, {r}_{i, j} \right)$ and $\left( {a}_{i, j + 1}, {c}_{i, j + 1}, {r}_{i, j + 1} \right)$
(resp.
$\left( {a}_{i, j}, {c}_{i, j}, {r}_{i, j} \right)$ and $\left( {a}_{i + 1, j}, {c}_{i + 1, j}, {r}_{i + 1, j} \right)$
).
First,
we shall show the two way determine the same direction.
For the edges
$W =\left\{ (i, j), (i, j - 1) \right\}$, $E = \left\{ (i, j), (i, j + 1) \right\}$,
$N =\left\{ (i, j), (i - 1, j) \right\}$ and $S = \left\{ (i, j), (i + 1, j) \right\}$,
the direction of them
are determined by
${r}_{i, j} - {a}_{i, j}$, ${r}_{i, j}$, ${c}_{i, j} - {a}_{i, j}$ and ${c}_{i, j}$ respectively.
We write the four value for $6$ possible choice 
down in the table \ref{table:direction}.
The horizontal edge 
$\left\{ (i, j), (i, j - 1) \right\}$ comes in $(i, j)$ 
if ${r}_{i, j} - {a}_{i, j} = 0$, and
$\left\{ (i, j), (i, j + 1) \right\}$ 
comes in $(i, j + 1)$
if ${r}_{i, j} = 0$.
The vertical edge 
$\left\{ (i, j), (i - 1, j) \right\}$ goes out of $(i, j)$ 
if ${c}_{i, j} - {a}_{i, j} = 0$., and
$\left\{ (i, j), (i + 1, j) \right\}$ 
goes out of $(i + 1, j)$
if ${c}_{i, j} = 0$.
\begin{table}[htbp]
	\centering
	\begin{tabular}{c|cccccc}
			&	$NE$	&	$NS$	&	$NW$	&	$ES$	&	$EW$	&	$SW$	\\
		\hline
		${r}_{i, j} - {a}_{i, j}$	&	$0$	&	$0$	&	$1$	&	$0$	&	$1$	&	$1$	\\
		${r}_{i, j}$	&	$0$	&	$1$	&	$1$	&	$0$	&	$0$	&	$1$	\\
		${c}_{i, j} - {a}_{i, j}$	&	$0$	&	$0$	&	$0$	&	$1$	&	$1$	&	$1$	\\
		${c}_{i, j}$	&	$0$	&	$1$	&	$0$	&	$1$	&	$0$	&	$1$
	\end{tabular}
	\caption{each value for the $6$ possible choices}
	\label{table:direction}
\end{table}
Since ${r}_{i, j}$ (resp. ${c}_{i, j}$) equals ${r}_{i, j + 1} - {a}_{i, j + 1}$ (resp. ${c}_{i + 1, j} - {a}_{i + 1, j}$),
the two ways determine the same.
Then,
the orientation of ${L}_{n}$ is obtaioned by ${ \left( {a}_{i, j} \right) }_{1 \leq i, j \leq n}$,
and it satisfies $2$-in-$2$-out.
Further more,
the orientation satisfies the boundary conditions for
${e}_{2}, {e}_{3}, \ldots , {e}_{n + 1}$
(resp. ${e}_{3n + 2}, {e}_{3n + 3}, \ldots , {e}_{4n}$ and ${e}_{1}$).
because
${r}_{k, 1}$ (resp. ${c}_{1, k}$) equals ${a}_{k, 1}$ (resp. ${a}_{1, k}$)
for $1 \leq k \leq n$.
On the other hand,
the orientation satisfies the boundary conditions for
${e}_{2n + 2}, {e}_{2n + 3}, \ldots , {e}_{3n + 1}$
(resp. ${e}_{n + 2}, {e}_{n + 3}, \ldots , {e}_{2n + 1}$)
because 
${r}_{k, n}$ (resp. ${c}_{n, k}$) equals $1$ for $1 \leq k \leq n$.
Therefore,
a state of six vertex model in $\mathcal{SV}(n)$ cane be determined by the triplet
of each entries for any ASM. 
\section{Height function}
Let $m$ and $n$ be positive integer.
In this section, 
we make a state of six vertex model on ${L}_{m, n}$
correspond to a matrix of size $(m + 1) \times (n + 1)$.
Especially,
we introduce a square matrix of size $n + 1$ 
which is called height function of degree $n$,
and
we shall construct a bijection 
between $\mathcal{SV}(n)$ and the set of all height function of degree $n$.
\subsection{Properties of a boundary condition of six vertex model}
As a preparation
to make a state of six vertex model correspond to a matrix,
we present the following lemma.
\begin{lemma}\label{lem:L_m_n}
	Let $m$ and $n$ be positive integers.
	For any state of six vertex model on ${L}_{m, n}$,
	the number of bounday edge which comes in the boundary vertex
	equals $m + n$.
\end{lemma}
Let $\varphi$ be a state of six vertex model on ${L}_{m, n}$.
First, we show the claim when $m = 1$.
\begin{enumerate}[(I)]
	\item	When $n = 1$,
		the claim is crealy correct from the definition.
	\item	When $n > 2$, we assume the claim holds up to $n - 1$.
				If the edge $\left\{ (1, n - 1), (1, n) \right\}$ goes out of (resp. comes in) $(1, n - 1)$,
				two (resp. one) of the three boundary edges 
				$\left\{ (1, n), (0, n) \right\}$, $\left\{ (1, n), (1, n + 1) \right\}$
				and $\left\{ (1, n), (2, n) \right\}$
				come in the bounday vertex.
				Moreover,
				$n - 1$ (resp. $n$) of the other $2n - 1$ boundary edges come in the boundary vertex
				from the hypothesis of induction.
				Therefore, 
				$n + 1$ of $2n + 2$ boundary edges comes in the boundary vertex.
\end{enumerate}
Then we showed the claim when $m = 1$.
Next,
we assume the claim hold up to $m - 1$ when $m > 2$.
Suppose that exactly $k$ out of $n$ vertical edges 
$\left\{ (m - 1, j), (m, j) \right\}$ ($1 \leq j \leq n$) come in the below vertex
in $\varphi$.
Applying the hypothesis of induction to
the state which is obtained by restricting $\varphi$ to
the $(m - 1) \times n$ grid which has $(1, 1)$ as the mostleft interior veretx in the top row,
it follows that
exactly $m + n - k - 1$ out of the other $2m + n - 2$ boundary edges come in the boundary vertex.
On the other hand,
the $1 \times n$ grid which has $(m, 1)$ as the mostleft interior vertex implies that
exactly $k + 1$ out of $n + 2$ boundary edges
$\left\{ (m, 1), (m, 0) \right\}$, $\left\{ (m, n), (m, n + 1) \right\}$
and $\left\{ (m, j), (m + 1, j) \right\}$ ($1 \leq j \leq n$)
come in the boundary vertex.
Therefore 
exactly $m + n$ boundary edges come in the boundary vertex
in $\varphi$.
Then we showed lemma~\ref{lem:L_m_n}.
\qed

Now we have following claim as a corollary of lemma~\ref{lem:L_m_n}.
\begin{cor}\label{cor:L_m_n}
	Let $n$ be a positive integer.
	For any state of sixvertex model in $\mathcal{SV}(n)$
	and any $1 \leq i \leq n$,
	exactly $i$ of $n$ vertical edges 
	$\left\{ (i, j), (i + 1, j) \right\}$ 
	($1 \leq j \leq n$)
	comes in the below vertex,
	and $n - i$ of $n$ horizontal edges
	$\left\{ (j, i), (j, i + 1) \right\}$
	($1 \leq j \leq n$)
	come in the right vertex.
\end{cor}
\subsection{The map from six vertex model to matricies}
Let us take a state $\varphi$
of six vertex model on ${L}_{m, n}$ arbitrary.
We shall make $\varphi$
of six vertex model on ${L}_{m, n}$
to correspond a $(m + 1) \times (n + 1)$ matrix
$
	{ \left( {h}_{i, j} \right) }_{
			0 \leq i \leq m, \, 
			0 \leq j \leq n
	}
$
which satisfies following conditions:
\begin{subequations}
\begin{align}
	&
	\left| {h}_{i, j} - {h}_{i, j - 1} \right| = 1
	&
	\left( 
		0 \leq i \leq m, \, 
		0 < j \leq n 
	\right)
	, \\
	&
	\left| {h}_{i, j} - {h}_{i - 1, j} \right| = 1
	&
	\left( 
		0 < i \leq m, \, 
		0 \leq j \leq n 
	\right)
	, \\
	&
	{h}_{0, 0} = 0
	.
\end{align}
\end{subequations}
\begin{definition}
Let $m$ and $n$ be positive integers.
We define the matrix ${ \left( {h}_{i, j} \right) }_{0 \leq i \leq m, \, 0 \leq j \leq n}$
by setting 
${h}_{i, j} - {h}_{i, j - 1}$ ($0 \leq i \leq m$, $0 < j \leq n$)
and
${h}_{i, j} - {h}_{i - 1, j}$ ($0 < i \leq m$, $0 \leq j \leq n$)
as follows:
\begin{subequations}
\begin{align}
	{h}_{i, j} - {h}_{i, j - 1}
	&=
	\begin{cases}
		1
		&
		\textrm{if } \left\{ (i, j), (i + 1, j) \right\} \textrm{ comes in } (i, j)
		, \\
		-1
		&
		\textrm{otherwise},
	\end{cases}
	\\
	{h}_{i, j} - {h}_{i - 1, j}
	&=
	\begin{cases}
		1
		&
		\textrm{if } \left\{ (i, j), (i, j  + 1) \right\} \textrm{ goes out of } (i, j)
		, \\
		-1
		&
		\textrm{otherwise}.
	\end{cases}
\end{align}
\end{subequations}
\end{definition}
Then we shall show the well-definedness in this way.
Now,
we set ${k}_{0}$, ${k}_{1}$, ${l}_{0}$ and ${l}_{1}$ as follows respectively:
\begin{subequations}
\begin{align}
	{k}_{0}
	&=
	\# \left\{ i \in [m] \, \middle| \, \{ (i, 0), (i, 1) \} \textrm{ comes in } (i, 0) \right\}
	, \\
	{k}_{1}
	&=
	\# \left\{ i \in [m] \, \middle| \, \{ (i, n), (i, n + 1) \} \textrm{ comes in } (i, n + 1) \right\}
	, \\
	{l}_{0}
	&=
	\# \left\{ j \in [n] \, \middle| \, \{ (0, j), (1, j) \} \textrm{ comes in } (0, j) \right\}
	, \\
	{l}_{1}
	&=
	\# \left\{ j \in [n] \, \middle| \, \{ (m, j), (m + 1, j) \} \textrm{ comes in } (m + 1, j) \right\}
	.
\end{align}
\end{subequations}
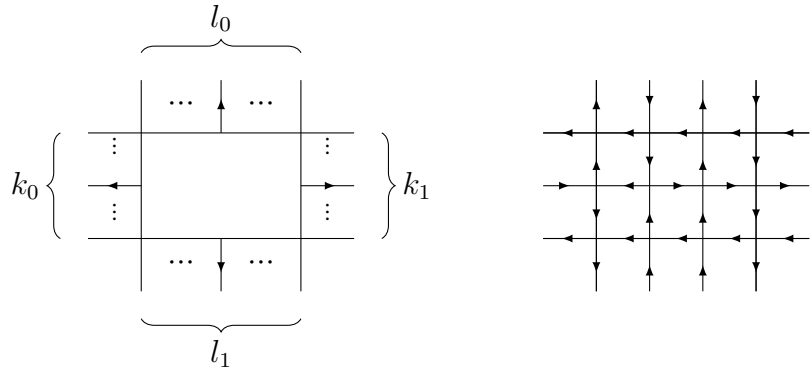
\begin{figure}[htbp]
	\centering
	\begin{tabular}{cc}
	\begin{minipage}{0.45\hsize}
	\centering
	\begin{tikzpicture}[scale=0.7]
		\coordinate (base_1) at (0, -1);
		\coordinate (base_2) at (1, 0);
		\draw ($(base_1)$) -- ($(base_1) + 5*(base_2)$);
		\draw ($3*(base_1)$) -- ($3*(base_1) + 5*(base_2)$);
		\draw ($(base_2)$) -- ($4*(base_1) + (base_2)$);
		\draw ($4*(base_2)$) -- ($4*(base_1) + 4*(base_2)$);
		\foreach \j in {2.5}{
			\draw[directed] ($(base_1) + \j*(base_2)$) -- ($\j*(base_2)$);
		}
		\foreach \j in {1.75, 3.25}{
			\node at ($0.5*(base_1) + \j*(base_2)$){$\cdots$};
		}
		\foreach \i in {2}{
			\draw[directed] ($\i*(base_1) + 4*(base_2)$) -- ($\i*(base_1) + 5*(base_2)$);
		}
		\foreach \i in {1.25, 2.5}{
			\node at ($\i*(base_1) + 4.5*(base_2)$){\scriptsize$\vdots$};
		}
		\foreach \j in {2.5}{
			\draw[directed] ($3*(base_1) + \j*(base_2)$) -- ($4*(base_1) + \j*(base_2)$);
		}
		\foreach \j in {1.75, 3.25}{
			\node at ($3.5*(base_1) + \j*(base_2)$){$\cdots$};
		}
		\foreach \i in {2}{
			\draw[directed] ($\i*(base_1) + 1*(base_2)$) -- ($\i*(base_1)$);
		}
		\foreach \i in {1.25, 2.5}{
			\node at ($\i*(base_1) + 0.5*(base_2)$){\scriptsize$\vdots$};
		}
		\draw [
			decorate,
			decoration={
				brace,
				amplitude=5pt,
				raise=2ex
			}
		]
			($(base_2)$) -- ($4*(base_2)$) 
			node[midway,yshift=2em]{${l}_{0}$};
		\draw [decorate,decoration={brace, amplitude=5pt, mirror, raise=2ex}]
			($4*(base_1) + (base_2)$) -- ($4*(base_1) + 4*(base_2)$) 
			node[midway,yshift=-2em]{${l}_{1}$};
		\draw [decorate,decoration={brace, amplitude=5pt, mirror, raise=2ex}]
			($(base_1)$) -- ($3*(base_1)$) 
			node[midway,xshift=-2em]{${k}_{0}$};
		\draw [decorate,decoration={brace, amplitude=5pt, raise=2ex}]
			($(base_1) + 5*(base_2)$) -- ($3*(base_1) + 5*(base_2)$) 
			node[midway,xshift=2em]{${k}_{1}$};
	\end{tikzpicture}
	\end{minipage}
	&
	\begin{minipage}{0.45\hsize}
	\centering
	\begin{tikzpicture}[scale=0.7]
		\coordinate (base_1) at (0, -1);
		\coordinate (base_2) at (1, 0);
		\draw ($(base_1)$) -- ($(base_1) + 5*(base_2)$);
		\draw ($3*(base_1)$) -- ($3*(base_1) + 5*(base_2)$);
		\draw ($(base_2)$) -- ($4*(base_1) + (base_2)$);
		\draw ($4*(base_2)$) -- ($4*(base_1) + 4*(base_2)$);
		\foreach \j in {1, 3}{
			\draw[directed] ($(base_1) + \j*(base_2)$) -- ($\j*(base_2)$);
		}
		\foreach \j in {2, 4}{
			\draw[reverse directed] ($(base_1) + \j*(base_2)$) -- ($\j*(base_2)$);
		}
		\foreach \i in {2}{
			\draw[directed] ($\i*(base_1) + 4*(base_2)$) -- ($\i*(base_1) + 5*(base_2)$);
		}
		\foreach \i in {1, 3}{
			\draw[reverse directed] ($\i*(base_1) + 4*(base_2)$) -- ($\i*(base_1) + 5*(base_2)$);
		}
		\foreach \j in {1, 4}{
			\draw[directed] ($3*(base_1) + \j*(base_2)$) -- ($4*(base_1) + \j*(base_2)$);
		}
		\foreach \j in {2, 3}{
			\draw[reverse directed] ($3*(base_1) + \j*(base_2)$) -- ($4*(base_1) + \j*(base_2)$);
		}
		\foreach \i in {1, 3}{
			\draw[directed] ($\i*(base_1) + 1*(base_2)$) -- ($\i*(base_1)$);
		}
		\foreach \i in {2}{
			\draw[reverse directed] ($\i*(base_1) + 1*(base_2)$) -- ($\i*(base_1)$);
		}
		\foreach \i/\ii in {1/2, 2/3, 3/4}{
			\draw[reverse directed] ($(base_1) + \i*(base_2)$) -- ($(base_1) + \ii*(base_2)$);
		}
		\foreach \i/\ii in {2/3, 3/4}{
			\draw[directed] ($2*(base_1) + \i*(base_2)$) -- ($2*(base_1) + \ii*(base_2)$);
		}
		\foreach \i/\ii in {1/2}{
			\draw[reverse directed] ($2*(base_1) + \i*(base_2)$) -- ($2*(base_1) + \ii*(base_2)$);
		}
		\foreach \i/\ii in {1/2, 2/3, 3/4}{
			\draw[reverse directed] ($3*(base_1) + \i*(base_2)$) -- ($3*(base_1) + \ii*(base_2)$);
		}
		\foreach \i/\ii in {2/3}{
			\draw[directed] ($\i*(base_1) + 1*(base_2)$) -- ($\ii*(base_1) + 1*(base_2)$);
		}
		\foreach \i/\ii in {1/2}{
			\draw[reverse directed] ($\i*(base_1) + 1*(base_2)$) -- ($\ii*(base_1) + 1*(base_2)$);
		}
		\foreach \i/\ii in {1/2}{
			\draw[directed] ($\i*(base_1) + 2*(base_2)$) -- ($\ii*(base_1) + 2*(base_2)$);
		}
		\foreach \i/\ii in {2/3}{
			\draw[reverse directed] ($\i*(base_1) + 2*(base_2)$) -- ($\ii*(base_1) + 2*(base_2)$);
		}
		\foreach \i/\ii in {1/2, 2/3}{
			\draw[reverse directed] ($\i*(base_1) + 3*(base_2)$) -- ($\ii*(base_1) + 3*(base_2)$);
		}
		\foreach \i/\ii in {1/2, 2/3}{
			\draw[directed] ($\i*(base_1) + 4*(base_2)$) -- ($\ii*(base_1) + 4*(base_2)$);
		}
	\end{tikzpicture}
	\end{minipage}
	\end{tabular}
	\caption{right;${k}_{0} = 2$, ${k}_{1} = 1$, ${l}_{0} = 2$, ${l}_{1} = 2$.}
\end{figure}
First, we show that
the values of 
${h}_{m, n}$ 
are consistent
when it is determined clockwise 
by ${l}_{0}$ and ${k}_{1}$
and 
when it is determined counterclockwise
by ${k}_{0}$ and ${l}_{1}$.
If we determine ${h}_{m, n}$ clockwise,
then
we have
\begin{align}\label{eq:height_CW}
	{h}_{m, n}
	&=
	2\left( {k}_{1} + {l}_{0} \right) - \left( m + n \right)
	.
\end{align}
On the other hand,
\begin{align}\label{eq:height_CCW}
	{h}_{m, n}
	&=
	-2 \left( {k}_{0} + {l}_{1} \right) + \left( m + n \right)
\end{align}
if we determine counterclockwise.
Since
$
	{k}_{0} + {k}_{1} + {l}_{0} + {l}_{1} 
	= 
	m + n
$,
from the lemma~\ref{lem:L_m_n}, 
the difference between
the RHS of \eqref{eq:height_CW} 
and 
the RHS of \eqref{eq:height_CCW}
equals $0$.
Then we showed 
the values of ${h}_{m, n}$
which is determined clockwise 
and 
which is determined counterclockwise
are consistent.
Moreover,
subsutitute
$m + n = {k}_{0} + {k}_{1} + {l}_{0} + {l}_{1}$
for \eqref{eq:height_CW} yields
\begin{align}\label{eq:height_4value}
	{h}_{m, n}
	&=
	{l}_{0} + {k}_{1} - {l}_{1} - {k}_{0}
	.
\end{align}

Next, let us take 
$1 \leq {i}_{0} \leq m$ 
and 
$1 \leq {j}_{0} \leq n$
arbitrarily.
We shall show that
the value of ${h}_{ {i}_{0}, {j}_{0} }$,
which is determined from the ${i}_{0} \times {j}_{0}$ grid
with $(1, 1)$ as the leftmost interior vertex of the top row,
and 
the value of ${h}_{m, n} - {h}_{ {i}_{0}, {j}_{0} }$,
which is determined from the $\left( m - {i}_{0} \right) \times \left( n - {j}_{0} \right)$ grid
with $(m, n)$ as the rightmost interior vertex of the bottom row,
does not contradict each other.
Now,
we set 
${k}_{i}$, ${l}_{i}$ ($2 \leq i \leq 5$) 
as follows respectively:
\begin{subequations}
\begin{align}
	{k}_{2}
	&=
	\# \left\{ i \in [{i}_{0}] \, \middle| \, \{ (i, 0), (i, 1) \} \textrm{ comes in } (i, 0) \right\}
	,\\
	{k}_{3}
	&=
	\# \left\{ i \in [{i}_{0}] \, \middle| \, \{ (i, {j}_{0}), (i, {j}_{0} + 1) \} \textrm{ comes in } (i, {j}_{0} + 1) \right\}
	, \\
	{k}_{4}
	&=
	\# \left\{ {i}_{0} < i \leq m \, \middle| \, \{ (i, {j}_{0}), (i, {j}_{0} + 1) \} \textrm{ comes in } (i, {j}_{0}) \right\}
	, \\
	{k}_{5}
	&=
	\# \left\{ {i}_{0} < i \leq m \, \middle| \, \{ (i, n), (i, n + 1) \} \textrm{ comes in } (i, n + 1) \right\}
	, \\
	{l}_{2}
	&=
	\# \left\{ j \in [{j}_{0}] \, \middle| \, \{ (0, j), (1, j) \} \textrm{ comes in } (0, j) \right\}
	, \\
	{l}_{3}
	&=
	\# \left\{ j \in [{j}_{0}] \, \middle| \, \{ ({i}_{0}, j), ({i}_{0} + 1, j) \} \textrm{ comes in } ({i}_{0} + 1, j) \right\}
	, \\
	{l}_{4}
	&=
	\# \left\{ {j}_{0} < j \leq n \, \middle| \, \{ ({i}_{0}, j), ({i}_{0} + 1, j) \} \textrm{ comes in } ({i}_{0}, j) \right\}
	, \\
	{l}_{5}
	&=
	\# \left\{ {j}_{0} < j \leq n \, \middle| \, \{ (m, j), (m + 1, j) \} \textrm{ comes in } (m + 1, j) \right\}
	.
\end{align}
\end{subequations}
\begin{figure}[htbp]
	\centering
	\begin{tikzpicture}[scale=0.5]
		\coordinate (base_1) at (0, -1);
		\coordinate (base_2) at (1, 0);
		\draw ($(base_1)$) -- ($(base_1) + 13*(base_2)$);
		\draw[very thick] ($6*(base_1)$) -- ($6*(base_1) + 13*(base_2)$);
		\draw ($7*(base_1)$) -- ($7*(base_1) + 13*(base_2)$);
		\draw ($11*(base_1)$) -- ($11*(base_1) + 13*(base_2)$);
		\draw ($(base_2)$) -- ($12*(base_1) + (base_2)$);
		\draw[very thick] ($7*(base_2)$) -- ($12*(base_1) + 7*(base_2)$);
		\draw ($8*(base_2)$) -- ($12*(base_1) + 8*(base_2)$);
		\draw ($12*(base_2)$) -- ($12*(base_1) + 12*(base_2)$);
		\foreach \j in {4}{
			\draw[directed] ($(base_1) + \j*(base_2)$) -- ($\j*(base_2)$);
		}
		\foreach \j in {2.5, 5.5}{
			\node at ($0.5*(base_1) + \j*(base_2)$){$\cdots$};
		}
		\foreach \i in {3.5}{
			\draw[directed] ($\i*(base_1) + 7*(base_2)$) -- ($\i*(base_1) + 8*(base_2)$);
		}
		\foreach \i in {2, 4.5}{
		\node at ($\i*(base_1) + 7.5*(base_2)$){$\vdots$};
		}
		\foreach \j in {4}{
			\draw[directed] ($6*(base_1) + \j*(base_2)$) -- ($7*(base_1) + \j*(base_2)$);
		}
		\foreach \j in {2.5, 5.5}{
			\node at ($6.5*(base_1) + \j*(base_2)$){$\cdots$};
		}
		\foreach \i in {3.5}{
			\draw[directed] ($\i*(base_1) + 1*(base_2)$) -- ($\i*(base_1)$);
		}
		\foreach \i in {2, 4.5}{
			\node at ($\i*(base_1) + 0.5*(base_2)$){$\vdots$};
		}
		\foreach \j in {10}{
			\draw[directed] ($7*(base_1) + \j*(base_2)$) -- ($6*(base_1) + \j*(base_2)$);
		}
		\foreach \j in {9, 11}{
			\node at ($6.5*(base_1) + \j*(base_2)$){$\cdots$};
		}
		\foreach \i in {9}{
			\draw[directed] ($\i*(base_1) + 12*(base_2)$) -- ($\i*(base_1) + 13*(base_2)$);
		}
		\foreach \i in {8, 10}{
			\node at ($\i*(base_1) + 12.5*(base_2)$){$\vdots$};
		}
		\foreach \j in {10}{
			\draw[directed] ($11*(base_1) + \j*(base_2)$) -- ($12*(base_1) + \j*(base_2)$);
		}
		\foreach \j in {9, 11}{
			\node at ($11.5*(base_1) + \j*(base_2)$){$\cdots$};
		}
		\foreach \i in {9}{
			\draw[directed] ($\i*(base_1) + 8*(base_2)$) -- ($\i*(base_1) + 7*(base_2)$);
		}
		\foreach \i in {8, 10}{
			\node at ($\i*(base_1) + 7.5*(base_2)$){$\vdots$};
		}
		\draw [decorate, decoration={brace, amplitude=5pt, raise=2ex}]
			($(base_2)$) -- ($7*(base_2)$) 
			node[midway,yshift=2em]{${l}_{2}$};
		\draw [decorate,decoration={brace, amplitude=5pt, raise=2ex}]
			($6*(base_1) + (base_2)$) -- ($6*(base_1) + 7*(base_2)$) 
			node[midway,yshift=2em]{${l}_{3}$};
		\draw [decorate, decoration={brace, amplitude=5pt, mirror, raise=2ex}]
			($7*(base_1) + 8*(base_2)$) -- ($7*(base_1) + 12*(base_2)$) 
			node[midway,yshift=-2em]{${l}_{4}$};
		\draw [decorate, decoration={brace, amplitude=5pt, mirror, raise=2ex}]
			($12*(base_1) + 8*(base_2)$) -- ($12*(base_1) + 12*(base_2)$) 
			node[midway,yshift=-2em]{${l}_{5}$};
		\draw [decorate,decoration={brace, amplitude=5pt, mirror, raise=2ex}]
			($(base_1)$) -- ($6*(base_1)$) 
			node[midway,xshift=-2em]{${k}_{2}$};
		\draw [decorate,decoration={brace, amplitude=5pt, raise=2ex}]
			($(base_1) + 8*(base_2)$) -- ($6*(base_1) + 8*(base_2)$) 
			node[midway,xshift=2em]{${k}_{3}$};
		\draw [decorate,decoration={brace, amplitude=5pt, mirror, raise=2ex}]
			($7*(base_1) + 7*(base_2)$) -- ($11*(base_1) + 7*(base_2)$) 
			node[midway,xshift=-2em]{${k}_{4}$};
		\draw [decorate,decoration={brace, amplitude=5pt, raise=2ex}]
			($7*(base_1) + 13*(base_2)$) -- ($11*(base_1) + 13*(base_2)$) 
			node[midway,xshift=2em]{${k}_{5}$};
		\node at ($6*(base_1) + 13*(base_2)$)[right]{$i = {i}_{0}$};
		\node at ($12*(base_1) + 7*(base_2)$)[below]{$j = {j}_{0}$};
	\end{tikzpicture}
	\caption{}
\end{figure}
Focusing on
the ${i}_{0} \times {j}_{0}$ grid, 
we have
\begin{align}\label{eq:height_above_left_grid}
	{h}_{ {i}_{0}, {j}_{0} }
	&=
	{l}_{2} + {k}_{3} - {l}_{3} - {k}_{2}
	.
\end{align}
On the other hand,
it follows from
the $\left( m - {i}_{0} \right) \times \left(n - {j}_{0} \right)$ grid 
that
\begin{align}\label{eq:height_below_right_grid}
	{h}_{m, n} - {h}_{ {i}_{0}, {j}_{0} }
	&=
	{l}_{4} + {k}_{5} - {l}_{5} - {k}_{4}
	.
\end{align}
Then
we just show 
the sum of
the RHS of \eqref{eq:height_above_left_grid}
and 
the RHS of \eqref{eq:height_below_right_grid}
equals 
the RHS of \eqref{eq:height_4value}.
Now, 
from the lemma~\ref{lem:L_m_n},
the boundary condition of 
the ${i}_{0} \times \left( n - {j}_{0} \right)$ grid 
with $(1, n)$ as the rightmost interior vertex of the top row
implies that
\begin{align}\label{eq:height_above_right_grid}
	{l}_{0} - {l}_{2} 
	+ {k}_{1} - {k}_{5}  
	- {l}_{4} 
	- {k}_{3} 
	&=
	0
	.
\end{align}
Moreover,
the boundary condition of
the $\left( m - {i}_{0} \right) \times {j}_{0}$ grid
with $(m, 1)$ as the leftmost interior vertex of the bottom row
show that
\begin{align}\label{eq:height_below_left_grid}
	- {l}_{3} 
	- {k}_{4}  
	+{l}_{1} - {l}_{5} 
	+ {k}_{0} - {k}_{2} 
	&=
	0
	.
\end{align}
We remark that
the sum of
the RHS of \eqref{eq:height_above_left_grid}
and 
the RHS of \eqref{eq:height_below_right_grid}
equals
$
	\left( \textrm{RHS of \eqref{eq:height_4value}} \right)
	- 
	\left( \textrm{LHS of \eqref{eq:height_above_right_grid}} \right)
	+ 
	\left( \textrm{LHS of \eqref{eq:height_below_left_grid}} \right)
$.
It follows that
the sum of
the RHS of \eqref{eq:height_above_left_grid}
and 
the RHS of \eqref{eq:height_below_right_grid}
equals
the RHS of \eqref{eq:height_4value}.
Therefore, 
it was shown that 
the value of ${h}_{ {i}_{0}, {j}_{0} }$,
which is determined from the ${i}_{0} \times {j}_{0}$ grid,
and 
the value of ${h}_{m, n} - {h}_{ {i}_{0}, {j}_{0} }$,
which is determined from the $\left( m - {i}_{0} \right) \times \left( n - {j}_{0} \right)$ grid,
does not contradict each other
for any
$1 \leq {i}_{0} \leq m$, $1 \leq {j}_{0} \leq n$.
\subsection{height function}
We define height function of degree $n$.
\begin{definition}
	The $(n + 1) \times (n + 1)$ matrix 
	$H = {\left( {h}_{i, j} \right)}_{0 \leq i, j \leq n}$ is called 
	\textsl{height function} of degree $n$
	if it satisfying the following conditions: 
	\begin{subequations}
		\begin{align}
			\label{eq:adj_condition_row}
			&
			\left|
				{h}_{i + 1, j} - {h}_{i, j}
			\right|
			= 1
				&	(0 \leq i < n, \; 0 \leq j \leq n), 
			\\
			\label{eq:adj_condition_column}
			&
			\left|
				{h}_{i, j + 1} - {h}_{i, j}
			\right|
			= 1
				&	(0 \leq i \leq n, \; 0 \leq j < n),
			\\
			\label{eq:bc_of_height}
			&
			{h}_{i, 0} = {h}_{0, i} = {h}_{n -i, n} = {h}_{n, n - i} 
				= i
					&	(0 \leq i \leq n).
		\end{align}
	\end{subequations}
\end{definition}
We call \eqref{eq:adj_condition_row} and \eqref{eq:adj_condition_column} adjacent conditons,
and \eqref{eq:bc_of_height} boundary condition of height function.
We denote the set of all height functions of degree $n$ as $\mathcal{H}_{n}$.
A partially order on $\mathcal{H}_{n}$
is defined as
$
	{ \left( {h}_{i, j} \right) }_{0 \leq i, j \leq n}
	\leq
	{ \left( {g}_{i, j} \right) }_{0 \leq i, j \leq n}
$
if
$
	{h}_{i, j} \leq {g}_{i, j}
$
for $1 \leq i, j < n$.
Now,
we hold the following proposition.
\begin{proposition}
	Let $n$ be a positive integer.
	We define a map 
	$
		f 
		\colon
		\mathcal{SV}(n) \rightarrow \mathcal{H}_{n}
		;
		\varphi \mapsto { \left( {h}_{i, j} \right) }_{0 \leq i, j \leq n}
	$
	by setting 
${h}_{i, j} - {h}_{i, j - 1}$ ($0 \leq i \leq n$, $0 < j \leq n$)
and
${h}_{i, j} - {h}_{i - 1, j}$ ($0 < i \leq n$, $0 \leq j \leq n$)
as follows:
\begin{subequations}
\begin{align}
	{h}_{i, j} - {h}_{i, j - 1}
	&=
	\begin{cases}
		1
		&
		\textrm{if } \left\{ (i, j), (i + 1, j) \right\} \textrm{ comes in } (i, j)
		, \\
		-1
		&
		\textrm{otherwise},
	\end{cases}
	\\
	{h}_{i, j} - {h}_{i - 1, j}
	&=
	\begin{cases}
		1
		&
		\textrm{if } \left\{ (i, j), (i, j  + 1) \right\} \textrm{ goes out of } (i, j)
		, \\
		-1
		&
		\textrm{otherwise}.
	\end{cases}
\end{align}
\end{subequations}
Then 
the map
$
	f 
	\colon
	\mathcal{SV}(n) \rightarrow \mathcal{H}_{n}
$
is bijective.
\end{proposition}
\begin{proof}
From the definition,
it is clear that $f$ is injective.
We just show surjectivity.
First,
we remark that
any height function
${ \left( {h}_{i, j} \right) }_{0 \leq i, j \leq n}$
gives a orientation of ${L}_{n}$
by setting the direction of 
$\left\{ (i, j), (i + 1, j) \right\}$
from the value ${h}_{i, j} - {h}_{i, j - 1}$
for $0 \leq i \leq n$, $1 \leq j \leq n$,
and
the direction of 
$\left\{ (i, j), (i, j + 1) \right\}$
from the value ${h}_{i, j} - {h}_{i - 1, j}$
for $1 \leq i \leq n$, $0 \leq j \leq n$.
Moreover,
it is clear that
the boundary condition of the orientation
is the open boundary condition.
Second,
the triplet
$
	\left(
		{h}_{i - 1, j}, {h}_{i, j}, {h}_{i, j - 1}
	\right)
$
can be one of the following six values:
\begin{subequations}
\begin{align}
	\label{eq:height_NE}
	&\left( {h}_{i - 1, j - 1} + 1, {h}_{i - 1, j - 1} + 2, {h}_{i - 1, j - 1} + 1 \right),	\\
	\label{eq:height_NS}
	&\left( {h}_{i - 1, j - 1} + 1, {h}_{i - 1, j - 1}, {h}_{i - 1, j - 1} + 1 \right),	\\
	\label{eq:height_NW}
	&\left( {h}_{i - 1, j - 1} + 1, {h}_{i - 1, j - 1}, {h}_{i - 1, j - 1} - 1 \right),	\\
	\label{eq:height_ES}
	&\left( {h}_{i - 1, j - 1} - 1, {h}_{i - 1, j - 1}, {h}_{i - 1, j - 1} + 1 \right),	\\
	\label{eq:height_EW}
	&\left( {h}_{i - 1, j - 1} - 1, {h}_{i - 1, j - 1}, {h}_{i - 1, j - 1} - 1 \right),	\\
	\label{eq:height_SW}
	&\left( {h}_{i - 1, j - 1} - 1, {h}_{i - 1, j - 1} - 2, {h}_{i - 1, j - 1} - 1 \right).
\end{align}
\end{subequations}
Then
exactly two of the four edges which are incident to $(i, j)$
goes out of $(i, j)$
in the orientation which is determined by a height function.
In fact,
when \eqref{eq:height_NE},
a state of $(i, j)$ is $NE$,
when \eqref{eq:height_NS},
a state of $(i, j)$ is $NS$,
when \eqref{eq:height_NW},
a state of $(i, j)$ is $NW$,
when \eqref{eq:height_ES},
a state of $(i, j)$ is $ES$,
when \eqref{eq:height_EW},
a state of $(i, j)$ is $EW$,
and
when \eqref{eq:height_SW},
a state of $(i, j)$ is $SW$
for any $1 \leq i, j \leq n$.
Therefore,
we obtain a state of six vertex model on ${L}_{n}$
which has the open boundary condition
from each height function of degree $n$.
\end{proof}
\subsection{The bijection between ASM and height function}
Up to here,
we showed a bijection
between $\mathcal{A}_{n}$ and $\mathcal{SV}(n)$,
and
a bijection
between $\mathcal{SV}(n)$ and $\mathcal{H}_{n}$.
Then
we hold the following proposition
from the two bijections.
\begin{proposition}
	Let $n$ be a positive integer.
	For any ASM ${ \left( {a}_{i, j} \right) }_{1 \leq i, j \leq n} \in \mathcal{A}_{n}$,
	we set a matrix ${ \left( {h}_{i, j} \right) }_{0 \leq i, j \leq n}$
	as following:
	\begin{subequations}
	\begin{align}\label{eq:ASM_to_height}
		&
		{h}_{i, j}
		:=
		i + j - 2 \displaystyle\sum_{1 \leq k \leq i} \displaystyle\sum_{1 \leq l \leq j} {a}_{i, j}
		&
		\left( 0 \leq i, j \leq n \right) .
	\end{align}
	Then a bijection between $\mathcal{A}_{n}$ and $\mathcal{H}_{n}$ is given in this way.
	Here, the inverse map is given as following:
	\begin{align}
		&
		{a}_{i, j}
		:=
		- \frac{1}{2}
		\left(
			{h}_{i - 1, j - 1} - {h}_{i - 1, j} + {h}_{i, j} - {h}_{i, j - 1}
		\right)
		&
		\left(
			1 \leq i, j \leq n
		\right) .
	\end{align}
	\end{subequations}
\end{proposition}
\begin{proof}
First, 
we recall that
a vertical edge
$\left\{ (i, j), (i + 1, j) \right\}$ goes out of $(i, j)$
if and only if
$
	{c}_{i, j}
	= 
	\sum_{1 \leq k \leq i} {a}_{i, j}
	= 
	1
$
for $1 \leq i, j \leq n$.
Then
we  have
\begin{align}
	&
	{h}_{i, j}
	=
	i + j 
	- 2 \# \left\{ 
		k \in [j] 
	\, \middle| \, 
		{c}_{i, j}
		= 1 
	\right\}
	&
	\left(
		1 \leq i \leq n,
		0 \leq j \leq n
	\right)
	.
\end{align}
Since
$
	\# \left\{ 
		k \in [j] 
	\, \middle| \, 
		{c}_{i, j}
		= 1 
	\right\}
	=
	\sum_{1 \leq l \leq j} {c}_{i, l}
$,
the equation
\eqref{eq:ASM_to_height} follows.
\end{proof}
\subsection{Properties of height functions}
In this section,
we shall state more properties of height functions.
Especially, we focus on the possible values of each entry of a height function.
First,
we denote $\mathcal{I}(n)$ as the set 
$
	\left\{
		(i, j) \in \mathbb{Z}^{2}
	\, \middle| \,
		1 \leq i, j < n
	\right\}
$.
Now,
we define a positive integer-valued function
$\operatorname{trc}	\colon	\mathcal{I}(n)	\rightarrow	\mathbb{Z}$.
\begin{definition}
	Let $n$ be a positive integer.
	We define 
	$
		\operatorname{trc}
		\colon
		\mathcal{I}(n)
		\rightarrow
		\mathbb{Z}
	$
	as follows:
	\begin{align}
		\operatorname{trc}(i, j)
		&:=
		\min \left\{ i, n - i, j, n - j \right\}
		&
		\left(
			1 \leq i, j < n
		\right)
		.
	\end{align}
	Then 
	we call $\operatorname{trc}(i, j)$ 
	the \textsl{track} of $(i, j)$ for each $(i, j) \in \mathcal{I}(n)$.
\end{definition}
We remark that
$\mathcal{I}(n)$ is 
decomposed 
into $\left\lfloor n/2 \right\rfloor$
disjoint union 
by the track
(i.e., 
$
	\mathcal{I}(n)
	:=
	\bigsqcup_{
		1 \leq l \leq \left\lfloor \frac{n}{2} \right\rfloor
	}
	\left\{
		(i, j) \in \mathcal{I}(n)
	\, \middle| \,
		\operatorname{trc}(i, j) = l
	\right\}
$).

Now,
we hold following propisition.
\begin{proposition}\label{prop:range_of_height}
	Let $n$ be a positive integer, 
	and 
	${ \left( {h}_{i, j} \right) }_{0 \leq i, j \leq n}$
	a height function of degree $n$.
	For $(i, j) \in \mathcal{I}(n)$,
	there are $(l + 1)$s possible values for
	the $(i, j)$-entry ${h}_{i, j}$
	if 
	track of $(i, j)$ equals $l$.
	The possible values are as follows:
	\begin{enumerate}[(i)]
		\item
		when $i = l$ and $l \leq j \leq n - l$
		(resp. $l \leq i \leq n - l$ and $j = l$),
		$
			{h}_{i, j}
			=
			(j - l) + 2k
		$
		(resp.
		$
			(i - l) + 2k
		$
		)
		($0 \leq k \leq l$),
		\item
		when $l \leq i \leq n - l$ and $j = n - l$
		(resp. $i = n - l$ and $l \leq j \leq n - l$),
		$
			{h}_{i, j}
			=
			(n - l - i) + 2k
		$
		(resp.
		$
			(n - l - j) + 2k
		$
		)
		($0 \leq k \leq l$).
	\end{enumerate}
\end{proposition}
\begin{proof}
Now,
we denote $\sum_{1 \leq i \leq i, 1 \leq l \leq j} {a}_{i, j}$ as ${s}_{i, j}$
($0 \leq i, j \leq n$)
for 
${ \left( {a}_{i, j} \right) }_{1 \leq i, j \leq n} \in \mathcal{A}_{n}$.
Here,
we remark that
${s}_{i, 0} = {s}_{0, i} = 0$ for $0 \leq i \leq n$.
Since
the total sum of a row of ASM equals $1$,
for $1 \leq i \leq n$,
we have the following equations:
\begin{subequations}
\begin{align}
	&
	0 \leq {s}_{i, 1} \leq {s}_{i, 2} \leq \cdots \leq {s}_{i, n} = i,
	\\
	&
	{s}_{i, j} - {s}_{i, j - 1}
	=
	0 \textrm{ or } 1
	&
	\left(
		1 \leq j \leq n
	\right)
	.
\end{align}
\end{subequations}
Then it follows that
\begin{align}\label{eq:cornersum_conditon}
	&
	\max \left\{ 0, i + j - n \right\} \leq 
	{s}_{i, j}
	\leq \min \left\{ i, j \right\}
	&
	\left(
		1 \leq i, j \leq n
	\right)
	.
\end{align}
In the same way for a column, 
exactly the same equation \eqref{eq:cornersum_conditon} holds.
Moreover,
the equation \eqref{eq:ASM_to_height} show that
\begin{enumerate}[(i)]
	\item\label{case:LL}
	when $i + j \leq n$ and $i \leq j$,
	$
		j - i \leq
		{h}_{i, j}
		\leq i + j
	$,
	\item\label{case:GL}
	when $i + j \geq n$ and $i \leq j$,
	$
		j - i \leq
		{h}_{i, j}
		\leq 2n - (i + j)
	$,
	\item\label{case:GG}
	when $i + j \geq n$ and $i \geq j$,
	$
		i - j \leq
		{h}_{i, j}
		\leq 2n - (i + j)
	$,
	\item\label{case:LG}
	when $i + j \leq n$ and $i \geq j$,
	$
		i - j \leq
		{h}_{i, j}
		\leq i + j
	$.
\end{enumerate}
We remark that
the track of $(i, j)$ equals
$i$ when \eqref{case:LL},
$(n - j)$ when \eqref{case:GL},
$(n - i)$ when \eqref{case:GG},
$j$ when \eqref{case:LG}.
It follows that
the difference between the RHS and the LHS
is equal to twice the track of $(i, j)$
for each case.
Since ${h}_{i, j}$ decrease $2$
every time ${s}_{i, j}$ increase $1$,
proposition \ref{prop:range_of_height} follows.
\end{proof}
Now,
we define $\mathbb{P}_{n}$ as the set
	$
		\left\{
			(i, j, k) \in \mathbb{Z}^{3}
		\, \middle| \,
			1 \leq i, j < n \textrm{ and }
			1 \leq k \leq i, j, n - i, n - j
		\right\}
	$. 
Let us set
$
	\iota \colon \mathcal{H}_{n} \rightarrow {2}^{ \mathbb{P}_{n} }
$
as a map such that
$(i, j, k) \in \iota \left( { \left( {h}_{i, j} \right) }_{0 \leq i, j \leq n} \right)$
if and only if
${h}_{i, j} \geq \left| i - j \right| + 2k$.
Then,
we have 
$
	\iota \left( { \left( {h}_{i, j} \right) }_{0 \leq i, j \leq n} \right)
	\subset
	\iota \left( { \left( {g}_{i, j} \right) }_{0 \leq i, j \leq n} \right)
$
if
$
	{ \left( {h}_{i, j} \right) }_{0 \leq i, j \leq n}
	\leq
	{ \left( {g}_{i, j} \right) }_{0 \leq i, j \leq n}
$.
Here,
we define a partially order on $\mathbb{P}_{n}$ by a cover relation.
The element
$(i, j, k)$ covers $(i', j', k')$ if it satisfies the following conditions:
\begin{subequations}
\begin{align}
	&
	\left| i - j \right| + 2k = \left| i' - j' \right| + 2k' + 1,
	\\
	&
	\left| i - i' \right| + \left| j - j' \right| = 1.
\end{align}
\end{subequations}
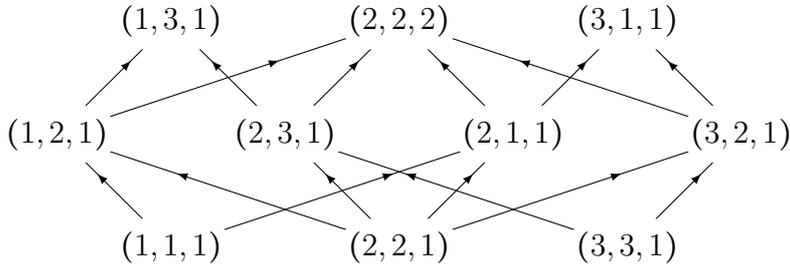
\begin{figure}[htbp]
	\centering
	\begin{tikzpicture}[xscale=3.0, yscale=1.5]
		\foreach \i/\j/\k in {1/1/1, 2/2/2, 3/3/3}{
			\coordinate (h_\i_\j_2) at ($(\k, 0) - (2, 0)$);
		}
		\foreach \i/\j/\k in {1/2/1, 2/3/2, 2/1/3, 3/2/4}{
			\coordinate (h_\i_\j_3) at ($(\k, 1) - (2.5, 0)$); 
		}
		\foreach \i/\j/\k in {1/3/1, 2/2/2, 3/1/3}{
			\coordinate (h_\i_\j_4) at ($(\k, 2) - (2, 0)$);
		}
		\foreach \j/\r/\jj/\rr in {1/2/2/3, 2/3/3/4}{
			\draw[directed] (h_1_\j_\r) -- (h_1_\jj_\rr);
		}
		\foreach \j/\r/\jj/\rr in {2/2/3/3}{
			\draw[directed] (h_2_\j_\r) -- (h_2_\jj_\rr);
		}
		\foreach \i/\r/\ii/\rr in {3/2/2/3, 2/3/1/4}{
			\draw[directed] (h_\i_3_\r) -- (h_\ii_3_\rr);
		}
		\foreach \i/\r/\ii/\rr in {2/2/1/3}{
			\draw[directed] (h_\i_2_\r) -- (h_\ii_2_\rr);
		}
		\foreach \j/\r/\jj/\rr in {1/3/2/4}{
			\draw[directed] (h_2_\j_\r) -- (h_2_\jj_\rr);
		}
		\foreach \i/\r/\ii/\rr in {3/3/2/4}{
			\draw[directed] (h_\i_2_\r) -- (h_\ii_2_\rr);
		}
		\foreach \j/\r/\jj/\rr in {3/2/2/3, 2/3/1/4}{
			\draw[directed] (h_3_\j_\r) -- (h_3_\jj_\rr);
		}
		\foreach \j/\r/\jj/\rr in {2/2/1/3}{
			\draw[directed] (h_2_\j_\r) -- (h_2_\jj_\rr);
		}
		\foreach \i/\r/\ii/\rr in {1/2/2/3, 2/3/3/4}{
			\draw[directed] (h_\i_1_\r) -- (h_\ii_1_\rr);
		}
		\foreach \i/\r/\ii/\rr in {2/2/3/3}{
			\draw[directed] (h_\i_2_\r) -- (h_\ii_2_\rr);
		}
		\foreach \j/\r/\jj/\rr in {3/3/2/4}{
			\draw[directed] (h_2_\j_\r) -- (h_2_\jj_\rr);
		}
		\foreach \i/\r/\ii/\rr in {1/3/2/4}{
			\draw[directed] (h_\i_2_\r) -- (h_\ii_2_\rr);
		}
		\foreach \j/\r in {1/2, 2/3}{
			\fill[fill=white] (h_1_\j_\r) circle (8pt);
			\node at (h_1_\j_\r){
				$(1, \j, 1)$
			};
		}
		\foreach \i/\r in {1/4, 2/3}{
			\fill[fill=white] (h_\i_3_\r) circle (8pt);
			\node at (h_\i_3_\r){
				$(\i, 3, 1)$
			};
		}
		\foreach \j/\r in {3/2, 2/3}{
			\fill[fill=white] (h_3_\j_\r) circle (8pt);
			\node at (h_3_\j_\r){
				$(3, \j, 1)$
			};
		}
		\foreach \i/\r in {3/4, 2/3}{
			\fill[fill=white] (h_\i_1_\r) circle (8pt);
			\node at (h_\i_1_\r){
				$(\i, 1, 1)$
			};
		}
		\fill[fill=white] (h_2_2_2) circle (8pt);
		\node at (h_2_2_2){
			$(2, 2, 1)$
		};
		\fill[fill=white] (h_2_2_4) circle (8pt);
		\node at (h_2_2_4){
			$(2, 2, 2)$
		};
	\end{tikzpicture}
	\caption{The Hasse diagram of $\mathbb{P}_{3}$.}
\end{figure}
As a observation,
$\iota \colon \mathcal{H}_{n} \rightarrow {2}^{ \mathbb{P}_{n} }$
gives a order ideal of $\mathbb{P}_{n}$.
Moreover,
the poset $\mathbb{P}_{n}$ is graded, 
the rank of $(i, j, k) \in \mathbb{P}_{n}$ is expressed as
$\left| i - j \right| + 2k - 2$,
and the generating funcition is
$
	\sum_{0 \leq r \leq n - 2}
		(n - r - 1)(r + 1) {q}^{r}
$.
\section{Perfect matching and link pattern}
Let $n$ be a positive integer,
and ${L}_{n}$ the graph which is introduced in section \ref{subsection:SVM_FPL}.
When a FPL on ${L}_{n}$ is given, 
we can obtain a pairing of boundary vertices
if we focus only on edges with either black or white color.
In this sectetion,
we define monochromatic path 
and perfect matching
as a preparation for describing
the behavior of FPLs using that pairing.
Let 
$
	\boldsymbol{p} = 
	\left(
		{v}_{0}, {v}_{1}, \ldots , {v}_{m}
	\right)
$ 
be a path in ${L}_{n}$ 
where ${v}_{0}, {v}_{1} , \ldots , {v}_{m}$ are all distinct,
and $\psi$ a FPL. . 
We call $\boldsymbol{p}$ a
\textsl{monochromatic path} of $\psi$ if 
all edges have the same color
(i.e., 
$
	\psi(\{ {v}_{i - 1}, {v}_{i} \} )
	=	\psi(\{ {v}_{0}, {v}_{1} \} )
$
for $1 \leq i \leq m$). 
Further, 
we call $\boldsymbol{p}$ a \textsl{black path} (resp. \textsl{white path})
or we say $\boldsymbol{p}$ has color $b$ (resp. $w$) if
all edges have color $b$
(resp.
$w$). 
In graph theory, a path $\boldsymbol{p}$ is called \textsl{cycle} 
if the two end vertices ${v}_{0}$ and ${v}_{m}$ coincide. 
When a monochromatic path $\boldsymbol{p}$ is a cycle, 
we call $\boldsymbol{p}$ a
\textsl{monochromatic cycle}. 
Next, we define matching on $[n]$.
Here, we remark that $[n]$ is the set 
$\{ 1, 2, \ldots , n \}$. 
\begin{definition}
Let $n, p$ be positive integers which satisfy 
$
	0 \leq 2p \leq n
$, 
and  
$
	\mu = 
		\left(
			[n], E(\mu)
		\right)
$
be a graph. 
We say $\mu$ is a
\textsl{$p$-matching on $[n]$} if it satisfies the following conditions:
	\begin{enumerate}[(i)]
		\item
			$
					{u}_{1}, {u}_{2}, \ldots , {u}_{p}, {v}_{1}, {v}_{2}, \ldots , {v}_{p}
			$	
			are elements of $[n]$ which are all distinct,
		\item
			$
				E(\mu) =  .
				\left\{
					\left\{ {u}_{i}, {v}_{i} \right\}
				\, \middle| \,
					1 \leq i \leq p
				\right\} .
			$
	\end{enumerate}
We say a vertex $j \in V(\mu)$ is \textsl{single} if 
$
	j \nin 
		\left\{
			{u}_{1}, {u}_{2}, \ldots , {u}_{p}, {v}_{1}, {v}_{2}, \ldots , {v}_{p}
		\right\}
$. 
Especially, $\mu$ is called \textsl{perfect} if $n = 2p$. 
\end{definition}
Note that 
for each boundary vertex $v \in {V}_{1}(n)$, 
there uniquely exists a distinct boundary vertex $w \in {V}_{1}(n)$ such that
there is a monochromatic path $\boldsymbol{p}$ 
which has $v$ and $w$ as its end vertices. 
Then, $\psi$ determines $2$ kind of perfect matchings on $[2n]$
by focusing only on the black pathes
or focusing only on the white pathes. 
In addition, any two monochromatic paths of $\psi$
which have the same color
and start from and end at boundary vertices 
does not have a common vertex. 
We are ready to explain perfect matchings of boundary vertices 
determined by a FPL,
we define non-crossing matching. 
\begin{definition}
Let $n, p$ be positive integers, and $\mu$ be $p$-matching on $[n]$. 
We say $\mu$ is \textsl{non-crossing} if neitheir of the followings happens: 
\begin{subequations}
	\begin{align}
		\label{eq:crossing1}
		&\mu \textrm{ contains a pair } \{ {u}_{1}, {v}_{1} \} , \{ {u}_{2}, {v}_{2} \}
			\textrm{ of edges such that } {u}_{1} < {u}_{2} < {v}_{1} < {v}_{2}, \\
		\label{eq:crossing2}
		&\mu \textrm{ contains an edge } \{ {u}_{1}, {v}_{1} \} 
			\textrm{ and a single vertex } j 
				\textrm{ such that } {u}_{1} < j < {v}_{1}. 
	\end{align}
\end{subequations}
\end{definition}
\begin{figure}[htbp]
	\begin{center}
	\begin{tabular}{lr}
	\begin{minipage}{0.45\hsize}
		\centering
		\begin{tikzpicture}[scale=0.75]
			\coordinate (u1) at (0, 0) node at (u1) [below]{${u}_{1}$};
			\coordinate (u2) at (1, 0) node at (u2) [below]{${u}_{2}$};
			\coordinate (v1) at (3, 0) node at (v1) [below]{${v}_{1}$};
			\coordinate (v2) at (4, 0) node at (v2) [below]{${v}_{2}$};
			\draw[thick] (u1) to [bend left=60, distance=1cm] (v1);
			\draw[thick] (u2) to [bend left=60, distance=1cm] (v2);
			\foreach \x in {u1, u2, v1, v2} {
				\fill (\x) circle (2pt);
			}
		\end{tikzpicture}
		\subcaption{
			The case \eqref{eq:crossing1}
		}
		\label{pic:crossing1}
	\end{minipage}
	&
	\begin{minipage}{0.45\hsize}
		\centering
		\begin{tikzpicture}[scale=0.75]
			\coordinate (u1) at (0, 0) node at (u1) [below]{${u}_{1}$};
			\coordinate (v1) at (3, 0) node at (v1) [below]{${v}_{1}$};
			\coordinate (j) at (1, 0) node at (j) [below]{$j$};
			\draw[thick] (u1) to [bend left=60, distance=1cm] (v1);
			\foreach \x in {u1, v1, j} {
				\fill (\x) circle (2pt);
			}
		\end{tikzpicture}
		\subcaption{
			The case \eqref{eq:crossing2}
		}
		\label{pic:crossing2}
	\end{minipage}
	\end{tabular}
	\end{center}
	\caption{The cases which never happen in a non-crossing matching}
\end{figure}
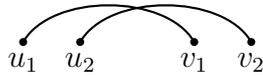
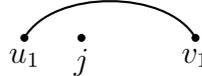
\subsection{Link patterns}
Now,
we shall introduce a link pattern.
Let us write 
non-crossing perfect matching on $[2n]$
with $2n$ vertices on the circumference
and $n$ edges drawn by arcs inside the circle,
and call this diagram a \textsl{link pattern of size $n$}. 
we denote $\mathcal{F}(2n)$ as the set of all link patterns of size $n$. 
 
\begin{figure}[htbp]
\begin{center}
\begin{tabular}{cc}
	\begin{tikzpicture}[scale=0.5]
	\draw (0, 0) circle [radius=2cm];
	\foreach \angle / \x in {45/7, 90/8, 135/1, 180/2, 225/3, 270/4, 315/5, 360/6} {
		\coordinate (v\x) at (\angle: 2);
		\fill (v\x) circle (2pt)
			node at ($(v\x) + (\angle: 10pt)$){$\x$};
	}
	\draw (v8) to [bend left=60, distance =0.5cm] (v1);
	\draw (v2) to [bend left=60, distance =1.5cm] (v5);
	\draw (v3) to [bend left=60, distance =0.5cm] (v4);
	\draw (v6) to [bend left=60, distance=0.5cm] (v7);
	\end{tikzpicture}
	&
	\begin{tikzpicture}[scale=0.75]
			\foreach \x in {1, 2, ..., 8}{
				\coordinate (v_\x) at (\x, 0); 
				\node at (v_\x) [below]{$\x$};
			}
			\draw[thick] (v_1) to [bend left=60, distance=1.5cm] (v_8);
			\draw[thick] (v_2) to [bend left=60, distance=1cm] (v_5);
			\draw[thick] (v_3) to [bend left=60, distance=0.5cm] (v_4);
			\draw[thick] (v_6) to [bend left=60, distance=0.5cm] (v_7);
			\foreach \x in {1, 2, ..., 8} {
				\fill (v_\x) circle (2pt);
			}
	\end{tikzpicture}
\end{tabular}
\end{center}
	\caption{An example of link pattern of size $4$}
	\label{pic:LP_example}
\end{figure}
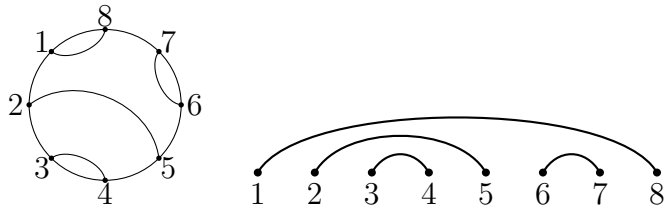
Then, we have 
$\# \mathcal{F}({2n}) = \frac{1}{n + 1} \binom{2n}{n}$, 
we denote ${C}_{n} = \frac{1}{n + 1} \binom{2n}{n}$
which is called the $n$th \textsl{Catalan number}. 
\par
\begin{figure}[H]
\begin{tabular}{ccc}
	\begin{minipage}{0.3\hsize}
	\centering
	\begin{tikzpicture}[scale = 0.6]
	\draw[->] (0, 5) -- (0, 0);
		\node at (0, 0)[below]{$i$};
	\draw[->] (0, 5) -- (5, 5);
		\node at (5, 5)[right]{$j$};
	\foreach \x in {1, ..., 4} {
		\coordinate (v\x) at (\x, 0) node at (v\x)[below]{
		};
	}
	\foreach \x in {5, ..., 8} {
		\coordinate (v\x) at ($(5, \x) + (0, -4)$) node at (v\x)[right]{
		};
	}
	\foreach \x in {9, ..., 12} {
		\coordinate (v\x) at ($(13, 5) - (\x, 0)$) node at (v\x)[above]{
		};
	}
	\foreach \x in {13, ..., 16} {
		\coordinate (v\x) at ($(0, 17) - (0, \x)$) node at (v\x)[left]{
		};
	}
	\foreach \y in {1, ..., 4} {
		\foreach \x in {1, ..., 4} {
			\coordinate (v\x_\y) at ($(\x, \y)$) node at (v\x_\y)[left]{
			};
		}
	}
	\foreach \x/\k in {1/6, 3/8}{
		\draw[thick] (\x, 1) -- (v\x);
			\node at (v\x)[below]{
				${e}_{\k}$
			};
	}
	\foreach \x/\k in {2/7, 4/9}{
		\draw[dashed] (\x, 1) -- (v\x);
			\node at (v\x)[below]{
				${e}_{\k}$
			};
	}
	\foreach \x/\k in {5/10, 7/12}{
		\draw[thick] (v\x) -- ($(4, \x) - (0, 4)$);
			\node at (v\x)[right]{
				${e}_{\k}$
			};
	}
	\foreach \x/\k in {6/11, 8/13}{
		\draw[dashed] (v\x) -- ($(4, \x) - (0, 4)$);
			\node at (v\x)[right]{
				${e}_{\k}$
			};
	}
	\foreach \x/\k in {9/14, 11/16}{
		\draw[thick] ($(13, 4) - (\x, 0)$) -- (v\x);
			\node at (v\x)[above]{
				${e}_{\k}$
			};
	}
	\foreach \x/\k in {10/15, 12/1}{
		\draw[dashed] ($(13, 4) - (\x, 0)$) -- (v\x);
			\node at (v\x)[above]{
				${e}_{\k}$
			};
	}
	\foreach \x/\k in {13/2, 15/4}{
		\draw[thick] (v\x) -- ($(1, 17) - (0, \x)$);
			\node at (v\x)[left]{
				${e}_{\k}$
			};
	}
	\foreach \x/\k in {14/4, 16/5}{
		\draw[dashed] (v\x) -- ($(1, 17) -(0, \x)$);
			\node at (v\x)[left]{
				${e}_{\k}$
			};
	}
	\foreach \x in {1, 2} {
		\draw[dashed] (v\x_4) -- ($(v\x_4) + (1, 0)$);
	}
	\foreach \x in {3} {
		\draw[thick] (v\x_4) -- ($(v\x_4) + (1, 0)$);
	}
	\foreach \x in {1, 2, 3} {
		\draw[dashed] (v\x_3) -- ($(v\x_3) + (1, 0)$);
	}
	\foreach \x in {1, 2} {
		\draw[dashed] (v\x_2) -- ($(v\x_2) + (1, 0)$);
	}
	\foreach \x in {3} {
		\draw[thick] (v\x_2) -- ($(v\x_2) + (1, 0)$);
	}
	\foreach \x in {2} {
		\draw[dashed] (v\x_1) -- ($(v\x_1) + (1, 0)$);
	}
	\foreach \x in {1, 3} {
		\draw[thick] (v\x_1) -- ($(v\x_1) + (1, 0)$);
	}
	\foreach \y in {1} {
		\draw[dashed] (v1_\y) -- ($(v1_\y) + (0, 1)$);
	}
	\foreach \y in {2, 3} {
		\draw[thick] (v1_\y) -- ($(v1_\y) + (0, 1)$);
	}
	\foreach \y in {1, 2, 3} {
		\draw[thick] (v2_\y) -- ($(v2_\y) + (0, 1)$);
	}
	
	\foreach \y in {1} {
		\draw[dashed] (v3_\y) -- ($(v3_\y) + (0, 1)$);
	}
	\foreach \y in {2, 3} {
		\draw[thick] (v3_\y) -- ($(v3_\y) + (0, 1)$);
	}
	\foreach \y in {1, 3} {
		\draw[dashed] (v4_\y) -- ($(v4_\y) + (0, 1)$);
	}
	\foreach \y in {2} {
		\draw[thick] (v4_\y) -- ($(v4_\y) + (0, 1)$);
	}
	\foreach \x in {1,3, 6, 8, 9, 11, 14, 16} {
		\fill (v\x)[fill=black, draw] circle (2pt); 
	}
	\foreach \x in {2, 4, 5, 7, 10, 12, 13, 15} {
		\fill (v\x)[fill=white, draw] circle (2pt);
	}
	\foreach \x in {1, 3} {
		\foreach \y in {1, 3} {
			\fill (v\x_\y)[fill=white, draw] circle (2pt);
		}
		\foreach \y in {2, 4} {
			\fill (v\x_\y)[fill=black, draw] circle (2pt);
		}
	}
	\foreach \x in {2, 4} {
		\foreach \y in {2, 4} {
			\fill (v\x_\y)[fill=white, draw] circle (2pt);
		}
		\foreach \y in {1, 3} {
			\fill (v\x_\y)[fill=black, draw] circle (2pt);
		}
	}
	\end{tikzpicture}
	\subcaption{$\psi \in \fpl(4, {\tau}_{-})$}
	\end{minipage}
&
	\begin{minipage}{0.3\hsize}
	\centering
	\begin{tikzpicture}[scale = 0.4]
	\foreach \x in {1, ..., 4} {
		\coordinate (v\x) at (\x, 0) node at (v\x)[below]{
		};
	}
	\foreach \x in {5, ..., 8} {
		\coordinate (v\x) at ($(5, \x) + (0, -4)$) node at (v\x)[right]{
		};
	}
	\foreach \x in {9, ..., 12} {
		\coordinate (v\x) at ($(13, 5) - (\x, 0)$) node at (v\x)[above]{
		};
	}
	\foreach \x in {13, ..., 16} {
		\coordinate (v\x) at ($(0, 17) - (0, \x)$) node at (v\x)[left]{
		};
	}
	\foreach \y in {1, ..., 4} {
		\foreach \x in {1, ..., 4} {
			\coordinate (v\x_\y) at ($(\x, \y)$) node at (v\x_\y)[left]{
			};
		}
	}
	\foreach \x in {1, 3}{
		\draw[thick] (\x, 1) -- (v\x)
		;
	}
	\foreach \x in {5, 7}{
		\draw[thick] (v\x) -- ($(4, \x) - (0, 4)$)
		;
	}
	\foreach \x in {9, 11}{
		\draw[thick] ($(13, 4) - (\x, 0)$) -- (v\x)
		;
	}
	\foreach \x in {13, 15}{
		\draw[thick] (v\x) -- ($(1, 17) - (0, \x)$)
		;
	}
	\foreach \x in {3} {
		\draw[thick] (v\x_4) -- ($(v\x_4) + (1, 0)$);
	}
	\foreach \x in {} {
		\draw[thick] (v\x_3) -- ($(v\x_3) + (1, 0)$);
	}
	\foreach \x in {3} {
		\draw[thick] (v\x_2) -- ($(v\x_2) + (1, 0)$);
	}
	\foreach \x in {1, 3} {
		\draw[thick] (v\x_1) -- ($(v\x_1) + (1, 0)$);
	}
	\foreach \y in {2, 3} {
		\draw[thick] (v1_\y) -- ($(v1_\y) + (0, 1)$);
	}
	\foreach \y in {1, 2, 3} {
		\draw[thick] (v2_\y) -- ($(v2_\y) + (0, 1)$);
	}
	\foreach \y in {2, 3} {
		\draw[thick] (v3_\y) -- ($(v3_\y) + (0, 1)$);
	}
	\foreach \y in {2} {
		\draw[thick] (v4_\y) -- ($(v4_\y) + (0, 1)$);
	}
	\foreach \x in {1,3,  ..., 15} {
		\fill (v\x)[fill=black, draw] circle (2pt); 
	}
	\draw (2.5, 2.5) circle [radius=4cm];
	\foreach \angle / \y in {0/6, 45/7, 90/8, 135/1, 180/2, 225/3, 270/4, 315/5} {
		\coordinate (w\y) at ($(\angle: 4cm) + (2.5, 2.5)$);
		\fill (w\y) circle (2pt) node at ($(w\y) + (\angle: 0.5cm)$){$\y$};
	}
	\foreach \x / \y in {1/3, 3/4, 5/5, 7/6, 9/7, 11/8, 13/1, 15/2} {
		\draw[thick] (v\x) -- (w\y);
	}
	\end{tikzpicture}
	\subcaption{${\pi}_{b, {\tau}_{-}}(\psi)$}
	\end{minipage}
&
	\begin{minipage}{0.3\hsize}
	\centering
	\begin{tikzpicture}[scale = 0.4]
	\foreach \x in {1, ..., 4} {
		\coordinate (v\x) at (\x, 0) node at (v\x)[below]{
		};
	}
	\foreach \x in {5, ..., 8} {
		\coordinate (v\x) at ($(5, \x) + (0, -4)$) node at (v\x)[right]{
		};
	}
	\foreach \x in {9, ..., 12} {
		\coordinate (v\x) at ($(13, 5) - (\x, 0)$) node at (v\x)[above]{
		};
	}
	\foreach \x in {13, ..., 16} {
		\coordinate (v\x) at ($(0, 17) - (0, \x)$) node at (v\x)[left]{
		};
	}
	\foreach \y in {1, ..., 4} {
		\foreach \x in {1, ..., 4} {
			\coordinate (v\x_\y) at ($(\x, \y)$) node at (v\x_\y)[left]{
			};
		}
	}
	\foreach \x in {2, 4}{
		\draw[dashed] (\x, 1) -- (v\x)
		;
	}
	\foreach \x in {6, 8}{
		\draw[dashed] (v\x) -- ($(4, \x) - (0, 4)$)
		;
	}
	\foreach \x in {10, 12}{
		\draw[dashed] ($(13, 4) - (\x, 0)$) -- (v\x)
		;
	}
	\foreach \x in {14, 16}{
		\draw[dashed] (v\x) -- ($(1, 17) -(0, \x)$)
		;
	}
	\foreach \x in {1, 2} {
		\draw[dashed] (v\x_4) -- ($(v\x_4) + (1, 0)$);
	}
	\foreach \x in {1, 2, 3} {
		\draw[dashed] (v\x_3) -- ($(v\x_3) + (1, 0)$);
	}
	\foreach \x in {1, 2} {
		\draw[dashed] (v\x_2) -- ($(v\x_2) + (1, 0)$);
	}
	\foreach \x in {2} {
		\draw[dashed] (v\x_1) -- ($(v\x_1) + (1, 0)$);
	}
	\foreach \y in {1} {
		\draw[dashed] (v1_\y) -- ($(v1_\y) + (0, 1)$);
	}
	\foreach \y in {1} {
		\draw[dashed] (v3_\y) -- ($(v3_\y) + (0, 1)$);
	}
	\foreach \y in {1, 3} {
		\draw[dashed] (v4_\y) -- ($(v4_\y) + (0, 1)$);
	}
	\foreach \x in {2, 4, ..., 16} {
		\fill (v\x)[fill=white, draw] circle (2pt);
	}
	\draw (2.5, 2.5) circle [radius=4cm];
	\foreach \angle / \y in {0/6, 45/7, 90/8, 135/1, 180/2, 225/3, 270/4, 315/5} {
		\coordinate (w\y) at ($(\angle: 4cm) + (2.5, 2.5)$);
		\fill (w\y) circle (2pt) node at ($(w\y) + (\angle: 0.5cm)$){$\y$};
	}
	\foreach \x / \y in {16/3, 2/4, 4/5, 6/6, 8/7, 10/8, 12/1, 14/2} {
		\draw[dashed] (v\x) -- (w\y);
	}
	\end{tikzpicture}
	\subcaption{${\pi}_{w, {\tau}_{-}}(\psi)$}
	\end{minipage}
\end{tabular}
\caption{An example of ${\pi}_{b, {\tau}_{-}}$ and ${\pi}_{w, {\tau}_{-}}$}
\label{pic:example_pi_minus_tikz}
\end{figure}
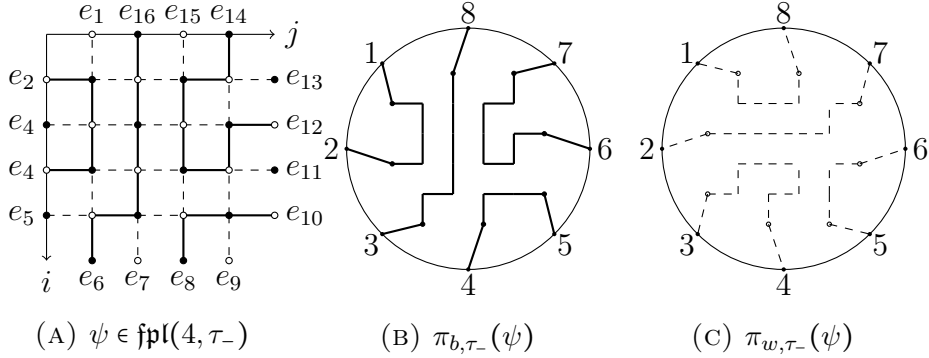

Here, we explain the $4$ maps from the set of FPLs to the set of link patterns.  
We define the following $4$ maps: 
	$
		{\pi}_{b, {\tau}_{-}}, 
		{\pi}_{w, {\tau}_{-}} : \mathfrak{fpl}(n, {\tau}_{-}) \rightarrow \mathcal{F}(2n)
	$
		and
	$
		{\pi}_{b, {\tau}_{+}}, 
		{\pi}_{w, {\tau}_{+}} : \mathfrak{fpl}(n, {\tau}_{+}) \rightarrow \mathcal{F}(2n)
	$.
\begin{itemize}
	\item The map 
		$
		{\pi}_{b, {\tau}_{-}} \; (\textrm{resp. } {\pi}_{w, {\tau}_{-}})
			 \colon \mathfrak{fpl}(n, {\tau}_{-}) \rightarrow \mathcal{F}(2n)
		$ 
		associates a link pattern with a given FPL as follows:
		if there is a black (resp. white) path which 
			connects ${e}_{2i}$ with ${e}_{2j}$ (resp. ${e}_{2i - 1}$ with ${e}_{2j - 1}$), 
			then we draw an edge between $i$ and $j$.  
\end{itemize}
The map 
$
	{\pi}_{b, {\tau}_{+}} \; (\textrm{resp. } {\pi}_{w, {\tau}_{+}})
		: \mathfrak{fpl}(n, {\tau}_{+}) \rightarrow \mathcal{F}(2n)
$
is obtained in the same way, by replacing 
${e}_{2i}$ (resp. ${e}_{2i - 1}$) with ${e}_{2i - 1}$ (resp. ${e}_{2i}$) and 
${e}_{2j}$ (resp. ${e}_{2j - 1}$) with ${e}_{2j - 1}$ (resp. ${e}_{2j}$). 
In Figure \ref{pic:example_pi_minus_tikz}, 
the black (resp. white) path of the FPL $\psi$ in (a) gives the link pattern (b) (resp. (c)).
Let $\mu \in \mathcal{F}({2n})$. 
Let ${\Psi}_{n, -}(\mu)$ (resp. ${\Psi}_{n, +}(\mu)$)
denote 
the cardinality of the set of FPLs with the boundary condition ${\tau}_{-}$ (resp. ${\tau}_{+}$)
and 
whose link pattern is $\mu$
(i.e., 
	$
		{\Psi}_{n, -}(\mu)	= 
			\# \left\{
				\psi \in \fpl(n, {\tau}_{-})
			\; \middle| \;
				{\pi}_{b, -}(\psi) = \mu
			\right\} 
	$, 
	$
		{\Psi}_{n, +}(\mu)	=
			\# \left\{
				\psi \in \fpl(n, {\tau}_{+})
			\; \middle| \;
				{\pi}_{b, +}(\psi) = \mu
			\right\} 
	$).
\subsection{Operators on link patterns}
\begin{figure}[H]
\centering
\begin{tabular}{ccc}
	\begin{minipage}{0.3\hsize}
		\centering
		\begin{tikzpicture}[scale=0.5]
			\draw (0, 0) circle [radius=2cm];
			\foreach \angle / \x in 
				{45/7, 90/8, 135/1, 180/2, 225/3, 270/4, 315/5, 360/6}
			{
				\coordinate (v\x) at (\angle: 2);
				\fill (v\x) circle (2pt)
				node at ($(v\x) + (\angle: 10pt)$){$\x$};
			}
			\draw (v8) to [bend left=60, distance =0.5cm] (v1);
			\draw (v2) to [bend left=60, distance =1.5cm] (v5);
			\draw (v3) to [bend left=60, distance =0.5cm] (v4);
			\draw (v6) to [bend left=60, distance=0.5cm] (v7);
		\end{tikzpicture}
	\end{minipage}
	&
		$\overset{\longmapsto}{ {e}_{2} }$
	&
	\begin{minipage}{0.3\hsize}
		\centering
		\begin{tikzpicture}[scale=0.5]
			\draw (0, 0) circle [radius=2cm];
			\foreach \angle / \x in 
				{45/7, 90/8, 135/1, 180/2, 225/3, 270/4, 315/5, 360/6}
			{
				\coordinate (v\x) at (\angle: 2);
				\fill (v\x) circle (2pt)
				node at ($(v\x) + (\angle: 10pt)$){$\x$};
			}
			\draw (v8) to [bend left=60, distance =0.5cm] (v1);
			\draw (v2) to [bend left=60, distance =0.5cm] (v3);
			\draw (v4) to [bend left=60, distance =0.5cm] (v5);
			\draw (v6) to [bend left=60, distance =0.5cm] (v7);
		\end{tikzpicture}
	\end{minipage}
\end{tabular}
\caption{Example in the case when $\{ j, j + 1 \} \nin E(\mu)$}
\label{pic:matchmaker_1}
\end{figure}
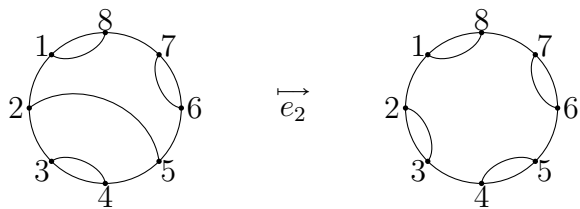

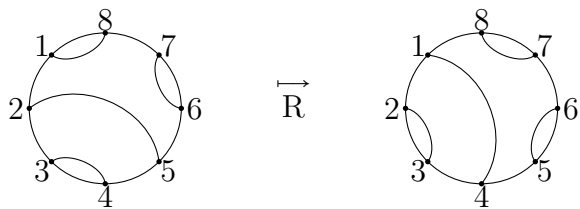
\begin{figure}[H]
	\centering
	\begin{tabular}{ccc}
	\begin{minipage}{0.3\hsize}
		\centering
		\begin{tikzpicture}[scale=0.5]
			\draw (0, 0) circle [radius=2cm];
			\foreach \angle / \x in 
				{45/7, 90/8, 135/1, 180/2, 225/3, 270/4, 315/5, 360/6}
			{
				\coordinate (v\x) at (\angle: 2);
				\fill (v\x) circle (2pt)
				node at ($(v\x) + (\angle: 10pt)$){$\x$};
			}
			\draw (v8) to [bend left=60, distance =0.5cm] (v1);
			\draw (v2) to [bend left=60, distance =1.5cm] (v5);
			\draw (v3) to [bend left=60, distance =0.5cm] (v4);
			\draw (v6) to [bend left=60, distance=0.5cm] (v7);
		\end{tikzpicture}
	\end{minipage}
	&
		$\overset{\longmapsto}{\operatorname{R}}$
	&
	\begin{minipage}{0.3\hsize}
		\centering
		\begin{tikzpicture}[scale=0.5]
			\draw (0, 0) circle [radius=2cm];
			\foreach \angle / \x in 
				{45/7, 90/8, 135/1, 180/2, 225/3, 270/4, 315/5, 360/6}
			{
				\coordinate (v\x) at (\angle: 2);
				\fill (v\x) circle (2pt)
				node at ($(v\x) + (\angle: 10pt)$){$\x$};
			}
			\draw (v7) to [bend left=60, distance =0.5cm] (v8);
			\draw (v1) to [bend left=60, distance =1.5cm] (v4);
			\draw (v2) to [bend left=60, distance =0.5cm] (v3);
			\draw (v5) to [bend left=60, distance =0.5cm] (v6);
		\end{tikzpicture}
	\end{minipage}
\end{tabular}
\caption{The example of $R$}
\label{pic:R_example}
\end{figure}

Here, we define operators on link patterns 
which are useful to examine the behavior of 
the ${\Psi}_{n, -}$. 
\begin{definition}
For a positive integer $n$ 
and a positive integer $j$ which satisfies $1 \leq j \leq 2n$, 
we define an operator 
${e}_{j} \colon \mathcal{F}({2n}) \rightarrow \mathcal{F}({2n})$ 
which is called \textsl{matchmaker} (\cite{dG}). 
For
$\mu \in \mathcal{F}({2n})$, 
we define ${e}_{j}\mu \in \mathcal{F}(2n)$ as follows:
\begin{enumerate}[(i)]
	\item If $\{ j, j + 1 \} \in E(\mu)$, 
		$
			{e}_{j}\mu := \mu
		$.
	\item If $\{ j, j + 1 \} \nin E(\mu)$,
		the edge set of ${e}_{j}\mu$ 
		is defined as
		$
				\left\{
					{\eta}_{1}, {\eta}_{2}, \cdots , {\eta}_{2n - 2}, 
					\{ j, j + 1 \} , \{ u, v \}
				\right\}
		$
		when the edge set of $\mu$ equals
		$
				\left\{
					{\eta}_{1}, {\eta}_{2}, \cdots , {\eta}_{2n - 2}, 
					\{ j, u \} , \{ j + 1, v \}
				\right\}
		$. 
\end{enumerate}
\end{definition}
Figure \ref{pic:matchmaker_1} illustlates the operator ${e}_{2}$.
Next, we define the rotation operator 
$\operatorname{R} \colon \mathcal{F}(2n) \rightarrow \mathcal{F}(2n)$. 
\begin{definition}
For
$\mu \in \mathcal{F}(2n)$, 
we define $\operatorname{R}\mu \in \mathcal{F}(2n)$ as follows: 
The edge set of $\operatorname{R}\mu$ is defined as
$
	\left\{
		\{ {u}_{1} - 1, {v}_{1} - 1 \} , \{ {u}_{2} - 1, {v}_{2} - 1 \} , \cdots , \{ {u}_{n} - 1, {v}_{n} - 1 \}
	\right\}
$ 
when the edge set of $\mu$ equals
$ 
	\left\{
		\{ {u}_{1}, {v}_{1} \} , \{ {u}_{2}, {v}_{2} \} , \cdots , \{ {u}_{n}, {v}_{n} \}
	\right\}
$.
\end{definition}
Figure \ref{pic:R_example} illustrates the rotation operator $\operatorname{R}$. 
\section{Dihedral symmetry and Gylation}
Let $n, l$ be positive integers. For ${\mu}_{b}, {\mu}_{w} \in \mathcal{F}(2n)$, 
we define a
${\Psi}_{n, -}\left( {\mu}_{b}, {\mu}_{w}; l \right)$
(resp.
${\Psi}_{n, +}\left( {\mu}_{b}, {\mu}_{w}; l \right)$) 
to refine ${\Psi}_{n, -}({\mu}_{b})$
(resp. ${\Psi}_{n, +}({\mu}_{b})$)
as follows.  
Let
${\Psi}_{n, -}\left( {\mu}_{b}, {\mu}_{w}; l \right)$ 
(resp.
${\Psi}_{n, +}\left( {\mu}_{b}, {\mu}_{w}; l \right)$) 
denote
the cardinality of 
the FPLs 
in $\fpl(n, {\tau}_{-})$
(resp. $\fpl(n, {\tau}_{+})$) 
whose link patterns determined by black (resp. white) paths equals ${\mu}_{b}$ (resp. ${\mu}_{w}$), 
and have $l$ monochromatic cycles. 
\subsection{Wieland's dihedral symmetry theorem}
Now, the following proposition
which is called \textsl{Wieland's dihedral symmetry theorem}
we have as a characterization of ${\Psi}_{n, -}$. 
\begin{proposition}\label{prop:Wieland}
	Let $n, l$ be positive integers. 
	For ${\mu}_{b}, {\mu}_{w} \in \mathcal{F}(2n)$, we hold 
	\begin{align}
		\label{eq:DS_m}
		{\Psi}_{n, -}\left( {\mu}_{b}, {\mu}_{w}; l \right)
			&=	{\Psi}_{n, -}\left( \operatorname{R}^{- 1} {\mu}_{b}, \operatorname{R} {\mu}_{w}; l \right). 
	\end{align}
\end{proposition} 
The proposition \ref{prop:Wieland} is proven 
by constructing the bijection which is called \textrm{gyration}. 
Next we explain the gyration. 
\subsection{Gyration}
We construct the gyration
$\operatorname{G} \colon \fpl(n, {\tau}_{-}) \rightarrow \fpl(n, {\tau}_{-})$
in two ways.
First, we define two maps
$\operatorname{G}_{0} \colon \fpl(n, {\tau}_{-}) \rightarrow \fpl(n, {\tau}_{-})$
and
$\operatorname{G}_{1} \colon \fpl(n, {\tau}_{-}) \rightarrow \fpl(n, {\tau}_{-})$,
then we construct $\operatorname{G}$ as 
$\operatorname{G}_{0} \operatorname{G}_{1}$. 
Second, we define two maps
$\operatorname{H}_{0} \colon \fpl(n, {\tau}_{+}) \rightarrow \fpl(n, {\tau}_{-})$
and
$\operatorname{H}_{1} \colon \fpl(n, {\tau}_{-}) \rightarrow \fpl(n, {\tau}_{+})$,
then we construct $\operatorname{G}$ as 
$\operatorname{H}_{0} \operatorname{H}_{1}$. 

Let $\alpha$ be a plaquette. 
We define a map 
$
	\operatorname{G}_{\alpha} \colon 
		\fpl(n, {\tau}_{-}) \rightarrow \fpl(n, {\tau}_{-})
$
which 
affects only the colors of 
edges in $\alpha$ 
and leaves the colors of other edges invariant.

Let $\alpha$ be an interior plaquette.
First of all, we define a map 
$\mathcal{N}_{\alpha} \colon \fpl(n, {\tau}_{-}) \rightarrow \{ 0, 1, - 1 \}$ 
to represent the state of $\alpha$
as follows.
\begin{definition}
	Let $n$ be a positive integer, and
	${\alpha}_{i, j}$ an interior plaquette. 
	Here, we label four edges of $\alpha$ as
	${\eta}_{1} = \left\{ \{ i, j \} , \{ i + 1, j \} \right\} $, 
	${\eta}_{2} = \left\{ \{ i + 1, j \} , \{ i + 1, j + 1 \} \right\} $, 
	${\eta}_{3} = \left\{ \{ i, j + 1 \} , \{ i + 1, j + 1 \} \right\} $ and
	${\eta}_{4} = \left\{ \{ i, j \} , \{ i, j + 1 \} \right\} $. 
	For $\psi \in \fpl(n, {\tau}_{-})$, 
	we define
	$\mathcal{N}_{{\alpha}_{i, j}} (\psi) \in \{ 0, 1, -1 \} $ 
	as follows: 
	\begin{equation}
		\mathcal{N}_{{\alpha}_{i, j}}(\psi) := 
			\begin{cases}
				1	&
					\textrm{if }
					\psi\left( 
						{\eta}_{1}
					\right)
					= 
						\psi\left( 
							{\eta}_{3} 
						\right)
					= w
					\textrm{ and }
					\psi\left( 
						{\eta}_{2} 
					\right)
					= 
						\psi\left( 
							{\eta}_{4} 
						\right)
					= b,
					\\
				- 1	&
					\textrm{if }
					\psi\left( 
						{\eta}_{1} 
					\right)
					= 
						\psi\left( 
							{\eta}_{3} 
						\right)
					= b
					\textrm{ and }
					\psi\left( 
						{\eta}_{2} 
					\right)
					= 
						\psi\left( 
							{\eta}_{4} 
						\right)
					= w, 
					\\
				0	&	\textrm{otherwise}.
			\end{cases}
	\end{equation}
\end{definition}
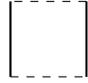
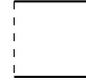
\begin{figure}[H]
\centering
\begin{tabular}{cc}
	\begin{minipage}{0.45\hsize}
		\centering
		\begin{tikzpicture}
			\coordinate (v1) at (0, 0);
			\coordinate (v2) at (1, 0);
			\coordinate (v3) at (1, 1);
			\coordinate (v4) at (0, 1);
			\draw[dashed] (v1) -- (v2);
			\draw[thick] (v2) -- (v3);
			\draw[dashed] (v3) -- (v4);
			\draw[thick] (v4) -- (v1);
		\end{tikzpicture}
		\subcaption{ 
			$\mathcal{N}_{{\alpha}_{i, j}} = 1$.}
		\label{pic:N_alpha_plus}
	\end{minipage}
&
	\begin{minipage}{0.45\hsize}
		\centering
		\begin{tikzpicture}
			\coordinate (v1) at (0, 0);
			\coordinate (v2) at (1, 0);
			\coordinate (v3) at (1, 1);
			\coordinate (v4) at (0, 1);
			\draw[thick] (v1) -- (v2);
			\draw[dashed] (v2) -- (v3);
			\draw[thick] (v3) -- (v4);
			\draw[dashed] (v4) -- (v1);
		\end{tikzpicture}
		\subcaption{
			$\mathcal{N}_{{\alpha}_{i, j}} = -1$.}
		\label{pic:N_alpha_minus}
	\end{minipage}	
\end{tabular}
\caption{The states of ${\alpha}_{i, j}$ which satisfies $\mathcal{N}_{{\alpha}_{i, j}} \neq 0$.}
\end{figure}
Then we define an operator 
$
	\operatorname{G}_{\alpha} \colon 
		\fpl(n, {\tau}_{-}) \rightarrow \fpl(n, {\tau}_{-})
$
as follows.
We remark that $\operatorname{G}_{\alpha}$ is defined for 
not only interior plaquettes but also boundary plaquettes. 
\begin{definition}
	Let $n$ be a positive integer, and $\alpha$ a plaquette.
	For $\psi \in \fpl(n, {\tau}_{-})$, 
	we define $\operatorname{G}_{\alpha}\psi \in \fpl(n, {\tau}_{-})$
	as follows:
	\begin{enumerate}[(i)]
		\item	$\operatorname{G}_{\alpha}$ does not affect the colors of edges not in $\alpha$
			(i.e., $\operatorname{G}_{\alpha}\psi(e) = \psi(e)$ 
				for $e \in E({L}_{n}) \setminus E(\alpha)$).   
		\item	If $\alpha$ is an interior plaquette and $\mathcal{N}_{\alpha}(\psi) \neq 0$
			(resp. $\mathcal{N}_{\alpha}(\psi) = 0$), 
			$\operatorname{G}_{\alpha}$ reverses (resp. keeps) the colors of  
				all edges in $\alpha$ .
		\item	
			If $\alpha$ is a boundary plaquette, 
			$\operatorname{G}_{\alpha}$ keeps the colors of 
				all edges in $\alpha$.  
	\end{enumerate} 
\end{definition} 
For two distinct plaquettes $\alpha$ and $\beta$
which have the same parity,
the composition of
$\operatorname{G}_{\alpha}$ and $\operatorname{G}_{\beta}$
commutes with each other
because
$\alpha$ and $\beta$ don't have a common edge
(i.e., 
$
	\operatorname{G}_{\alpha} \operatorname{G}_{\beta}
		=	\operatorname{G}_{\beta} \operatorname{G}_{\alpha}
$).
Then we can define
$
	\operatorname{G}_{0} \; (\textrm{resp. } \operatorname{G}_{1} )
		\colon \fpl(n, {\tau}_{-}) \rightarrow \fpl(n, {\tau}_{-})
$
as
$\prod_{\alpha: \even} \operatorname{G}_{\alpha}$
(resp. $\prod_{\alpha: \odd} \operatorname{G}_{\alpha}$)
without ambiguity.
\input tikz_gyration_example_ver3.tex
Up to here, 
we have constructed 
$\operatorname{G}_{0}$ and $\operatorname{G}_{1}$
by repeating the operations to FPL.
But in the following, 
we construct
$\operatorname{H}_{0}$ and $\operatorname{H}_{1}$
which repeat
the operations to a coloring
which is not a FPL.
First
we denote the set of all maps
from $E({L}_{n})$ to $\{ b, w \} $
by $\mathrm{Map}(E({L}_{n}), \{ b, w \} )$.
For any map $f \colon E({L}_{n}) \rightarrow \{ b, w \} $,
we also call $f$ a  
\textsl{coloring of $E({L}_{n})$}.

Now
we define
$
	\operatorname{C}_{\alpha} \colon 
		\mathrm{Map}(E({L}_{n}), \{ b, w \} ) \rightarrow \mathrm{Map}(E({L}_{n}), \{ b, w \} )
$
as the map which
reverse the colors of all edges in $\alpha$
and leaves the colors of other edges invariant.

Note that
$\mathcal{N}_{\alpha}$ can be extended to
a map from $\mathrm{Map}(E({L}_{n}), \{ b, w \} )$ to $\{ 0, 1, -1 \}$
when $\alpha$ is an interior plaquette,
and $\operatorname{G}_{\alpha}$ to
a map from $\mathrm{Map}(E({L}_{n}), \{ b, w \} )$ to 
$\mathrm{Map}(E({L}_{n}), \{ b, w \} )$.
Then we define an operator 
$
	\operatorname{H}_{\alpha} \colon
		\mathrm{Map}(E({L}_{n}), \{ b, w \} ) \rightarrow \mathrm{Map}(E({L}_{n}), \{ b, w \} )
$
as $\operatorname{G}_{\alpha} \operatorname{C}_{\alpha}$.

Let $\alpha$ and $\beta$ be 
plaquettes which have the same parity.
Here
$\alpha$ and $\beta$ do not have to be distinct.
Notice that
the composition of
$\operatorname{G}_{\alpha}$ and $\operatorname{C}_{\beta}$
is commutative
(i.e., 
$
	\operatorname{G}_{\alpha} \operatorname{C}_{\beta}
	= \operatorname{C}_{\beta} \operatorname{G}_{\alpha}
$).
Hence, 
the composition of $\operatorname{H}_{\alpha}$ and $\operatorname{H}_{\beta}$
is commutative
when $\alpha$ and $\beta$ is distinct.
Therefore we can define 
$\prod_{\alpha \colon \even} \operatorname{H}_{\alpha}$ and 
$\prod_{\alpha \colon \odd} \operatorname{H}_{\alpha}$
without ambiguity.
Moreover,
$\prod_{\alpha \colon \even} \operatorname{H}_{\alpha}$ 
(resp. $\prod_{\alpha \colon \odd} \operatorname{H}_{\alpha}$)
equals 
$
\prod_{\alpha \colon \even} \operatorname{G}_{\alpha} 
	\prod_{\alpha \colon \even} \operatorname{C}_{\alpha}
$
(resp.
$
\prod_{\alpha \colon \odd} \operatorname{G}_{\alpha} 
	\prod_{\alpha \colon \odd} \operatorname{C}_{\alpha}
$).

Note that
$\operatorname{G}_{\alpha}$, $\operatorname{C}_{\alpha}$ and $\operatorname{H}_{\alpha}$
are involutions. 
Since $\operatorname{G}_{0}$, $\operatorname{G}_{1}$, 
$\prod_{\alpha \colon \even} \operatorname{H}_{\alpha}$ and
$\prod_{\alpha \colon \odd} \operatorname{H}_{\alpha}$
are compositions of these commutative involutions, 
they are involutions.

We remark that
$E({L}_{n})$ can be expressed
$\bigsqcup_{\alpha \colon \even} E(\alpha)$ or
$\bigsqcup_{\alpha \colon \odd} E(\alpha)$.
Therefore 
$\bigsqcup_{\alpha \colon \even} \operatorname{C}_{\alpha}$
(resp. $\bigsqcup_{\alpha \colon \odd} \operatorname{C}_{\alpha}$)
reverse the colors of all edges in ${L}_{n}$.
Hence
we can define a bijection
$
	\operatorname{H}_{0} \colon
		\fpl(n, {\tau}_{+}) \rightarrow \fpl(n, {\tau}_{-})
$
(resp.
$
	\operatorname{H}_{1} \colon
		\fpl(n, {\tau}_{-}) \rightarrow \fpl(n, {\tau}_{+})
$)
by restricting
$\prod_{\alpha \colon \even} \operatorname{H}_{\alpha}$
(resp. $\prod_{\alpha \colon \odd} \operatorname{H}_{\alpha}$)
to $\fpl(n, {\tau}_{+})$ (resp. $\fpl(n, {\tau}_{-})$).
We use $\operatorname{H}_{0}$ and $\operatorname{H}_{1}$ to prove the following lemma.
\begin{lemma}\label{lem:H0H1}
	Let $n$ be a positive integer.
	For ${\mu}_{b}, {\mu}_{w} \in \mathcal{F}(2n)$, we have following two equations:
	\begin{subequations}
	\begin{align}
		\label{eq:H1_lem}
		{\Psi}_{n, -} \left( {\mu}_{b}, {\mu}_{w}; l \right) 
		&= 
		{\Psi}_{n, +} \left( \operatorname{R}^{-1} {\mu}_{b}, \operatorname{R} {\mu}_{w}; l \right) , 
		\\
		\label{eq:H0_lem}
		{\Psi}_{n, +} \left( {\mu}_{b}, {\mu}_{w}; l \right) 
		&= 
		{\Psi}_{n, -} \left( {\mu}_{b}, {\mu}_{w}; l \right) .
	\end{align}
	\end{subequations}
\end{lemma}
We take several steps to prove lemma \ref{lem:H0H1}.
First, we define vertices which are called fixed vertice of $\operatorname{H}_{1}$ 
(resp. $\operatorname{H}_{0})$
to prove \eqref{eq:H1_lem} (resp. \eqref{eq:H0_lem}).
Second, we divide a monochromatic path at each fixed vertex into short pathes.
Then we show how the short pathes are affected by 
$\operatorname{H}_{1}$ (resp. $\operatorname{H}_{0})$.

Let $\psi$ be a FPL in $\fpl(n, {\tau}_{-})$ 
and $v$ an interior vertex.
We label the two edges which are adjacent to $v$ and have color $b$ in $\psi$ as $e$ and $e'$.
Now, 
we remark that there are two plaquettes which contain $v$ and whose parity is odd (resp. even).
Then
we call $v$ a \textsl{fixed vertex of $\operatorname{H}_{1}$ (resp. $\operatorname{H}_{0}$) in $\psi$}
if $e$ and $e'$ are the edges of distinct plaquettes 
whose parity is odd (resp. even) and which contain $v$.
For any ${\psi}' \in \fpl(n, {\tau}_{+})$,
we also define a \textsl{fixed vertex of $\operatorname{H}_{1}$ (resp. $\operatorname{H}_{0}$) in ${\psi}'$}
in the same way.

Let $v = (i, j)$ be an interior veretex, and $k \in \{ 0, 1 \}$.
Figure \ref{fig:fixed_H0_even} 
(resp. Figure \ref{fig:fixed_H1_even})
shows whether $v$ is fixed vertex or not
if $i + j$ and $k$ have the same parity
(resp. don't have the same parity).
We illustrate $v$ as $\bullet$
and we color the plaquette which has the same parity of $k$
in Figure \ref{fig:fixed_H0_even} 
(resp. Figure \ref{fig:fixed_H1_even}).
\begin{figure}[htbp]
\begin{tabular}{ccc}
\begin{minipage}{0.3\hsize}
\centering
\begin{tikzpicture}
	\coordinate (O) at (0, 0);
	\coordinate (A) at (0, 1);
	\coordinate (B) at (1, 0);
	\coordinate (C) at (0, -1); 
	\coordinate (D) at (-1, 0);
	\fill[fill=lightgray, opacity=0.5] (0, 0) -- (0, 1) -- (-1, 1) -- (-1, 0) --cycle;
	\fill[fill=lightgray, opacity=0.5] (0, 0) -- (0, -1) -- (1, -1) -- (1, 0) -- cycle;
	\draw[dashed] (O) -- (A);
	\draw[dashed] (O) -- (B);
	\draw[thick] (O) -- (C);
	\draw[thick] (O) -- (D);
	\fill (O) [fill=black, draw] circle (2pt);
	\foreach \x in {A, B, C, D} {
		\fill (\x) [fill=white, draw] circle (2pt);
	}
\end{tikzpicture}
\subcaption{fixed vertex of $\operatorname{H}_{k}$}
\end{minipage}
&
\begin{minipage}{0.3\hsize}
\centering
\begin{tikzpicture}
	\coordinate (O) at (0, 0);
	\coordinate (A) at (0, 1);
	\coordinate (B) at (1, 0);
	\coordinate (C) at (0, -1); 
	\coordinate (D) at (-1, 0);
	\fill[fill=lightgray, opacity=0.5] (0, 0) -- (0, 1) -- (-1, 1) -- (-1, 0) --cycle;
	\fill[fill=lightgray, opacity=0.5] (0, 0) -- (0, -1) -- (1, -1) -- (1, 0) -- cycle;
	\draw[dashed] (O) -- (A);
	\draw[thick] (O) -- (B);
	\draw[dashed] (O) -- (C);
	\draw[thick] (O) -- (D);
	\fill (O) [fill=black, draw] circle (2pt);
	\foreach \x in {A, B, C, D} {
		\fill (\x) [fill=white, draw] circle (2pt);
	}
\end{tikzpicture}
\subcaption{fixed vertex of $\operatorname{H}_{k}$}
\end{minipage}
&
\begin{minipage}{0.3\hsize}
\centering
\begin{tikzpicture}
	\coordinate (O) at (0, 0);
	\coordinate (A) at (0, 1);
	\coordinate (B) at (1, 0);
	\coordinate (C) at (0, -1); 
	\coordinate (D) at (-1, 0);
	\fill[fill=lightgray, opacity=0.5] (0, 0) -- (0, 1) -- (-1, 1) -- (-1, 0) --cycle;
	\fill[fill=lightgray, opacity=0.5] (0, 0) -- (0, -1) -- (1, -1) -- (1, 0) -- cycle;
	\draw[dashed] (O) -- (A);
	\draw[thick] (O) -- (B);
	\draw[thick] (O) -- (C);
	\draw[dashed] (O) -- (D);
	\fill (O) [fill=black, draw] circle (2pt);
	\foreach \x in {A, B, C, D} {
		\fill (\x) [fill=white, draw] circle (2pt);
	}
\end{tikzpicture}
\subcaption{not fixed vertex of $\operatorname{H}_{k}$}
\end{minipage}
\\
\begin{minipage}{0.3\hsize}
\centering
\begin{tikzpicture}
	\coordinate (O) at (0, 0);
	\coordinate (A) at (0, 1);
	\coordinate (B) at (1, 0);
	\coordinate (C) at (0, -1); 
	\coordinate (D) at (-1, 0);
	\fill[fill=lightgray, opacity=0.5] (0, 0) -- (0, 1) -- (-1, 1) -- (-1, 0) --cycle;
	\fill[fill=lightgray, opacity=0.5] (0, 0) -- (0, -1) -- (1, -1) -- (1, 0) -- cycle;
	\draw[thick] (O) -- (A);
	\draw[dashed] (O) -- (B);
	\draw[dashed] (O) -- (C);
	\draw[thick] (O) -- (D);
	\fill (O) [fill=black, draw] circle (2pt);
	\foreach \x in {A, B, C, D} {
		\fill (\x) [fill=white, draw] circle (2pt);
	}
\end{tikzpicture}
\subcaption{not fixed vertex of $\operatorname{H}_{k}$}
\end{minipage}
&
\begin{minipage}{0.3\hsize}
\centering
\begin{tikzpicture}
	\coordinate (O) at (0, 0);
	\coordinate (A) at (0, 1);
	\coordinate (B) at (1, 0);
	\coordinate (C) at (0, -1); 
	\coordinate (D) at (-1, 0);
	\fill[fill=lightgray, opacity=0.5] (0, 0) -- (0, 1) -- (-1, 1) -- (-1, 0) --cycle;
	\fill[fill=lightgray, opacity=0.5] (0, 0) -- (0, -1) -- (1, -1) -- (1, 0) -- cycle;
	\draw[thick] (O) -- (A);
	\draw[dashed] (O) -- (B);
	\draw[thick] (O) -- (C);
	\draw[dashed] (O) -- (D);
	\fill (O) [fill=black, draw] circle (2pt);
	\foreach \x in {A, B, C, D} {
		\fill (\x) [fill=white, draw] circle (2pt);
	}
\end{tikzpicture}
\subcaption{fixed vertex of $\operatorname{H}_{k}$}
\end{minipage}
&
\begin{minipage}{0.3\hsize}
\centering
\begin{tikzpicture}
	\coordinate (O) at (0, 0);
	\coordinate (A) at (0, 1);
	\coordinate (B) at (1, 0);
	\coordinate (C) at (0, -1); 
	\coordinate (D) at (-1, 0);
	\fill[fill=lightgray, opacity=0.5] (0, 0) -- (0, 1) -- (-1, 1) -- (-1, 0) --cycle;
	\fill[fill=lightgray, opacity=0.5] (0, 0) -- (0, -1) -- (1, -1) -- (1, 0) -- cycle;
	\draw[thick] (O) -- (A);
	\draw[thick] (O) -- (B);
	\draw[dashed] (O) -- (C);
	\draw[dashed] (O) -- (D);
	\fill (O) [fill=black, draw] circle (2pt);
	\foreach \x in {A, B, C, D} {
		\fill (\x) [fill=white, draw] circle (2pt);
	}
\end{tikzpicture}
\subcaption{fixed vertex of $\operatorname{H}_{k}$}
\end{minipage}
\end{tabular}
\caption{If $i + j$ and $k$ have the same parity}
\label{fig:fixed_H0_even}
\end{figure}
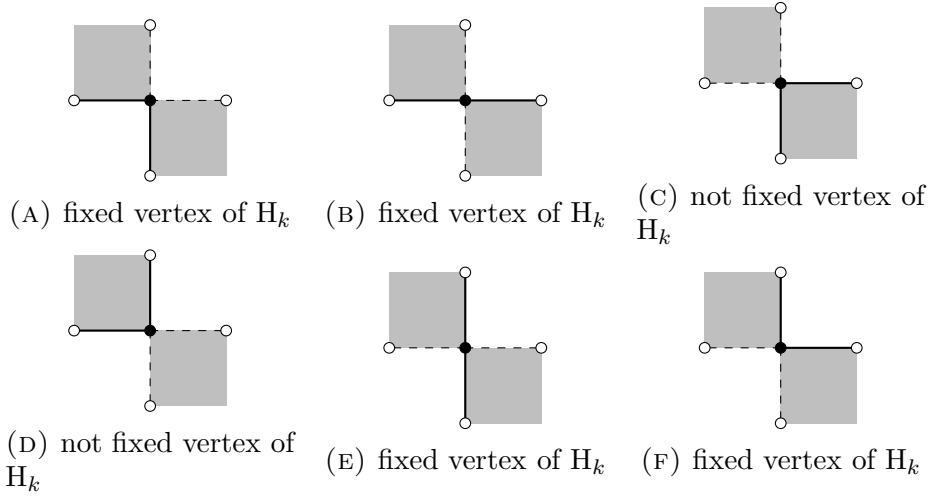
\begin{figure}[htbp]
\begin{tabular}{ccc}
\begin{minipage}{0.3\hsize}
\centering
\begin{tikzpicture}
	\coordinate (O) at (0, 0);
	\coordinate (A) at (0, 1);
	\coordinate (B) at (1, 0);
	\coordinate (C) at (0, -1); 
	\coordinate (D) at (-1, 0);
	\fill[fill=lightgray, opacity=0.5] (0, 0) -- (0, 1) -- (1, 1) -- (1, 0) --cycle;
	\fill[fill=lightgray, opacity=0.5] (0, 0) -- (0, -1) -- (-1, -1) -- (-1, 0) -- cycle;
	\draw[dashed] (O) -- (A);
	\draw[dashed] (O) -- (B);
	\draw[thick] (O) -- (C);
	\draw[thick] (O) -- (D);
	\fill (O) [fill=black, draw] circle (2pt);
	\foreach \x in {A, B, C, D} {
		\fill (\x) [fill=white, draw] circle (2pt);
	}
\end{tikzpicture}
\subcaption{not fixed vertex of $\operatorname{H}_{k}$}
\end{minipage}
&
\begin{minipage}{0.3\hsize}
\centering
\begin{tikzpicture}
	\coordinate (O) at (0, 0);
	\coordinate (A) at (0, 1);
	\coordinate (B) at (1, 0);
	\coordinate (C) at (0, -1); 
	\coordinate (D) at (-1, 0);
	\fill[fill=lightgray, opacity=0.5] (0, 0) -- (0, 1) -- (1, 1) -- (1, 0) --cycle;
	\fill[fill=lightgray, opacity=0.5] (0, 0) -- (0, -1) -- (-1, -1) -- (-1, 0) -- cycle;
	\draw[dashed] (O) -- (A);
	\draw[thick] (O) -- (B);
	\draw[dashed] (O) -- (C);
	\draw[thick] (O) -- (D);
	\fill (O) [fill=black, draw] circle (2pt);
	\foreach \x in {A, B, C, D} {
		\fill (\x) [fill=white, draw] circle (2pt);
	}
\end{tikzpicture}
\subcaption{fixed vertex of $\operatorname{H}_{k}$}
\end{minipage}
&
\begin{minipage}{0.3\hsize}
\centering
\begin{tikzpicture}
	\coordinate (O) at (0, 0);
	\coordinate (A) at (0, 1);
	\coordinate (B) at (1, 0);
	\coordinate (C) at (0, -1); 
	\coordinate (D) at (-1, 0);
	\fill[fill=lightgray, opacity=0.5] (0, 0) -- (0, 1) -- (1, 1) -- (1, 0) --cycle;
	\fill[fill=lightgray, opacity=0.5] (0, 0) -- (0, -1) -- (-1, -1) -- (-1, 0) -- cycle;
	\draw[dashed] (O) -- (A);
	\draw[thick] (O) -- (B);
	\draw[thick] (O) -- (C);
	\draw[dashed] (O) -- (D);
	\fill (O) [fill=black, draw] circle (2pt);
	\foreach \x in {A, B, C, D} {
		\fill (\x) [fill=white, draw] circle (2pt);
	}
\end{tikzpicture}
\subcaption{fixed vertex of $\operatorname{H}_{k}$}
\end{minipage}
\\
\begin{minipage}{0.3\hsize}
\centering
\begin{tikzpicture}
	\coordinate (O) at (0, 0);
	\coordinate (A) at (0, 1);
	\coordinate (B) at (1, 0);
	\coordinate (C) at (0, -1); 
	\coordinate (D) at (-1, 0);
	\fill[fill=lightgray, opacity=0.5] (0, 0) -- (0, 1) -- (1, 1) -- (1, 0) --cycle;
	\fill[fill=lightgray, opacity=0.5] (0, 0) -- (0, -1) -- (-1, -1) -- (-1, 0) -- cycle;
	\draw[thick] (O) -- (A);
	\draw[dashed] (O) -- (B);
	\draw[dashed] (O) -- (C);
	\draw[thick] (O) -- (D);
	\fill (O) [fill=black, draw] circle (2pt);
	\foreach \x in {A, B, C, D} {
		\fill (\x) [fill=white, draw] circle (2pt);
	}
\end{tikzpicture}
\subcaption{fixed vertex of $\operatorname{H}_{k}$}
\end{minipage}
&
\begin{minipage}{0.3\hsize}
\centering
\begin{tikzpicture}
	\coordinate (O) at (0, 0);
	\coordinate (A) at (0, 1);
	\coordinate (B) at (1, 0);
	\coordinate (C) at (0, -1); 
	\coordinate (D) at (-1, 0);
	\fill[fill=lightgray, opacity=0.5] (0, 0) -- (0, 1) -- (1, 1) -- (1, 0) --cycle;
	\fill[fill=lightgray, opacity=0.5] (0, 0) -- (0, -1) -- (-1, -1) -- (-1, 0) -- cycle;
	\draw[thick] (O) -- (A);
	\draw[dashed] (O) -- (B);
	\draw[thick] (O) -- (C);
	\draw[dashed] (O) -- (D);
	\fill (O) [fill=black, draw] circle (2pt);
	\foreach \x in {A, B, C, D} {
		\fill (\x) [fill=white, draw] circle (2pt);
	}
\end{tikzpicture}
\subcaption{fixed vertex of $\operatorname{H}_{k}$}
\end{minipage}
&
\begin{minipage}{0.3\hsize}
\centering
\begin{tikzpicture}
	\coordinate (O) at (0, 0);
	\coordinate (A) at (0, 1);
	\coordinate (B) at (1, 0);
	\coordinate (C) at (0, -1); 
	\coordinate (D) at (-1, 0);
	\fill[fill=lightgray, opacity=0.5] (0, 0) -- (0, 1) -- (1, 1) -- (1, 0) --cycle;
	\fill[fill=lightgray, opacity=0.5] (0, 0) -- (0, -1) -- (-1, -1) -- (-1, 0) -- cycle;
	\draw[thick] (O) -- (A);
	\draw[thick] (O) -- (B);
	\draw[dashed] (O) -- (C);
	\draw[dashed] (O) -- (D);
	\fill (O) [fill=black, draw] circle (2pt);
	\foreach \x in {A, B, C, D} {
		\fill (\x) [fill=white, draw] circle (2pt);
	}
\end{tikzpicture}
\subcaption{not fixed vertex of $\operatorname{H}_{k}$}
\end{minipage}
\end{tabular}
\caption{
	If $i + j$ and $k$ do not have the same parity
}
\label{fig:fixed_H1_even}
\end{figure}

In a FPL, 
the state of the interior plaquette ${\alpha}_{i, j}$ can be in 16 different situations. 
Figure \ref{fig:fixed_vertex_on_N_nonzero} and \ref{fig:fixed_vertex_on_N_zero}
illustrates the 16 different situations and their fixed vertex of $\operatorname{H}_{k}$
when $i + j$ and $k$ have the same parity as $\circ$.
\begin{figure}[H]
	\centering
	\begin{tabular}{cc}
		\begin{tikzpicture}
			\coordinate (A) at (0, 0);
			\coordinate (B) at (0, 1);
			\coordinate (C) at (1, 1);
			\coordinate (D) at (1, 0);
			\foreach \x/\y in {A/B, C/D} {
				\draw[thick] (\x) -- (\y);
			}
			\foreach \x/\y in {B/C, D/A} {
				\draw[dashed] (\x) -- (\y);
			}
			\foreach \x in {A, B, C, D} {
				\fill[fill=white, draw] (\x) circle (2pt);
			}
		\end{tikzpicture}
		&
		\begin{tikzpicture}
			\coordinate (A) at (0, 0);
			\coordinate (B) at (0, 1);
			\coordinate (C) at (1, 1);
			\coordinate (D) at (1, 0);
			\foreach \x/\y in {A/B, C/D} {
				\draw[dashed] (\x) -- (\y);
			}
			\foreach \x/\y in {B/C, D/A} {
				\draw[thick] (\x) -- (\y);
			}
			\foreach \x in {A, B, C, D} {
				\fill[fill=white, draw] (\x) circle (2pt);
			}
		\end{tikzpicture}
	\end{tabular}
	\caption{$\mathcal{N}_{\alpha} \neq 0$}
	\label{fig:fixed_vertex_on_N_nonzero}
\end{figure}
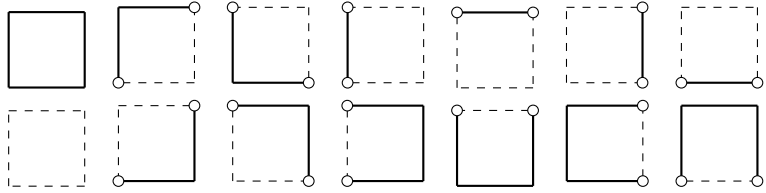
\begin{figure}[H]
	\centering
	\begin{tabular}{ccccccc}
		\begin{tikzpicture}
			\coordinate (A) at (0, 0);
			\coordinate (B) at (0, 1);
			\coordinate (C) at (1, 1);
			\coordinate (D) at (1, 0);
			\foreach \x/\y in {A/B, B/C, C/D, D/A} {
				\draw[thick] (\x) -- (\y);
			}
		\end{tikzpicture}
	&
		\begin{tikzpicture}
			\coordinate (A) at (0, 0);
			\coordinate (B) at (0, 1);
			\coordinate (C) at (1, 1);
			\coordinate (D) at (1, 0);
			\foreach \x/\y in {A/B, B/C} {
				\draw[thick] (\x) -- (\y);
			}
			\foreach \x/\y in {C/D, D/A} {
				\draw[dashed] (\x) -- (\y);
			}
			\foreach \x in {A, C} {
				\fill[fill=white, draw] (\x) circle (2pt);
			}
		\end{tikzpicture}
	&
		\begin{tikzpicture}
			\coordinate (A) at (0, 0);
			\coordinate (B) at (0, 1);
			\coordinate (C) at (1, 1);
			\coordinate (D) at (1, 0);
			\foreach \x/\y in {A/B, D/A} {
				\draw[thick] (\x) -- (\y);
			}
			\foreach \x/\y in {B/C, C/D} {
				\draw[dashed] (\x) -- (\y);
			}
			\foreach \x in {B, D} {
				\fill[fill=white, draw] (\x) circle (2pt);
			}
		\end{tikzpicture}
	&
	\begin{tikzpicture}
			\coordinate (A) at (0, 0);
			\coordinate (B) at (0, 1);
			\coordinate (C) at (1, 1);
			\coordinate (D) at (1, 0);
			\foreach \x/\y in {A/B} {
				\draw[thick] (\x) -- (\y);
			}
			\foreach \x/\y in {B/C, C/D, D/A} {
				\draw[dashed] (\x) -- (\y);
			}
			\foreach \x in {A, B} {
				\fill[fill=white, draw] (\x) circle (2pt);
			}
		\end{tikzpicture}
	&
		\begin{tikzpicture}
			\coordinate (A) at (0, 0);
			\coordinate (B) at (0, 1);
			\coordinate (C) at (1, 1);
			\coordinate (D) at (1, 0);
			\foreach \x/\y in {B/C} {
				\draw[thick] (\x) -- (\y);
			}
			\foreach \x/\y in {A/B, C/D, D/A} {
				\draw[dashed] (\x) -- (\y);
			}
			\foreach \x in {B, C} {
				\fill[fill=white, draw] (\x) circle (2pt);
			}
		\end{tikzpicture}
	&
		\begin{tikzpicture}
			\coordinate (A) at (0, 0);
			\coordinate (B) at (0, 1);
			\coordinate (C) at (1, 1);
			\coordinate (D) at (1, 0);
			\foreach \x/\y in {C/D} {
				\draw[thick] (\x) -- (\y);
			}
			\foreach \x/\y in {A/B, B/C, D/A} {
				\draw[dashed] (\x) -- (\y);
			}
			\foreach \x in {C, D} {
				\fill[fill=white, draw] (\x) circle (2pt);
			}
		\end{tikzpicture}
	&
		\begin{tikzpicture}
			\coordinate (A) at (0, 0);
			\coordinate (B) at (0, 1);
			\coordinate (C) at (1, 1);
			\coordinate (D) at (1, 0);
			\foreach \x/\y in {D/A} {
				\draw[thick] (\x) -- (\y);
			}
			\foreach \x/\y in {A/B, B/C, C/D} {
				\draw[dashed] (\x) -- (\y);
			}
			\foreach \x in {A, D} {
				\fill[fill=white, draw] (\x) circle (2pt);
			}
		\end{tikzpicture}
	\\
	\begin{tikzpicture}
			\coordinate (A) at (0, 0);
			\coordinate (B) at (0, 1);
			\coordinate (C) at (1, 1);
			\coordinate (D) at (1, 0);
			\foreach \x/\y in {A/B, B/C, C/D, D/A} {
				\draw[dashed] (\x) -- (\y);
			}
		\end{tikzpicture}
	&
		\begin{tikzpicture}
			\coordinate (A) at (0, 0);
			\coordinate (B) at (0, 1);
			\coordinate (C) at (1, 1);
			\coordinate (D) at (1, 0);
			\foreach \x/\y in {A/B, B/C} {
				\draw[dashed] (\x) -- (\y);
			}
			\foreach \x/\y in {C/D, D/A} {
				\draw[thick] (\x) -- (\y);
			}
			\foreach \x in {A, C} {
				\fill[fill=white, draw] (\x) circle (2pt);
			}
		\end{tikzpicture}
	&
		\begin{tikzpicture}
			\coordinate (A) at (0, 0);
			\coordinate (B) at (0, 1);
			\coordinate (C) at (1, 1);
			\coordinate (D) at (1, 0);
			\foreach \x/\y in {A/B, D/A} {
				\draw[dashed] (\x) -- (\y);
			}
			\foreach \x/\y in {B/C, C/D} {
				\draw[thick] (\x) -- (\y);
			}
			\foreach \x in {B, D} {
				\fill[fill=white, draw] (\x) circle (2pt);
			}
		\end{tikzpicture}
	&
	\begin{tikzpicture}
			\coordinate (A) at (0, 0);
			\coordinate (B) at (0, 1);
			\coordinate (C) at (1, 1);
			\coordinate (D) at (1, 0);
			\foreach \x/\y in {A/B} {
				\draw[dashed] (\x) -- (\y);
			}
			\foreach \x/\y in {B/C, C/D, D/A} {
				\draw[thick] (\x) -- (\y);
			}
			\foreach \x in {A, B} {
				\fill[fill=white, draw] (\x) circle (2pt);
			}
		\end{tikzpicture}
	&
		\begin{tikzpicture}
			\coordinate (A) at (0, 0);
			\coordinate (B) at (0, 1);
			\coordinate (C) at (1, 1);
			\coordinate (D) at (1, 0);
			\foreach \x/\y in {B/C} {
				\draw[dashed] (\x) -- (\y);
			}
			\foreach \x/\y in {A/B, C/D, D/A} {
				\draw[thick] (\x) -- (\y);
			}
			\foreach \x in {B, C} {
				\fill[fill=white, draw] (\x) circle (2pt);
			}
		\end{tikzpicture}
	&
		\begin{tikzpicture}
			\coordinate (A) at (0, 0);
			\coordinate (B) at (0, 1);
			\coordinate (C) at (1, 1);
			\coordinate (D) at (1, 0);
			\foreach \x/\y in {C/D} {
				\draw[dashed] (\x) -- (\y);
			}
			\foreach \x/\y in {A/B, B/C, D/A} {
				\draw[thick] (\x) -- (\y);
			}
			\foreach \x in {C, D} {
				\fill[fill=white, draw] (\x) circle (2pt);
			}
		\end{tikzpicture}
	&
		\begin{tikzpicture}
			\coordinate (A) at (0, 0);
			\coordinate (B) at (0, 1);
			\coordinate (C) at (1, 1);
			\coordinate (D) at (1, 0);
			\foreach \x/\y in {D/A} {
				\draw[dashed] (\x) -- (\y);
			}
			\foreach \x/\y in {A/B, B/C, C/D} {
				\draw[thick] (\x) -- (\y);
			}
			\foreach \x in {A, D} {
				\fill[fill=white, draw] (\x) circle (2pt);
			}
		\end{tikzpicture}
	\end{tabular}
	\caption{$\mathcal{N}_{\alpha} = 0$}
	\label{fig:fixed_vertex_on_N_zero}
\end{figure}
Let ${\alpha}_{i, j}$ be an odd interior plaquette, and $\psi \in \fpl(n, {\tau}_{-})$.
If ${\alpha}_{i, j}$ has a fixed vertex $v$ of $\operatorname{H}_{1}$ in $\psi$,
there uniquely exists a vertex $w \in V(\alpha)$
such that $w$ is a fixed vertex of $\operatorname{H}_{1}$ in $\psi$,
and it is connected to $v$ by 
black (resp. white) path
which does not pass through other fixed veretices of $\operatorname{H}_{1}$ in $\psi$.
Even after we operate 
$\operatorname{H}_{1}$,
$v$ and $w$ are fixed vertices of $\operatorname{H}_{1}$ in $\operatorname{H}_{1}\psi$,
and they are connected by black (resp. white) path
which does not pass through other fixed veretices of $\operatorname{H}_{1}$ in $\operatorname{H}_{1}\psi$.

Now, let $\boldsymbol{p} = ({v}_{0}, {v}_{1}, \ldots , {v}_{m})$ be a black (resp. white) path in $\psi$,
and ${v}_{0}$, ${v}_{m}$ fixed vertices of $\operatorname{H}_{1}$ in $\psi$.
We divide $\boldsymbol{p}$ at each fixed vertices of $\operatorname{H}_{1}$ in $\psi$
into short black (resp. white) pathes.
When $\boldsymbol{p}$ is divided into $k$ pathes, 
we set ${v}_{0}$ as ${w}_{0}$, ${v}_{m}$ as ${w}_{k}$, and
the $k$ short pathes as
$({w}_{0}, \ldots , {w}_{1})$, $({w}_{1}, \ldots , {w}_{2})$, $\ldots$, $({w}_{k - 1}, \ldots , {w}_{k})$.
Since ${w}_{i}$ and ${w}_{i + 1}$ are connected by balck (resp. white) path in $\operatorname{H}_{1}\psi$
for $0 \leq i < k$,
there is a black (resp. white) path which start from ${v}_{0}$ to ${v}_{m}$
and pass through ${w}_{1}, {w}_{2}, \ldots , {w}_{k - 1}$
in $\operatorname{H}_{1}\psi$. 
Moreover
such a black (resp. white) path
does not pass through any other fixed vertices of $\operatorname{H}_{1}$.

Next we focus on boundary plaquette. 
Let ${\alpha}_{i, j}$
be an odd boundary plaquette, and $\psi \in \fpl(n, {\tau}_{-})$.
In the plaquette ${\alpha}_{i, j}$, 
there is exactly one fixed vertex of $\operatorname{H}_{1}$ in $\psi$.
Now 
we label the fixed vertex as $w$,
and the two boundary edges in ${\alpha}_{i, j}$ as ${e}_{2k}$ and ${e}_{2k + 1}$.
Here, ${e}_{4n + 1}$ means ${e}_{1}$.
There is a black (resp. white) path
which connects $w$ and the boundary vertex
that is adjacent to ${e}_{2k}$ (resp. ${e}_{2k + 1}$),
and this 
path does not pass through 
a fixed vertex of $\operatorname{H}_{1}$ in $\psi$ other than $w$. 
Even after we operate $\operatorname{H}_{1}$, 
$w$ is a fixed vertex of $\operatorname{H}_{1}$ in $\operatorname{H}_{1}\psi$, 
and the other vertices in ${\alpha}_{i, j}$ are not 
fixed vertices. 
Moreover
there is a black (resp. white) path
which connects $w$ and the boundary vertex
that is adjacent to ${e}_{2k + 1}$ (resp. ${e}_{2k}$),
and this 
path does not pass through 
a fixed vertex of $\operatorname{H}_{1}$ in $\operatorname{H}_{1}\psi$ other than $w$.
Figure \ref{fig:fixed_boundary_H1}
illustrates the state of ${\alpha}_{i, j}$ in $\psi \in \fpl(n, {\tau}_{-})$
before and after we operate $\operatorname{H}_{1}$.
\input{tikz_fixed_boundary_H1_ver3.tex}
\par
Let ${v}_{i}$ be the boundary vertex 
which is adjacent to boundary edge ${e}_{i}$
for $1 \leq i \leq 2n$, 
and $\psi \in \fpl(n, {\tau}_{-})$.
We set $\boldsymbol{p}$ as a black (resp. white) path in $\psi$
which start from ${v}_{2k}$ (resp. ${v}_{2k + 1}$)
and end at ${v}_{2l}$ (resp. ${v}_{2l + 1}$).
We also set $w$ and $w'$ as the vertices in $\boldsymbol{p}$ 
which satisfy following condition:
when we divide $\boldsymbol{p}$ at each fixed vertices 
of $\operatorname{H}_{1}$ in $\psi$
into short black (resp. white) pathes, 
the short 
path which strat from ${v}_{2k}$ (resp. ${v}_{2k + 1}$) 
is ended at $w$, 
and the short 
path which end at ${v}_{2l}$ (resp. ${v}_{2l + 1}$)
is started from $w'$. 
From the above,
we have a black (resp. white) path $\boldsymbol{q}$ in $\operatorname{H}_{1}\psi$
which start from $w$ and end at $w'$.
Now, we glue the following three black (resp. white) pathes
in $\operatorname{H}_{1}\psi$:
the 
path which strat from ${v}_{2k + 1}$ (resp. ${v}_{2k}$) and end at $w$, 
$\boldsymbol{q}$, and
the 
path which strat from $w'$ and end at ${v}_{2l + 1}$ (resp. ${v}_{2l}$).
Then
we get the black (resp. white) path in $\operatorname{H}_{1}\psi$
which start from ${v}_{2k + 1}$ (resp. ${v}_{2k}$), 
pass through $w$ and $w'$, 
and end at ${v}_{2l + 1}$ (resp. ${v}_{2l}$).
\begin{figure}[htbp]
	\centering
	\begin{tikzpicture}
		\foreach \x in {1, 2, 3, 4}{
			\foreach \y in {1, 2, 3, 4}{
				\coordinate (v_\x_\y) at ($(\x, \y) + (1, 0)$);
			}
		}
		\foreach \x/\k in {1/1, 2/2, 3/3, 4/4}{
			\coordinate (b_\k) at ($(\x, 0) + (1, 0)$);
			\coordinate (v_\x_0) at (b_\k);
		}
		\foreach \y/\k in {1/5, 2/6, 3/7, 4/8}{
			\coordinate (b_\k) at ($(5, \y) + (1, 0)$);
			\coordinate (v_5_\y) at (b_\k);
		}
		\foreach \x/\k in {4/9, 3/10, 2/11, 1/12}{
			\coordinate (b_\k) at ($(\x, 5) + (1, 0)$);
			\coordinate (v_\x_5) at (b_\k);
		}
		\foreach \y/\k in {4/13, 3/14, 2/15, 1/16}{
			\coordinate (b_\k) at ($(0, \y) + (1, 0)$);
			\coordinate (v_0_\y) at (b_\k);
		}
		\foreach \x/\xx in {1/2, 3/4}{
			\foreach \y/\yy in {0/1, 2/3, 4/5}{
				\fill[fill = lightgray, opacity=0.5]
					(v_\x_\y) -- (v_\xx_\y) -- (v_\xx_\yy) -- (v_\x_\yy) -- cycle;
			}
		}
		\foreach \x/\xx in {0/1, 2/3, 4/5}{
			\foreach \y/\yy in {1/2, 3/4}{
				\fill[fill = lightgray, opacity=0.5]
					(v_\x_\y) -- (v_\xx_\y) -- (v_\xx_\yy) -- (v_\x_\yy) -- cycle;
			}
		}
			\draw[very thick] (v_3_4) -- (b_10);
			\draw[dashed] (v_4_4) -- (v_3_4);
			\draw[dashed] (v_4_4) -- (b_9);
			\coordinate (w_0) at (v_3_4);
			\draw[very thick] (v_3_4) -- (v_2_4);
			\draw[very thick] (v_2_3) -- (v_2_4);
			\draw[very thick] (v_2_3) -- (v_3_3);
			\draw[dashed] (v_3_3) -- (v_3_4);
			\coordinate (w_1) at (v_3_3);
			\draw[very thick] (v_3_2) -- (v_3_3);
			\draw[dashed] (v_3_3) -- (v_4_3);
			\draw[thick] (v_4_3) -- (v_4_2);
			\draw[dashed] (v_4_2) -- (v_3_2);
			\coordinate (w_2) at (v_3_2);
			\draw[very thick] (v_3_2) -- (v_2_2);
			\draw[dashed] (v_2_2) -- (v_2_1);
			\draw[thick] (v_2_1) -- (v_3_1);
			\draw[dashed] (v_3_1) -- (v_3_2);
			\coordinate (w_3) at (v_2_2);
			\draw[very thick] (v_1_2) -- (v_2_2);
			\draw[dashed] (v_2_2) -- (v_2_3);
			\draw[dashed] (v_1_3) -- (v_1_2);
			\draw[dashed] (v_1_3) -- (v_2_3);
			\coordinate (w_4) at (v_1_2);
			\draw[very thick] (v_1_2) -- (v_1_1);
			\draw[very thick] (b_16) -- (v_1_1);
			\draw[dashed] (b_15) -- (v_1_2);
		\foreach \k in {0, 1, 2, 3, 4}{
			\fill[fill=white, draw] (w_\k) circle (3pt);
			\node at (w_\k)[below left]{${w}_{\k}$};
		}
		\node at (b_10)[above]{${v}_{15}$};
		\node at (b_16)[left]{${v}_{5}$};
		\foreach \x in {1, 2, 3, 4}{
			\foreach \y in {1, 2, 3, 4}{
				\coordinate (u_\x_\y) at ($(\x, \y) - (6, 0)$);
			}
		}
		\foreach \x/\k in {1/1, 2/2, 3/3, 4/4}{
			\coordinate (d_\k) at ($(\x, 0) - (6, 0)$);
			\coordinate (u_\x_0) at (d_\k);
		}
		\foreach \y/\k in {1/5, 2/6, 3/7, 4/8}{
			\coordinate (d_\k) at ($(5, \y) - (6, 0)$);
			\coordinate (u_5_\y) at (d_\k);
		}
		\foreach \x/\k in {4/9, 3/10, 2/11, 1/12}{
			\coordinate (d_\k) at ($(\x, 5) - (6, 0)$);
			\coordinate (u_\x_5) at (d_\k);
		}
		\foreach \y/\k in {4/13, 3/14, 2/15, 1/16}{
			\coordinate (d_\k) at ($(0, \y) - (6, 0)$);
			\coordinate (u_0_\y) at (d_\k);
		}
		\foreach \x/\xx in {1/2, 3/4}{
			\foreach \y/\yy in {0/1, 2/3, 4/5}{
				\fill[fill = lightgray, opacity=0.5]
					(u_\x_\y) -- (u_\xx_\y) -- (u_\xx_\yy) -- (u_\x_\yy) -- cycle;
			}
		}
		\foreach \x/\xx in {0/1, 2/3, 4/5}{
			\foreach \y/\yy in {1/2, 3/4}{
				\fill[fill = lightgray, opacity=0.5]
					(u_\x_\y) -- (u_\xx_\y) -- (u_\xx_\yy) -- (u_\x_\yy) -- cycle;
			}
		}
			\draw[dashed] (u_3_4) -- (d_10);
			\draw[very thick] (u_4_4) -- (u_3_4);
			\draw[very thick] (u_4_4) -- (d_9);
			\coordinate (m_0) at (u_3_4);
			\draw[dashed] (u_3_4) -- (u_2_4);
			\draw[dashed] (u_2_3) -- (u_2_4);
			\draw[dashed] (u_2_3) -- (u_3_3);
			\draw[very thick] (u_3_3) -- (u_3_4);
			\coordinate (m_1) at (u_3_3);
			\draw[very thick] (u_3_2) -- (u_3_3);
			\draw[dashed] (u_3_3) -- (u_4_3);
			\draw[thick] (u_4_3) -- (u_4_2);
			\draw[dashed] (u_4_2) -- (u_3_2);
			\coordinate (m_2) at (u_3_2);
			\draw[very thick] (u_3_2) -- (u_2_2);
			\draw[dashed] (u_2_2) -- (u_2_1);
			\draw[thick] (u_2_1) -- (u_3_1);
			\draw[dashed] (u_3_1) -- (u_3_2);
			\coordinate (m_3) at (u_2_2);
			\draw[dashed] (u_1_2) -- (u_2_2);
			\draw[very thick] (u_2_2) -- (u_2_3);
			\draw[very thick] (u_1_3) -- (u_1_2);
			\draw[very thick] (u_1_3) -- (u_2_3);
			\coordinate (m_4) at (u_1_2);
			\draw[dashed] (u_1_2) -- (u_1_1);
			\draw[dashed] (d_16) -- (u_1_1);
			\draw[very thick] (d_15) -- (u_1_2);
		\foreach \k in {0, 1, 2, 3, 4}{
			\fill[fill=white, draw] (m_\k) circle (3pt);
			\node at (m_\k)[below left]{${w}_{\k}$};
		}
		\node at (d_9)[above]{${v}_{14}$};
		\node at (d_15)[left]{${v}_{4}$};
	\end{tikzpicture}
	\caption{
		left:the black path in $\psi \in \mathfrak{fpl}(4, {\tau}_{-})$,
		right:the black path in $\operatorname{H}_{1}\psi \in \mathfrak{fpl}(4, {\tau}_{+})$.
	}
\end{figure}
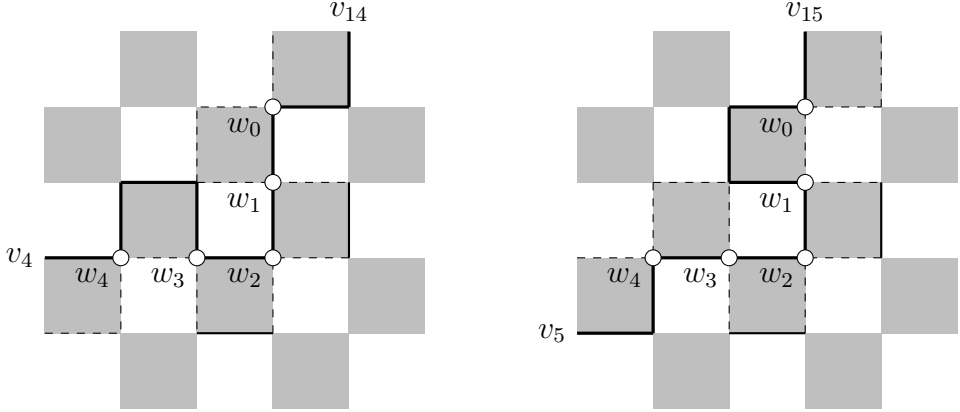
Therefore, 
$
	\{ k + 1, l + 1 \} \in {\pi}_{b, +}(\operatorname{H}_{1}\psi)
$
(resp.
$
	\{ k, l \} \in {\pi}_{w, +}(\operatorname{H}_{1}\psi)
$)
if
$
	\{ k, l \} \in {\pi}_{b, -}(\psi)
$
(resp.
$
	\{ k + 1, l + 1 \} \in {\pi}_{w, -}(\psi)
$).
Then we have two following equations:
\begin{subequations}
\begin{align}
	\label{eq:pi_b_H1}
	{\pi}_{b, +}(\operatorname{H}_{1}\psi)
		&=	\operatorname{R}^{-1} {\pi}_{b, -}(\psi),	\\
	\label{eq:pi_w_H1}
	{\pi}_{w, +}(\operatorname{H}_{1}\psi)
		&=	\operatorname{R} {\pi}_{w, -}(\psi).
\end{align}
\end{subequations}

Next, we focus on monochromatic cycles.
Let $\psi$ be a FPL in $\fpl(n, {\tau}_{-})$, and
$\boldsymbol{c}$ a black (resp. white) cycle in $\psi$. 
Note that $\boldsymbol{c}$ has length at least $4$. 
First we consider the case 
when the length of $\boldsymbol{c}$ is $4$. 
Now
there is an interior plaquette ${\alpha}_{i, j}$
such that $\boldsymbol{c}$ is consist of $4$ edges in ${\alpha}_{i, j}$
(i.e., 
$
	\boldsymbol{c}
		=	\left(
			(i, j), (i + 1, j), (i + 1, j + 1), (i, j + 1), (i, j)
		\right)
$).
If the plaquette ${\alpha}_{i, j}$ is even, 
the $4$ vertices are fixed vertices of $\operatorname{H}_{1}$ in $\psi$, 
and thus there is a black (resp. white) cycle
$\boldsymbol{c}'$ in $\operatorname{H}_{1}\psi$
such that $\boldsymbol{c}'$ pass through 
the $4$ verices of ${\alpha}_{i, j}$.
If the plaquette ${\alpha}_{i, j}$ is odd, 
since the color of the all edges in ${\alpha}_{i, j}$ is reversed, 
there is a white (resp. black) cycle 
$\boldsymbol{c}'$in $\operatorname{H}_{1}\psi$
such that $\boldsymbol{c}'$ is consist of $4$ edges in ${\alpha}_{i, j}$. 
\par
Second we consider the case
when the length of $\boldsymbol{c}$ is greater than $4$. 
Then
the cycle $\boldsymbol{c}$ run across multiple plaquettes, 
and pass through some fixed vertices of $\operatorname{H}_{1}$ in $\psi$. 
Therefore 
there is a black (resp. white) cycle
$\boldsymbol{c}'$ in $\operatorname{H}_{1}\psi$
such that $\boldsymbol{c}$ and $\boldsymbol{c}'$ 
pass through the same fixed vertices of $\operatorname{H}_{1}$.

From the above, 
we can construct the bijection between
the set of
monochromatic cycles in $\psi$
and the set of
monochromatic cycles in $\operatorname{H}_{1}\psi$. 
Together with \eqref{eq:pi_b_H1} and \eqref{eq:pi_w_H1}, 
we showed \eqref{eq:H1_lem}.
\par
In the same way, 
each black (resp. white) path which pass through 
some fixed vertices of $\operatorname{H}_{0}$ in $\operatorname{H}_{1}\psi$
also pass through same fixed vertices of $\operatorname{H}_{0}$ 
after we operate $\operatorname{H}_{0}$.
Moreover,
each even boundary plaquette has exactly one fixed vertex of $\operatorname{H}_{0}$.
If we label the two boundary edge in the plaquette as ${e}_{2k}$ and ${e}_{2k - 1}$, 
then the fixed vertex of $\operatorname{H}_{0}$
is connected to the boundary vertex ${v}_{2k - 1}$ (resp. ${v}_{2k}$)
by black (resp. white) path in $\operatorname{H}_{1}\psi$. 
After we operate $\operatorname{H}_{0}$,
the fixed vertex of $\operatorname{H}_{0}$
is connected to the boundary vertex ${v}_{2k}$ (resp. ${v}_{2k - 1}$)
by black (resp. white) path in $\operatorname{H}_{0}\operatorname{H}_{1}\psi$. 
Hence
we have two following equations:
\begin{subequations}
	\begin{align}
		\label{eq:pi_b_H0}
		{\pi}_{b, -}\left( \operatorname{H}_{0}\operatorname{H}_{1}\psi \right)
		&=
		{\pi}_{b, +}\left( \operatorname{H}_{1}\psi \right)
		\\
		\label{eq:pi_w_H0}
		{\pi}_{w, -}\left( \operatorname{H}_{0}\operatorname{H}_{1}\psi \right)
		&=
		{\pi}_{w, +}\left( \operatorname{H}_{1}\psi \right)
	\end{align}
\end{subequations}
On the other hand,
$\operatorname{H}_{0}\operatorname{H}_{1}\psi$
have the same number of monochromatic cycles as
$\operatorname{H}_{1}\psi$.
Therefore we hold \eqref{eq:H0_lem},
and lemma \ref{lem:H0H1} is prooved.
Then proposition \ref{prop:Wieland} follows from lemma \ref{lem:H0H1}.
\qed
\subsection{Periodicity of the gyration}
We remark that
the operation of repeating the gyration is periodic.
Let $\psi$ be a FPL in $\fpl(n, {\tau}_{-})$.
If we repeat the gyration to $\psi$ $n$ times,
we will get $\psi$ at least once.
Now
we set
$
	\fpl\left( n, {\tau}_{-}; \mathcal{O}(\psi) \right)
	:=
	\left\{
		\operatorname{G}^{k} \psi 
	\, \middle| \,
		0 \leq k < n
	\right\}
$.
Then
we hold the following propositions.
\begin{proposition}\label{prop:periodicity_N_alpha}
	Let $\psi$ be a FPL in $\fpl(n, {\tau}_{-})$ and $\alpha$ an interior plaquette.
	Then
	we have a following equation:
	\begin{equation}\label{eq:periodicity_N_alpha}
		\sum_{\varphi \in \fpl \left( n, {\tau}_{-}; \mathcal{O}(\psi) \right)}
			\mathcal{N}_{\alpha}(\varphi)
		=
		0.
	\end{equation} 
\end{proposition}
We also take several steps to prove proposition \ref{prop:periodicity_N_alpha}.
First,
we focus on the case
when $\alpha$ has an edge in common with a boundary plaquette.
Let $\alpha$ be an odd plaquette, and
$\beta$ the boundary plaquette
which has an edge $e$ in common with $\alpha$.
For each $k \geq 0$,
$\operatorname{G}^{k}\psi(e) \neq \operatorname{G}^{k + 1}\psi(e)$
if and only if
$\operatorname{G}_{1}$
reverse the color of $e$
in $\operatorname{G}^{k}\psi$
because 
$\operatorname{G}_{0}$
does not reverse the color of $e$
in  $\operatorname{G}_{1}\operatorname{G}^{k}\psi$.
We remark that
if $\psi$ satisfies the following conditions:
$
	\operatorname{G}^{k}\psi(e)
	\neq 
	\operatorname{G}^{k + 1}\psi(e)
$,
$
	\operatorname{G}^{k + 1}\psi(e)
	=
	\operatorname{G}^{k + 2}\psi(e)
	=
	\cdots
	=
	\operatorname{G}^{l}\psi(e)
$,
and
$
	\operatorname{G}^{l}\psi(e)
	\neq
	\operatorname{G}^{l + 1}\psi(e)
$,
then
$
	\mathcal{N}_{\alpha}\left( \operatorname{G}^{k}\psi \right)
	\neq
	0
$,
$
	\mathcal{N}_{\alpha}\left( \operatorname{G}^{k + 1}\psi \right)
	=
	\mathcal{N}_{\alpha}\left( \operatorname{G}^{k + 2}\psi \right)
	=
	\cdots
	=
	\mathcal{N}_{\alpha}\left( \operatorname{G}^{l - 1}\psi \right)
	=
	0
$,
and
$
	\mathcal{N}_{\alpha}\left( \operatorname{G}^{l}\psi \right)
	=
	- \mathcal{N}_{\alpha}\left( \operatorname{G}^{k}\psi \right)
$.
In other word,
the sequence
$
	{
		\left\{
			\mathcal{N}_{\alpha}\left( \operatorname{G}^{k}\psi \right)
		\right\}
	}_{
		k \geq 0
	}
$
alternates between $1$ and $-1$ except for 0.
Moreover,
when we denote $\# \fpl \left( n, {\tau}_{-}; \mathcal{O}(\psi) \right)$ as $m$,
$1$ and $-1$ appear the same number of times
between the $0$th term and the 
$m$th term in 
$
	{
		\left\{
			\mathcal{N}_{\alpha}\left( \operatorname{G}^{k}\psi \right)
		\right\}
	}_{
		k \geq 0
	}
$
because
the color of $e$ is reversed an even number times
while the gyration is repaeted 
$m$ times.
Since
the left hand side of \eqref{eq:periodicity_N_alpha}
is equal to 
the sum of the $0$th through $m$th terms of
$
	{
		\left\{
			\mathcal{N}_{\alpha}\left( \operatorname{G}^{k}\psi \right)
		\right\}
	}_{
		k \geq 0
	}
$,
then
\eqref{eq:periodicity_N_alpha} 
follows.
When $\alpha$ is even,
we can show the same by replacing
$\operatorname{G}_{1}$ with $\operatorname{G}_{0}$
and $\operatorname{G}$ with $\operatorname{G}^{-1}$.

Next, 
we consider two distinct interior plaquette 
they have a common edge $e$.
Let $\alpha$ (resp. $\beta$) be an odd (resp. even) interior plaquette which has the edge $e$.
We set two sequences
$
	{
		\left\{
			{\mu}_{n}
		\right\}
	}_{
		n \geq 0
	}
$
and
$
	{
		\left\{
			{\nu}_{n}
		\right\}
	}_{
		n \geq 0
	}
$
as following:
\begin{align}
	{\mu}_{n}
	&:=
	\begin{cases}
		\mathcal{N}_{\alpha} \left( \operatorname{G}^{k}\psi \right)	
		&	\left( n = 2k \right)	\\
		\mathcal{N}_{\beta} \left( \operatorname{G}_{1}\operatorname{G}^{k}\psi \right)	
		&	\left( n = 2k + 1 \right)
	\end{cases}
	, \\
	{\nu}_{n}
	&:=
	\begin{cases}
		\operatorname{G}^{k}\psi(e)	
		&	\left( n = 2k \right)	\\
		\operatorname{G}_{1}\operatorname{G}^{k}\psi(e)	
		&	\left( n = 2k + 1 \right)
	\end{cases}
	.
\end{align}
Note that
${\nu}_{n} \neq {\nu}_{n + 1}$
if and only if
${\mu}_{n} \neq 0$.
Similar to previous discussion,
$
	{
		\left\{
			{\mu}_{n}
		\right\}
	}_{
		n \geq 0
	}
$
alternates between $1$ and $-1$ except for 0,
and the sum of the $0$th through $(2m -1)$th terms of 
$
	{
		\left\{
			{\mu}_{n}
		\right\}
	}_{
		n \geq 0
	}
$
is equal to $0$.
Since
$
	\operatorname{G}_{1}\operatorname{G}^{k}\psi
$
is equal to
$
	\operatorname{G}^{-k}\operatorname{G}_{1}\psi
$,
we hold the following equation:
\begin{equation}
	\label{eq:alpha_and_beta}
	\sum_{0 \leq k \leq m - 1}
		\mathcal{N}_{\alpha}\left( \operatorname{G}^{k}\psi \right)
	+
	\sum_{0 \leq k \leq m - 1}
		\mathcal{N}_{\beta}\left( \operatorname{G}^{-k}\operatorname{G}_{1}\psi \right)
	=
	0.
\end{equation}
Since
the cardinarity of
$
	\fpl\left(
		n, {\tau}_{-}; \mathcal{O}\left( \operatorname{G}_{1}\psi \right)
	\right)
$
is equal to $m$,
the equation
\eqref{eq:alpha_and_beta} is equivalent to the following equation:
\begin{equation}
	\label{eq:alpha_and_beta_2}
	\sum_{
		\varphi \in
		\fpl\left(
			n, {\tau}_{-}; \mathcal{O}\left( \psi \right)
		\right)
	}
	\mathcal{N}_{\alpha}\left( \varphi \right)
	+
	\sum_{
		{\varphi}' \in
		\fpl\left(
			n, {\tau}_{-}; \mathcal{O}\left( \operatorname{G}_{1}\psi \right)
		\right)
	}
	\mathcal{N}_{\beta}\left( {\varphi}' \right)
	=
	0.
\end{equation}
Combined with the previous discussion,
proposition \ref{prop:periodicity_N_alpha} follows.
\qed
\section{The vector space which has link patterns as basis}
We consider  the $\mathbb{C}$-vector space which has link pattern  as basis
to refine a enumertation of FPLs. 
Let $n$ be a positive integer.
We denote ${\mathbb{C}}^{\mathcal{F}(2n)}$ 
as the $\mathbb{C}$-vector space which is spaned by $\mathcal{F}(2n)$,
and
we will write each element $x \in  {\mathbb{C}}^{\mathcal{F}(2n)}$ 
as $\left| x \right\rangle$ for emphasis. 
We also denote ${\mathbb{C}}^{\fpl(n, \tau)}$ 
as the $\mathbb{C}$-vector space which is spaned by $\fpl(n, \tau)$
when a boundary condition $\tau \in { \{ b, w \} }^{4n}$ is given,   
and we will write each element $y \in  {\mathbb{C}}^{\fpl(n, \tau)}$ 
as $\left\| y \right\rrangle$ for emphasis. 
In particular, 
we denote   
$
		\sum_{\psi \in \fpl(n, {\tau}_{-})}
			\left| {\pi}_{b, -}(\psi) \right\rangle 
	\in {\mathbb{C}}^{\mathcal{F}(2n)}
$
as
$\left| {s}_{n} \right\rangle$,
and
$
		\sum_{\psi \in \fpl(n, {\tau}_{-})}
			\left\| \psi \right\rrangle 
	\in {\mathbb{C}}^{\fpl(n, {\tau}_{-})}
$
as
$
	\left\| {s}_{n, {\tau}_{-}} \right\rrangle
$
.
\subsection{Operators on the vector space}
We define some operators on the vector space which are introduced above.
First,
we define 
${\Pi}_{-} \colon {\mathbb{C}}^{\fpl(n, {\tau}_{-})} \rightarrow {\mathbb{C}}^{\mathcal{F}(2n)}$ 
(Resp. ${\Pi}_{+} \colon {\mathbb{C}}^{\fpl(n, {\tau}_{+})} \rightarrow {\mathbb{C}}^{\mathcal{F}(2n)}$) 
as follows: 
\begin{subequations}
	\begin{align}
		{\Pi}_{-}\left( 
			\displaystyle\sum_{\psi \in \fpl(n, {\tau}_{-})} 
				{c}_{\psi} \left\| \psi \right\rrangle
		\right) &:=
			\displaystyle\sum_{\psi \in \fpl(n, {\tau}_{-})} 
				{c}_{\psi} \left| {\pi}_{b, -}(\psi) \right\rangle
		, \\
		{\Pi}_{+}\left( 
			\displaystyle\sum_{\psi \in \fpl(n, {\tau}_{+})} 
				{c}_{\psi} \left\| \psi \right\rrangle
		\right) &:=
			\displaystyle\sum_{\psi \in \fpl(n, {\tau}_{+})} 
				{c}_{\psi} \left| {\pi}_{b, +}(\psi) \right\rangle
		.
	\end{align}
\end{subequations}
Especially, we hold 
$
	{\Pi}_{-}(\left\| {s}_{n, {\tau}_{-}} \right\rrangle) = \left| {s}_{n} \right\rangle
$. 
Second, 
we define 
$
	\tilde{\mathcal{N}}_{\alpha}
	\colon
	{\mathbb{C}}^{\fpl(n, \tau)} \rightarrow {\mathbb{C}}^{\fpl(n, \tau)}
$
for each interior plaquette $\alpha$
as following:
\begin{equation}
	\tilde{\mathcal{N}}_{\alpha}\left(
		\displaystyle\sum_{\psi \in \fpl(n, \tau)} 
			{b}_{\psi} \left\| \psi \right\rrangle
	\right) :=
		\displaystyle\sum_{\psi \in \fpl(n, \tau)}
			\mathcal{N}_{\alpha} {b}_{\psi} \left\| \psi \right\rrangle.
\end{equation}
Next,
we make operators on $\mathcal{F}(2n)$ 
extend linearly. 
When the operator $X \colon \mathcal{F}(2n) \rightarrow \mathcal{F}(2n)$ is given, 
we define $\hat{X} \colon {\mathbb{C}}^{\mathcal{F}(2n)} \rightarrow {\mathbb{C}}^{\mathcal{F}(2n)}$ 
as follows:
\begin{equation}
	\hat{X}\left(
		\displaystyle\sum_{\mu \in \mathcal{F}(2n)} 
			{a}_{\mu} \left| \mu \right\rangle
	\right) :=
		\displaystyle\sum_{\mu \in \mathcal{F}(2n)}
			{a}_{\mu} \left| X(\mu) \right\rangle.
\end{equation}
Similarly, we make operators on $\fpl(n, \tau)$ extend linearly. 
When a boundary condition $\tau \in { \{ b, w \} }^{4n}$ and 
the operator $Y \colon \fpl(n, \tau) \rightarrow \fpl(n, \tau)$ is given, 
we define $\hat{Y} \colon {\mathbb{C}}^{\fpl(n, \tau)} \rightarrow {\mathbb{C}}^{\fpl(n, \tau)}$ 
as follows:
\begin{equation}
	\hat{Y}\left(
		\displaystyle\sum_{\psi \in \fpl(n, \tau)} 
			{b}_{\psi} \left\| \psi \right\rrangle
	\right) :=
		\displaystyle\sum_{\psi \in \fpl(n, \tau)}
			{b}_{\psi} \left\| Y(\psi) \right\rrangle.
\end{equation}
\def\Sym{\operatorname{Sym}}
We define operators 
$
	\Sym \colon {\mathbb{C}}^{\mathcal{F}(2n)} \rightarrow {\mathbb{C}}^{\mathcal{F}(2n)}
$
and  
$
	\mathrm{H}_{n} \colon {\mathbb{C}}^{\mathcal{F}(2n)} \rightarrow {\mathbb{C}}^{\mathcal{F}(2n)}
$
respectively as  
$
	\mathrm{Sym} := \sum_{k = 0}^{2n - 1} 
		{\left( \hat{\operatorname{R}} \right) }^{k}
$
,
$
	\mathrm{H}_{n} := \sum_{k = 1}^{2n}
		\hat{{e}_{k}}
$
.
Then we call $\Sym$ \textsl{symmetrise operator}, 
and we call $\mathrm{H}_{n}$ \textsl{Hamiltonian}. 
Now,
the following claim hold as a corollary of proposition \ref{prop:periodicity_N_alpha}.
\begin{proposition}
	Let $\alpha$ be an interior plaquette.
	We have the following equation:
	\begin{equation}
		\label{eq:Sym_pi_N_s}
		\Sym {\Pi}_{-}\hat{ \mathcal{N}_{\alpha} }
			\left\| {s}_{n, {\tau}_{-}} \right\rrangle
		=
		0.
	\end{equation}
\end{proposition}
\begin{proof}
First,
we define a equivalence on $\mathcal{F}(2n)$.
Let $\mu$ and $\nu$ be link patterns of size $n$.
We define $\mu$ and $\nu$ are equivalent
if there exists a non-negative integer $k$
such that $\nu = \operatorname{R}^{k} \mu$,
and we write $\mu \sim \nu$.
Then we denote
$
	\mathcal{F}(2n)/ \sim
$
as
$
	\mathcal{F}^{*}(2n) 
$
,
and we denote
the complete set of 
$
	\mathcal{F}^{*}(2n) 
$
as
$
	\left\{
		{\mu}_{1}, {\mu}_{2}, \ldots , {\mu}_{m}
	\right\}
$.

Next,
we shall transforme the left hand side of \eqref{eq:Sym_pi_N_s}.
\begin{align}
	&
	\Sym {\Pi}_{-}\hat{ \mathcal{N}_{\alpha} }
		\left\| {s}_{n, {\tau}_{-}} \right\rrangle
	\nonumber \\
	=&
	\displaystyle\sum_{k = 0}^{2n - 1}
		\displaystyle\sum_{\psi \in \fpl(n, {\tau}_{-})}
			\mathcal{N}_{\alpha}\left( \psi \right)
			\left| \operatorname{R}^{k} \pi(\psi) \right\rangle
	\nonumber \\
	=&
	\displaystyle\sum_{i = 1}^{m}
		\displaystyle\sum_{
			\substack{
				\psi \in \fpl(n, {\tau}_{-})	\\
				\pi(\psi) \sim {\mu}_{i}
			}
		}
			\mathcal{N}_{\alpha}(\psi)
	\displaystyle\sum_{k = 0}^{2n - 1}
		\left| \operatorname{R}^{k} {\mu}_{i} \right\rangle
\end{align}
We remark that 
the set of FPLs the link pattern of which equivalent to ${\mu}_{i}$
is represented as disjonit union of some orbits of FPL
for each $1 \leq i \leq m$.
Therefore
$
	\displaystyle\sum_{
			\substack{
				\psi \in \fpl(n, {\tau}_{-})	\\
				\pi(\psi) \sim {\mu}_{i}
			}
		}
			\mathcal{N}_{\alpha}(\psi)
	=
	0
$
follows to proposition \ref{prop:periodicity_N_alpha},
thus \eqref{eq:Sym_pi_N_s} follows.
\end{proof}
\section{Futurework}
Althought we deal with the ordinary type of ASM in this paper,
we are interested in half-turn ASM.
We will define a poset similar to $\mathbb{P}_{n}$
for half-turn ASM. 
Further more,
we aim to make half-turn ASM correspond to a root systems.

\vskip20pt
\noindent
{\Large\textbf{Acknowledgment}}
\par
I would like to express my deep gratitude to Prof.~Masao~Ishikawa and Yoshiki~Takayama 
for helpful suggestions and advices.


\begin{thebibliography}{99}
\bibitem{Br}
Davis M. Bressoud, {
\em Proofs and Confirmations --- The Story of the Alternating Sign Matrix
Conjecture},
Cambridge University Press (1999).
\bibitem{ASM_B_order}
Richard~A.~Brualdi and Michael~W.~Schroeder,
``Alternating sign matrices and their Bruhat order'',
{\em Discrete Math.}, {\bf 340} (2017) 1996 -- 2019.
\bibitem{RS2018}
		{Luigi~Cantini and Andrea~Sportiello}, 
		{``Proof of Razumov-Stroganov conjecture''}, 
		{\em J. Combin. Theory Ser. A}, {\bf 5}
		(2011) 1549 --1574
\bibitem{Ku}
G.~Kuperberg,
``Another proof of the alternating sign matrix conjecture'',
{\em Internat. of Math. Res. Notices}, {\bf 1996} (1996) 139 -- 150.
\bibitem{Ku2}
G.~Kuperberg,
``Classes of Alternating-Sign Matrices under One Roof'',
{\em Ann. of Math.}, {\bf 156} (2002) 835 -- 866.
\bibitem{Lascoux}
A.~Lascoux and M.~Sch\"{u}tzenberger, 
“Treillis et bases des groupes de coxeter”, 
{\em Electr. J. Comb.},
3:342—351, 1996.
\bibitem{MRR1}
W.H. Mills, David P. Robbins, and Howard Rumsey, Jr.,
``Alternating Sign Matrices and Descending Plane Partitions'',
{\em J. Combin. Theory Ser. A}, {\bf 34} (1983) 340 -- 359.
\bibitem{O}
S.~Okada,
``Enumeration of Symmetry Classes of Alternating Sign Matrices and
Characters of Classical Groups'',
{\em J. Algebraic Combin.}, {\bf 23} (2006) 43 -- 69.
\bibitem{ec1}
		{Richard~P.~Stanley}, 
		{\em Enumerative~Combinatorics~Volume~1~second~edition}, 
		{Cambridge University Press}
		(2011) 
\bibitem{SW}
J.~Striker and N.~Williams,
``Promotion and rowmotion'', {\em European J. Combin.}, {\bf 33} (2012)
1919 -- 1942.
\bibitem{St}
J~Striker,
``The toggle group, homomesy, and the
Razumov-Stroganov correspondence'',
{\em Electr. J. Comb.}, {\bf 22(2)} (2015) \#P2.57
	\bibitem{Wieland_2000}
		{B.~Wieland}, 
		{``Large Dihedral Symmetry of the set of Alternating Sign Matrices}, 
		{\em Electr. J. Comb.}, {\bf 7} (2000) R37
\bibitem{DZ}
Doron~Zeilberger, 
``Proof of the Alternating Sign Matrix Conjecture''
{\em Elec. J. Comb.}, {\bf 3} (1996) R13
\end{thebibliography}
\end{document}